\documentclass[11pt]{article}
\usepackage[a4paper,margin=27mm]{geometry}
\usepackage[T1]{fontenc}
\usepackage{lmodern,amsmath,amssymb,amsthm,mathtools,mathrsfs}
\usepackage{microtype,enumitem,booktabs,array,longtable}
\usepackage{tikz,embedfile}
\usetikzlibrary{arrows.meta,calc,positioning}
\usepackage[colorlinks=true,linkcolor=blue!45!black,citecolor=blue!45!black,urlcolor=blue!45!black]{hyperref}
\usepackage[nameinlink,noabbrev]{cleveref}
\hypersetup{
 pdftitle={Escape--Covering Duality},
 pdfauthor={Zhipeng Deng},
 pdfsubject={Balanced supports, Minkowski worm constants, and restricted symplectic inequalities},
 pdfkeywords={convex geometry, balanced support, Minkowski worm, symplectic capacity, Wetzel problem}
}

\definecolor{keywordcolor}{rgb}{0.7,0.1,0.1}
\definecolor{commentcolor}{rgb}{0.4,0.4,0.4}
\definecolor{symbolcolor}{rgb}{0.0,0.1,0.6}
\definecolor{sortcolor}{rgb}{0.1,0.5,0.1}

\usepackage{listings}
\lstdefinelanguage{lean}{morekeywords={import,open,namespace,end,theorem,lemma,def,abbrev,structure,variable,by,exact,apply,have,let,if,then,else,fun,where,noncomputable},sensitive=true,morecomment=[l]{--},morestring=[b]"}

\lstdefinestyle{LeanCode}{
  language=lean,
  basicstyle=\ttfamily\fontsize{8pt}{8pt}\selectfont,
  keywordstyle=[1]{\color{keywordcolor}},
  keywordstyle=[2]{\color{sortcolor}},
  commentstyle=\itshape\color{commentcolor},
  stringstyle=\ttfamily,
  breaklines=true,
  breakatwhitespace=false,
  columns=fullflexible,
  keepspaces=true,
  showstringspaces=false,
  tabsize=2,
  extendedchars=false,
  captionpos=b
}

\usepackage{newunicodechar}
\newunicodechar{⦃}{\ensuremath{\{\!\{}}
\newunicodechar{⦄}{\ensuremath{\}\!\}}}
\newunicodechar{⋂}{\ensuremath{\bigcap}}
\newunicodechar{ₗ}{\ensuremath{_{l}}}
\newunicodechar{ᵢ}{\ensuremath{_{i}}}
\newunicodechar{ℝ}{\ensuremath{\mathbb{R}}}
\newunicodechar{ₜ}{\ensuremath{_{t}}}

\lstdefinestyle{codeappendix}{
  basicstyle=\ttfamily\scriptsize,
  stepnumber=1,
  numbersep=6pt,
  showstringspaces=false,
  breaklines=true,
  breakatwhitespace=false,
  columns=fullflexible,
  keepspaces=true,
  frame=single,
  xleftmargin=1.2em,
  framexleftmargin=0.8em,
  aboveskip=0.8em,
  belowskip=0.8em,
  literate={ℝ}{{$\mathbb{R}$}}1 {∧}{{$\wedge$}}1 {≤}{{$\le$}}1 {⊢}{{$\vdash$}}1 {⟨}{{$\langle$}}1 {⟩}{{$\rangle$}}1 {·}{{$\cdot$}}1 {×}{{$\times$}}1
}
\lstdefinestyle{LeanCode}{
  style=codeappendix,
  language={}
}

\lstdefinestyle{LeanCode}{
  style=codeappendix,
  language=lean,
  literate=
   {∀}{{$\forall$}}1 {∈}{{$\in$}}1 {→}{{$\to$}}1
   {≤}{{$\leq$}}1 {∧}{{$\wedge$}}1 {ℝ}{{$\mathbb{R}$}}1
   {⟨}{{$\langle$}}1 {⟩}{{$\rangle$}}1 {ℕ}{{$\mathbb{N}$}}1
   {·}{{$\cdot$}}1 {∑}{{$\sum$}}1 {∃}{{$\exists$}}1
   {↔}{{$\leftrightarrow$}}1 {≥}{{$\geq$}}1 {⊆}{{$\subseteq$}}1
   {←}{{$\leftarrow$}}1 {∩}{{$\cap$}}1 {₀}{{$_0$}}1
   {¬}{{$\neg$}}1 {≠}{{$\neq$}}1,
  keywordstyle=\color{keywordcolor},
  commentstyle=\itshape\color{commentcolor}
}

\lstdefinestyle{PythonCode}{
  style=codeappendix,
  language=Python
}

\numberwithin{equation}{section}
\newtheorem{theorem}{Theorem}[section]
\newtheorem{lemma}[theorem]{Lemma}
\newtheorem{proposition}[theorem]{Proposition}
\newtheorem{corollary}[theorem]{Corollary}
\theoremstyle{definition}

\newtheorem{remark}[theorem]{Remark}

\newcommand{\R}{\mathbb R}

\newcommand{\Sph}{\mathbb S}

\newcommand{\cC}{\mathcal C}

\newcommand{\HH}{\mathcal H}
\newcommand{\Int}{\operatorname{int}}
\newcommand{\conv}{\operatorname{conv}}
\newcommand{\diam}{\operatorname{diam}}
\newcommand{\dist}{\operatorname{dist}}
\newcommand{\len}{\operatorname{len}}
\newcommand{\vol}{\operatorname{vol}}
\newcommand{\area}{\operatorname{area}}
\newcommand{\per}{\operatorname{per}}

\newcommand{\supp}{\operatorname{supp}}
\newcommand{\SO}{\operatorname{SO}}
\newcommand{\OO}{\operatorname{O}}
\newcommand{\GL}{\operatorname{GL}}
\newcommand{\op}{\mathrm{op}}
\newcommand{\dH}{d_{\mathrm H}}
\newcommand{\Wet}{\mathcal W}

\newcommand{\eps}{\varepsilon}

\setlist{itemsep=2pt,topsep=4pt}
\title{Promotion and Extension of General Escape-Covering Duality Relationship}
\author{Zhipeng Deng}
\date{}

\begin{document}
\maketitle
\begin{abstract}
Bellman's lost-in-a-forest problem, Moser's worm problem, Lebesgue's universal covering problems,
and the classical "fits and covers" questions can be organized by a common obstruction
principle: an object fits a region if some admissible group action places it in the region, whereas
an escape object is one for which every placement fails to remain in the interior. 

In this paper, we develop a general theory of escape–covering duality for any dimension and explore several applications. For convex translation covers of all closed planar curves of Euclidean
length one, we prove
\[
 \Wet\ge f(\sqrt3/9)+c_1>0.15545785378555,c_1=0.000009944944227372...
\]
Here $f(\sqrt3/9)$ is Wetzel's classical lower bound. 

For a rational polytopal Minkowski gauge and prescribed rational
facet directions, we express universal closed-loop cover feasibility
as a rational polyhedron in the facet offsets. After removing
translations, minimum volume is attained at a vertex. Refining the
facet directions gives explicit two-sided rational bounds for the
unrestricted fixed-gauge worm constant. Four exact benchmarks recover
known sharp values. The construction uses an elementary
circuit--circumradius form of the established non-fitting-polygon
and capacity theory. We supply its support and ledger derivation,
the finite complete test family in arbitrary dimension, a cover
vertex bound, and the resulting algebraicity statement.

The surrounding escape--covering framework includes endpoint placement,
attainment, gauge stability, and finite-to-continuum error bounds.
For nonsymmetric gauges the objects are parametrized closed curves,
so orientation-dependent length is retained. Through the published
capacity identity, the improved Wetzel estimate gives a product
volume--capacity factor below $1.023781939711102$ for every planar
ellipsoidal gauge.
\end{abstract}



\setcounter{tocdepth}{2}
\tableofcontents

\section{Introduction}\label{sec:introduction}

Bellman’s lost-in-a-forest problem, Moser’s worm problem, Lebesgue's universal covering problems,
and the classical “fits and covers” questions can be organized by a common obstruction
principle: an object fits a region if some admissible group action places it in the region, whereas
an escape object is one for which every placement fails to remain in the interior. Bellman's lost-in-a-forest problem and its relation to universal covering
predate this article; see Finch--Wetzel \cite{FinchWetzel}.
We have the following relation to the present
work. Reference~\cite{Deng2412} transforms boundary configurations
into constrained path optimization; \cite{Deng2606} develops its
TSPN equivalence, convergence, and further escape variants.
References~\cite{Deng2607,Deng2608} study strip and triangle support
formulations. Reference~\cite{Deng2609} develops balanced support
calibrations and temporal ledgers for triangular covers. Those
formulations and elementary ingredients are credited here, with
the necessary proofs included. The present development concerns
the centered closed-gauge circuit formula, explicit facet-offset
optimization and rational brackets, and the quantitative Wetzel
width-deficit argument with its tangency refinement. 

Three distinctions are essential in this paper. Failure of placement in a body's
\emph{interior} permits placement in the closed body at the critical
threshold. Open curves with rotations and closed loops with translations
are different object--motion systems. For a nonsymmetric gauge,
length belongs to the parametrized loop: the opposite traversals in
\eqref{eq:asymmetric-example} have lengths three and six although
their traces coincide. These conventions match the worm formulation
in \cite{Rudolf}. Exact finite optimization is connected to the
continuum only through proved feasibility or limiting statements.

The auxiliary framework supplies endpoint placement, attainment,
metric hitting representations, Hausdorff stability, orientation
sampling with scale repair, and finite-witness convergence. These
results support the two principal estimates above. The symplectic
applications explicitly use the external capacity identity
\cite{RudolfCapacity}. The unrestricted Viterbo conjecture has a
counterexample \cite{HKO}; we consider restricted gauges and relaxed
ratios.

The smallest convex set containing a translate of every closed planar
curve of length one is unknown. Its area $\Wet$ is the translation-only
closed-curve constant studied by Wetzel \cite{Wetzel1973}; it differs
from Moser's problem for open arcs and congruent placement. Our main
numerical result is
\begin{equation}\label{eq:intro-wetzel}
 \Wet\ge f(r_0)+c_1>0.15545785378555,\qquad r_0=\sqrt3/9.
\end{equation}
The exact algebraic increment $c_1$ is specified in
\cref{thm:wetzel-improvement}. Wetzel's comparison body has a uniform
width deficit near the critical inradius. One deficient support line
is also tangent to the inball. Retaining this information gives
$c_1=0.000009944944227372\ldots$, more than $29\%$ larger than the
symmetric two-cap increment $c_0=0.000007681852018924\ldots$ derived
below. The proof separates geometric inequalities from certified
rational evaluation. It does not determine the minimum cover.

Our second focus is the outer optimization of Minkowski worm covers.
For a rational polytopal gauge $T$ and prescribed rational facet
directions, universal unit-loop covers are exactly the offsets in an
explicit rational polyhedron. A positive concave volume-root function
attains its minimum at a vertex after translation is removed
(\cref{lem:positive-concave,thm:fixed-normals}). For a normal mesh of
size $\eta$, its exactly computed optimum $V_T(U)$ satisfies
\begin{equation}\label{eq:intro-facet}
 \frac{V_T(U)}{(1+3R\eta/\rho)^d}\le W_T\le V_T(U),
\end{equation}
where $\rho,R$ are explicit inball and radius bounds
(\cref{cor:facet-hierarchy}). This is a complete rational certification
method, without a polynomial-time claim. \Cref{prop:facet-benchmarks}
works through four small cases with known sharp answers.

The computational bridge is the formula
\begin{equation}\label{eq:intro-circuits}
 E_T^{\rm tr}(K)=\min_a\frac{\sum_i h_K(a_i)}
 {\operatorname{rad}_T\{-a_1,-a_1-a_2,\ldots,-\sum_i a_i\}},
\end{equation}
over ordered balanced lists of at most $d+1$ vectors with
$\sum_i|a_i|=1$ (\cref{thm:circuit-radius}). Its proof optimizes a
directed ledger over its center and eliminates the contact vertices.
For polytopes it reduces the denominator to a circumradius linear
program with $d+1$ variables. The finite non-fitting-polygon
characterization and the capacity it computes are established results;
the role of this formula here is the explicit elimination and its
consequences for cover optimization. Balanced translation tests are
classical separation/Helly--Klee machinery, not a new containment
principle \cite{Klee,BezdekBezdek,AkopyanEtAl,KruppRudolf}.

\begin{table}[htbp]
\centering\small
\caption{Provenance and scope of the principal results.}
\label{tab:provenance}
\begin{tabular}{@{}>{\raggedright\arraybackslash}p{0.24\linewidth}>{\raggedright\arraybackslash}p{0.34\linewidth}>{\raggedright\arraybackslash}p{0.35\linewidth}@{}}
\toprule
Topic & Established input & Role in this article\\
\midrule
Translation and short polygons & Helly--Klee; Bezdek--Bezdek;
asymmetric billiard theory \cite{Klee,BezdekBezdek,AkopyanEtAl,KruppRudolf}
& Self-contained support and ledger proofs; explicit elimination.\\[4pt]
Finite complete test family & Planar family in \cite[Theorem~3.3]{BMP},
which also notes extension to higher dimensions
& Full $d$-dimensional proof, explicit cover vertex bound, algebraicity.\\[4pt]
Fixed facet directions & Circuit representation and
Brunn--Minkowski concavity & Rational offset polyhedron, vertex
minimum, quantitative rational mesh bounds.\\[4pt]
Wetzel lower bound & Wetzel's inradius and three-cap comparison
\cite{Wetzel1973} & Width deficit, tangency-sensitive added area,
explicit larger lower bound.\\[4pt]
Sharp model gauges & Simplex and parallelotope cases
\cite[Theorems~5.1--5.2]{Balitskiy2020}; planar results \cite{RudolfR4,BMP}
& Recovery and exact computational benchmarks; no new sharp values claimed.\\[4pt]
Capacity & Published capacity identity and symmetric
diagonal formula \cite{RudolfCapacity,ArtsteinKarasevOstrover}
& Improved ellipsoidal factor.\\
\bottomrule
\end{tabular}
\end{table}

The paper is organized as follows. Section~\ref{sec:foundations}
establishes the threshold, endpoint, metric, and quotient framework.
Section~\ref{sec:certificates} gives the balanced-support and ledger
proofs, circuit formulas, finite tests, and certified escape
approximations. Section~\ref{sec:cover-optimization} treats gauge
stability, cover attainment, polytopal minimizers, exact facet-offset
optimization, benchmarks, and rational approximation bounds.
Section~\ref{sec:wetzel} proves \eqref{eq:intro-wetzel} and the
finite-witness consequences. Section~\ref{sec:consequences} gives
geometric realizations, known sharp model cases, restricted capacity
estimates and metric search variants.
Section~\ref{sec:outlook} states the unresolved extremal questions.
Appendix~\ref{app:dependencies} maps the main theorem dependencies;
Appendix~\ref{app:arithmetic} explains the exact arithmetic and
reproducibility files.

\section{General duality and geometric foundations}\label{sec:foundations}
\subsection{Thresholds, endpoint placement, and homogeneity}\label{sec:threshold}
Throughout the Euclidean results $d\ge1$. A convex body is a compact
convex set with nonempty interior; lower-dimensional convex sets are
specified separately. We write $B_2^d=\{x:|x|\le1\}$,
$h_A(u)=\max_{a\in A}u\cdot a$ for a nonempty compact set $A$, and
$\mu_C(v)=\inf\{t>0:v\in tC\}$ for a convex body $C$ with
$0\in\Int C$. Volume is Lebesgue measure and $\dH$ is Hausdorff
distance, using the stated metric. Inball radii are positive, and net
meshes and Hausdorff error tolerances are nonnegative. All indexed
finite target samples are nonempty.
An object $A$ has a nonempty compact trace $\Gamma(A)$ in a space $X$ and a
size $\tau(A)\in[0,\infty]$. A rectifiable curve is a continuous map on a
compact interval, with traversal multiplicity counted in its length.
Injectivity is not imposed. In particular, a closed curve can traverse a
segment twice. For asymmetric length gauges the parametrized curve, rather
than its unoriented trace, is the object.

Let a group $G$ act by homeomorphisms on $X$. An object fits strictly in $F$
if $g\Gamma(A)\subset\Int F$ for some $g\in G$, and fits closedly if the
same inclusion holds with $F$ in place of $\Int F$. Put
\begin{align}
 E(F)&=\inf\{\tau(A):A\text{ has no strict placement in }F\},\label{eq:threshold}\\
 \rho^\circ(F)&=\sup\{s\ge0:\tau(A)<s\Longrightarrow
          A\text{ has a strict placement in }F\text{ for all }A\}.
\end{align}
We use $\inf\varnothing=+\infty$.

\begin{theorem}[Threshold identity]\label{thm:threshold}
For every object and placement system, $E(F)=\rho^\circ(F)$ in the
extended nonnegative reals.
\end{theorem}
\begin{proof}
Every object smaller than $E(F)$ fits strictly, since otherwise its size
would contradict the infimum defining $E(F)$. If $s>E(F)$ is finite,
the defining property of an infimum supplies a non-fitting object of size
less than $s$. Thus precisely the subcritical values, with the appropriate
endpoint convention, contribute to the same supremum. If $E(F)=\infty$,
every finite-size object fits and both sides are infinite. The case $E(F)=0$
follows by the same argument for every $s>0$.
\end{proof}

The identity is order-theoretic. Compactness is needed for the following
endpoint conclusion. To keep sequential arguments explicit, all groups in
the metric results are assumed locally compact and metrizable, with
continuous proper isometric actions. An action is proper when
$(g,x)\mapsto(x,gx)$ is a proper map.

\begin{theorem}[Closed placement at the endpoint]\label{thm:endpoint}
Let $F$ be compact in a proper metric space. Suppose $0<E(F)<\infty$.
If $\tau(A)=E(F)$ and objects $A_n$ satisfy
\[
 \tau(A_n)<E(F),\qquad \dH(\Gamma(A_n),\Gamma(A))\longrightarrow0,
\]
then $A$ has a closed placement in $F$.
\end{theorem}
\begin{proof}
Choose strict placements $g_n\Gamma(A_n)\subset F$. Choose
$a_n\in\Gamma(A_n)$ converging to $a\in\Gamma(A)$. Properness places all
$g_n$ in a compact transporter of the compact set
$\{a,a_1,a_2,\ldots\}$ into $F$. A subsequence converges to $g$.
The isometric action and uniform continuity on compact sets imply
\[
 \dH(g_n\Gamma(A_n),g\Gamma(A))\longrightarrow0.
\]
The limit is contained in the closed set $F$.
\end{proof}

\begin{corollary}[Critical objects and their subcritical approximants]
\label{cor:critical}
For the object class of all rooted rectifiable open paths, with length
as size, the hypotheses of \cref{thm:endpoint} imply that every path
of length at most $E(F)$ has a closed placement in $F$.
In Euclidean homogeneous systems the same conclusion holds for closed
loops and compact diameter objects whenever contractions are admissible.
\end{corollary}
\begin{proof}
Objects of size strictly below the threshold fit in the interior by
\cref{thm:threshold}. For a critical open path of length $L=E(F)$,
use its arclength parametrization on $[0,L]$ and truncate to
$[0,(1-1/n)L]$. The truncated length is less than $L$, and every
deleted point is within $L/n$ of the retained endpoint. Thus the
traces converge in Hausdorff distance and \cref{thm:endpoint} applies.
For loops and diameter objects, contraction about a trace point by
$1-1/n$ preserves the object class, scales size by that factor, and
converges uniformly on the compact trace. The endpoint theorem again
applies. In particular a critical non-fitting object can touch the
boundary in a closed placement without admitting a strict placement.
\end{proof}

\begin{lemma}[Infima and mutual downward cofinality]\label{lem:cofinal}
Let $S,T\subset\R$ be nonempty and bounded below. If every $s\in S$
has some $t\in T$ with $t\le s$, and every $t\in T$ has some $s\in S$
with $s\le t$, then $\inf S=\inf T$. Moreover, for any real $c$,
\[
 c\le\inf S\quad\Longleftrightarrow\quad c\le s\quad(s\in S).
\]
\end{lemma}
\begin{proof}
The first cofinality condition gives $\inf T\le t\le s$ for every
$s\in S$, so $\inf T\le\inf S$. The other condition gives the
reverse inequality. The final equivalence says precisely that $\inf S$
is the greatest lower bound: one direction uses $\inf S\le s$, and
the other uses the defining maximality among lower bounds. These two
inequalities also prove uniqueness of a greatest lower bound.
\end{proof}

\begin{theorem}[Homogeneous quotient]\label{thm:quotient}
Suppose the nonempty object and cover classes are preserved by positive
dilations, object sizes scale by $\lambda$, and the finite nonnegative
cover cost $\Phi(F)\in[0,\infty)$ scales by $\lambda^q$, $q>0$. Assume simultaneous dilation preserves strict and closed fitting,
and $0<E(F)<\infty$ for every admissible cover. Assume also:
\begin{enumerate}[label=(\roman*)]
\item every object of size at most $E(F)$ has a closed placement in $F$;
\item a closed universal unit cover strictly covers every object of size
less than one.
\end{enumerate}
Then the infimum of costs of closed universal unit covers is
\begin{equation}\label{eq:general-quotient}
 \mathfrak C=\inf_F\frac{\Phi(F)}{E(F)^q}.
\end{equation}
\end{theorem}
\begin{proof}
By (i), $F/E(F)$ is a universal unit cover and has cost
$\Phi(F)/E(F)^q$. Therefore $\mathfrak C$ is at most the infimum on the
right. If $U$ is a universal unit cover, (ii) and
\cref{thm:threshold} imply $E(U)\ge1$. Its quotient is consequently at
most $\Phi(U)$. Taking the infimum over such $U$ proves the reverse
inequality. Equivalently, the sets of normalized costs and universal-cover
costs are mutually downward cofinal, so \cref{lem:cofinal} applies.
\end{proof}

\begin{lemma}[Convex strictification]\label{lem:strictification}
Let $K\subset\R^d$ be convex with nonempty interior, and let the placement
group consist of translations and a compact set of linear maps. Suppose the object family and size are homogeneous.
If $K$ is a closed universal unit cover, every object of size less than one
has a strict placement, provided size-zero objects are singletons.
\end{lemma}
\begin{proof}
Choose $y\in\Int K$ and $\rho>0$ with $y+\rho B_2^d\subset K$.
For $z\in K$,
$0\le t<1$, and $|v|<\rho$, convexity gives
\[
 (1-t)y+tz+(1-t)v=(1-t)(y+v)+tz\in K.
\]
Thus the open ball of radius $(1-t)\rho$ about $(1-t)y+tz$
lies in $K$, proving $y+t(K-y)\subset\Int K$. For an object of size
$0<\ell<1$, expand by $1/\ell$, place the expanded object in $K$, and
contract this placement toward $y$ by $\ell$. The two scale factors cancel
and leave an admissible placement of the original object. Singletons fit
in the interior directly.
\end{proof}

\subsection{Reversal of motion and metric hitting paths}\label{sec:metric}
We use the standard length-space terminology of \cite{BBI}.
\begin{proposition}[Motion reversal, boundary hitting, and compact states]
\label{prop:motion}
Let a group act by homeomorphisms. For arbitrary sets $A,B,F$,
\begin{align}
 gA\subset F&\quad\Longleftrightarrow\quad A\subset g^{-1}F,
       \label{eq:reversal-subset}\\
 gA\cap B\ne\varnothing&\quad\Longleftrightarrow\quad
       A\cap g^{-1}B\ne\varnothing.\label{eq:reversal}
\end{align}
If $F$ is closed, $\gamma$ is continuous, and $g\gamma(0)\in F$,
then $g\Gamma(\gamma)\not\subset\Int F$ if and only if
$\Gamma(\gamma)\cap g^{-1}\partial F\ne\varnothing$.
For compact $F$ and a proper action on a proper metric space, the
configuration set $\cC_F=\{g:go\in F\}$ is compact, and the family
$N_g=g^{-1}\partial F$, $g\in\cC_F$, has a closed graph.
\end{proposition}
\begin{proof}
A bijection satisfies $gx\in F$ if and only if $x\in g^{-1}F$;
applying this pointwise proves \eqref{eq:reversal-subset}. A point
in $gA\cap B$ is $ga$ for some $a\in A$ with $a\in g^{-1}B$,
and applying $g$ proves the converse implication in \eqref{eq:reversal}.

Write $\beta=g\circ\gamma$. If $\beta(0)\in\partial F$, there is
already a boundary hit. Otherwise $\beta(0)\in\Int F$. When
$\beta$ fails interior containment, the nonempty closed set
$\{t\in[0,1]:\beta(t)\notin\Int F\}$ has a least member $t_*>0$.
Continuity and closedness of $F$ imply $\beta(t_*)\in F$, since
$\beta(t)\in\Int F$ for $t<t_*$. Hence
$\beta(t_*)\in F\setminus\Int F=\partial F$. Conversely a boundary
hit directly precludes interior containment. No convexity or boundary
smoothness is needed.

The proper map $(g,x)\mapsto(x,gx)$ has compact inverse image of
$\{o\}\times F$. Its projection onto the group is exactly $\cC_F$,
so this configuration set is compact. Finally, if $g_n\to g$,
$x_n\to x$, and $x_n\in N_{g_n}$, continuity gives
$g_nx_n\to gx$. Since $g_nx_n\in\partial F$ and the boundary is
closed, $gx\in\partial F$, or $x\in N_g$. This proves closedness
of the graph. Empty targets would preclude a feasible hitting path;
the next theorem explicitly assumes nonempty targets and feasibility.
\end{proof}

Here escape means reaching the boundary. The finite-length feasibility
assumption below is separate from compactness of the configurations.

\begin{theorem}[Closed-graph metric TSPN and convergence]\label{thm:metric}
Let $X$ be a proper geodesic metric space, $o\in X$, and $C$ a nonempty
compact metric space. Let $N_c\subset X$ be nonempty closed targets with closed
graph
\[
 \{(c,x):x\in N_c\}\subset C\times X.
\]
Assume one finite-length path from $o$ meets every $N_c$. Let $V$ be the
minimum-length infimum for this continuum hitting problem. Then:
\begin{enumerate}[label=(\alph*)]
\item $V$ is attained.
\item For $C_m=\{c_1,\ldots,c_m\}$, the finite optimum equals
\begin{equation}\label{eq:tspn}
 V_{C_m}=\min_{\pi,p_i\in N_{c_i}}
 \left[d(o,p_{\pi(1)})+
       \sum_{j=2}^{m}d(p_{\pi(j-1)},p_{\pi(j)})\right],
\end{equation}
where $\pi$ ranges over permutations. The minimum is attained.
\item If finite sets $C_m$ are nested with dense union, then
$V_{C_m}\uparrow V$. Arclength-normalized finite minimizers are precompact
in the uniform topology, and every convergent subsequence with indices
tending to infinity has a continuum-minimizing limit.
\end{enumerate}
\end{theorem}
\begin{proof}
Normalize a minimizing sequence on $[0,1]$ by arclength. A common length
bound $L$ gives common Lipschitz constant $L$ and confines every trace to
$\overline B(o,L)$, which is compact. Arzel\`a--Ascoli supplies a uniform
limit. For each fixed target, choose a hitting time and pass to a
subsequence of times. The limit path still hits the target. Length is
lower semicontinuous because each finite partition sum passes to the
limit, and length is the supremum of these sums. This proves (a).

For a finite family, order any chosen hitting times along a feasible
curve. The distances between successive selected points are bounded by
the intervening subarc lengths, giving one inequality in
\eqref{eq:tspn}. Concatenated minimizing geodesics give the reverse
inequality. An objective-bounded minimizing sequence of selected points
stays in a common ball about $o$, even when targets are unbounded. The
closed target intersections with this ball are compact, so the finite
minimum is attained.

Nested constraints give $V_{C_m}\uparrow\bar V\le V$. Normalize finite
minimizers and take a uniform limit along a subsequence $m_k\to\infty$.
The mesh of $C_m$ tends to zero by compactness and density. Given $c\in C$,
choose $c_k\in C_{m_k}$ converging to $c$ and corresponding hitting
times $t_k$. A subsequence of $t_k$ converges, and uniform convergence
gives $\gamma_{m_k}(t_k)\to\gamma(t)$. The closed target graph implies
$\gamma(t)\in N_c$. Hence
\[
 V\le\len(\gamma)\le\liminf_kV_{C_{m_k}}=\bar V\le V.
\]
This proves (c), including its assertion about every cluster point.
\end{proof}

\begin{corollary}[Variational formulation]\label{cor:gamma}
Under \cref{thm:metric}, let $C_m$ be nested finite samples with dense
union, as in part~(c), and define on continuous rooted paths
\[
 F_m(\gamma)=\len(\gamma)+I_m(\gamma),\qquad
 F(\gamma)=\len(\gamma)+I(\gamma),
\]
where each indicator is zero for the respective feasible paths and
$+\infty$ otherwise. Then $F_m$ $\Gamma$-converges to $F$ in the uniform
topology.
\end{corollary}
\begin{proof}
For a uniformly convergent sequence with finite $\liminf F_m$, pass
to a subsequence realizing this lower limit and having finite energies.
Its paths satisfy each fixed sampled constraint eventually. The
closed-graph argument in \cref{thm:metric} gives continuum feasibility
of the limit. Lower semicontinuity of length proves the $\liminf$
inequality. For a continuum-feasible $\gamma$, the constant sequence
is a recovery sequence since every sampled constraint is satisfied.
For $F(\gamma)=\infty$ the recovery inequality imposes no restriction.
\end{proof}
\begin{remark}
Bounded length alone does not make arbitrary parametrizations
equicontinuous. The compactness assertion uses arclength normalization;
it is not an assertion of uniform-topology equicoercivity without that
normalization. A sampled feasible incumbent bounds the sampled problem
from above, but need not bound the continuum problem from above.
\end{remark}

\subsection{Euclidean length and diameter dualities}\label{sec:euclidean}
Let $H\le\OO(d)$ be a compact subgroup. Define $E_H(K)$ as the least
Euclidean length of an open curve having no placement $z+R\Gamma\subset
\Int K$ with $R\in H$. Define $D_H(K)$ analogously using the diameter of
a nonempty compact set. We write $M_{d,H}$ and $L_{d,H}$ for the minimum
volume infima of convex covers of, respectively, all open curves of
length at most one and all diameter-at-most-one compact sets under these
motions.
The physical heading convention is $H=\SO(d)$. Allowing reflections gives
$H=\OO(d)$; the notation keeps these conventions distinct.

\begin{lemma}[Convexification and strict containment]\label{lem:convex-hull}
For every nonempty compact set $A\subset\R^d$, its convex hull is
compact and $\diam(\conv A)=\diam(A)$. If $K$ is a convex body,
then $A\subset\Int K$ if and only if $\conv A\subset\Int K$.
Convexification commutes with affine maps.
\end{lemma}
\begin{proof}
Carath\'eodory's theorem expresses each point of $\conv A$ as a
convex combination of at most $d+1$ points of $A$. It therefore
expresses $\conv A$ as the image of a compact product with a compact
simplex, proving compactness. For $x=\sum_i\alpha_i a_i$ and
$y=\sum_j\beta_j b_j$ in the hull,
\[
 |x-y|=\left|\sum_{i,j}\alpha_i\beta_j(a_i-b_j)\right|
 \le\sum_{i,j}\alpha_i\beta_j|a_i-b_j|\le\diam A.
\]
The reverse diameter inequality follows from $A\subset\conv A$.
If $A\subset\Int K$, compactness supplies $\eps>0$ such that
$A+\eps B_2^d\subset K$. Convexifying gives
$\conv A+\eps B_2^d\subset K$, so the hull is also in the interior.
The converse is immediate. Affine maps preserve convex combinations.
\end{proof}

\begin{theorem}[Euclidean endpoints and quotient identities]\label{thm:euclidean}
For every convex body $K$ with nonempty interior, both thresholds are
positive, finite, and attained. In the diameter problem the minimizing
obstruction may be chosen convex. Every critical object admits a closed
placement. If $K$ contains a Euclidean ball of radius $r$, then
\[
 2r\le E_H(K)\le\diam K,\qquad r\le D_H(K)\le\diam K.
\]
Moreover,
\begin{equation}\label{eq:euclidean-quotients}
 M_{d,H}=\inf_K\frac{\vol_d(K)}{E_H(K)^d},\qquad
 L_{d,H}=\inf_K\frac{\vol_d(K)}{D_H(K)^d}.
\end{equation}
\end{theorem}
\begin{proof}
A curve of length $\ell<2r$, centered at an arclength midpoint,
lies in the open ball of radius $r$: both remaining subarcs have length
at most $\ell/2$. A compact set of diameter less than $r$ lies in
that ball after translating one of its points to the center. A segment
of length $\diam K$ cannot fit strictly: if both endpoints were in
$\Int K$, moving them slightly apart would contradict the definition
of diameter. This proves the displayed bounds. Normalize non-fitting curves to start at the
origin and parameterize by arclength. Bounded length gives a compact
uniformly convergent subsequence. Strict fitting of its limit would
persist under sufficiently small uniform perturbations, contradicting
non-fitting of the sequence. Length lower semicontinuity then proves
attainment.

For diameter objects, \cref{lem:convex-hull} preserves diameter and
fitting in a convex interior. Translate each convex competitor to contain zero
and use Blaschke selection in a common ball. Diameter is continuous and
strict fitting is stable, proving attainment of $D_H(K)$.
Endpoint closed placements follow by contractions and compactness of $H$
and of translating vectors. Finally, convex strictification and
\cref{thm:quotient} prove \eqref{eq:euclidean-quotients}.
\end{proof}

The curve bounds $2r\le E_H(F)\le\diam F$ also hold for a
compact, possibly nonconvex Euclidean forest $F$ containing an $r$-ball.
The same midpoint and segment arguments prove them. Thus a finite escape
path exists for every such forest; \cref{prop:motion,thm:metric} then give
the exact boundary-hitting representation and recovery from dense finite
configuration sets. Convexity enters the quotient and diameter reduction,
not these open-path existence and hitting statements.

\subsection{Minkowski loops and the translative quotient}\label{sec:minkowski}
Let $T\subset\R^d$ be a convex body with $0\in\Int T$ and
\[
 T^\circ=\{v:v\cdot t\le1\ \forall t\in T\}.
\]
For a nonsymmetric gauge, traversal orientation is part of the object.
Two parametrized closed curves can have the same trace and different
lengths; see \eqref{eq:asymmetric-example}. Reparametrizations preserving
orientation preserve length. Throughout this subsection translations are
the allowed placements, while traces are used only for containment.
The length of an absolutely continuous parametrized curve is
\begin{equation}\label{eq:gauge-length}
 \ell_T(q)=\int_0^1h_T(\dot q(t))\,dt
          =\int_0^1\mu_{T^\circ}(\dot q(t))\,dt.
\end{equation}
The equality follows directly from the definition of the polar.
For an arbitrary rectifiable parametrized curve, use its arclength
representative with the original traversal order in \eqref{eq:gauge-length};
pauses have no effect. All loops in the Minkowski covering classes are
rectifiable, and a zero-length loop is constant because $h_T(v)\ge c|v|$
for some $c>0$. In particular, length ratios below are taken only for
nonconstant rectifiable loops.
We reserve $E_T^{\rm tr}(K)$ for the least $T$-length of a
\emph{closed} loop whose trace has no strict translate into $K$, and
$W_T$ for the least-volume infimum of convex translation covers for
all such loops of length at most one. A superscript is retained on
$E_T^{\rm tr}$ to avoid confusion with rotational open-path $E_H$.

\begin{theorem}[Minkowski quotient and finite polygon reduction]\label{thm:minkowski}
For each convex body $K$, $E_T^{\rm tr}(K)$ is positive, finite, and
attained by a polygonal loop with at most $d+1$ vertices. Every loop
of $T$-length at most this threshold has a closed translate in $K$, and
\begin{equation}\label{eq:minkowski-quotient}
 W_T=\inf_K\frac{\vol_d(K)}{E_T^{\rm tr}(K)^d}.
\end{equation}
\end{theorem}
\begin{proof}
There are $0<c\le C$ with $c|v|\le h_T(v)\le C|v|$.
An inball gives a positive lower bound for the threshold; a long
doubled segment gives a finite upper bound. Translate each loop so that
$q(0)=0$, then parametrize on $[0,1]$ by Euclidean arclength. Norm
equivalence gives a uniform Lipschitz bound and confines the traces to a
common ball, so Arzel\`a--Ascoli gives uniform precompactness.
For completeness, gauge length is lower semicontinuous under uniform
convergence: it equals the supremum of the polygonal sums
$\sum_jh_T(q(t_j)-q(t_{j-1}))$. The inequality from sums to the integral
follows by convexity; equality follows by approximating the derivative
in $L^1$ by step functions and using the Lipschitz continuity of $h_T$.
Non-fitting is closed under this convergence, so the threshold is
attained. Apply \cref{thm:helly} to obtain the polygon.

Contract a critical loop by factors increasing to one. Each contraction
has a strict translate, and compactness of its initial image in $K$
bounds the translating vectors. Passing to the limit proves endpoint
coverage. Convex strictification and \cref{thm:quotient} then give
\eqref{eq:minkowski-quotient}.
\end{proof}

\begin{corollary}[Pointwise escape inequalities]\label{cor:pointwise}
For a fixed gauge $T$ and a real constant $c\ge0$,
\begin{equation}\label{eq:pointwise}
 W_T\ge c\quad\Longleftrightarrow\quad
 \vol_d(K)\ge c\,E_T^{\rm tr}(K)^d
 \quad\hbox{for every convex body }K.
\end{equation}
\end{corollary}
\begin{proof}
The quotients in \eqref{eq:minkowski-quotient} are positive. Apply the
greatest-lower-bound equivalence in \cref{lem:cofinal} to that quotient
set and multiply by the positive number $E_T^{\rm tr}(K)^d$.
\end{proof}

\begin{remark}[Norm-preserving linear motions]\label{rem:norm-motions}
Let $H\le\GL(d)$ be compact with $R(T^\circ)=T^\circ$ for every
$R\in H$. If translations and these maps are allowed, the analogous
closed-loop threshold $E_{T,H}(K)$ and universal-cover constant
$W_{T,H}$ satisfy
\[
 W_{T,H}=\inf_K\frac{\vol_d(K)}{E_{T,H}(K)^d}.
\]
Indeed $h_T(Rv)=h_T(v)$, so these motions preserve length. Small loops
fit in an inball. Compactness of $H\subset\GL(d)$ bounds its least
singular values away from zero, so a sufficiently long doubled segment
cannot fit after any allowed motion. Thus the threshold is positive and
finite. After translating the initial points to zero, norm equivalence
bounds the Euclidean lengths of a minimizing sequence.
Arclength compactness and stability of strict fitting give attainment.
Critical contractions admit placements with bounded translations and a
convergent subsequence in $H$, giving closed endpoint placement.
\Cref{lem:strictification,thm:quotient} prove the formula. The
$d+1$-vertex conclusion of \cref{thm:minkowski} is a translative
statement and is not used for this larger motion group.
\end{remark}

\section{Balanced supports, directed ledgers, and escape computation}\label{sec:certificates}
\subsection{Balanced supports and finite obstruction sets}\label{sec:support}
The containment principle below is a classical consequence of convex
separation, subdifferentials, and Carath\'eodory's theorem, in the
Helly--Klee translation tradition \cite{Klee,Rockafellar,Schneider}.
We prove its signed-margin form to fix the strict-containment
conventions and the certificates used later.
For a nonempty compact set $A$, write
$h_A(u)=\max_{x\in A}u\cdot x=h_{\conv A}(u)$.
A normalized balanced collection consists of unit vectors $u_i$ and
weights $\lambda_i\ge0$ such that
\begin{equation}\label{eq:balance}
 \sum_i\lambda_i=1,\qquad\sum_i\lambda_i u_i=0.
\end{equation}
Zero weights can be deleted. Below, at most $d+1$ directions are needed.

\begin{theorem}[Signed translation margin]\label{thm:balance}
Let $K$ be a convex body and $A$ nonempty and compact. Set
\[
 m_K(A)=\min_{z\in\R^d}\max_{u\in\Sph^{d-1}}
       \{z\cdot u+h_A(u)-h_K(u)\}.
\]
Then
\begin{equation}\label{eq:margin-dual}
 m_K(A)=\max_{\substack{2\le s\le d+1\\(u,\lambda)\text{ satisfies }\eqref{eq:balance}}}
       \sum_{i=1}^s\lambda_i\bigl(h_A(u_i)-h_K(u_i)\bigr).
\end{equation}
A translate of $A$ lies in $K$ exactly when $m_K(A)\le0$ and in $\Int K$
exactly when $m_K(A)<0$.
\end{theorem}
\begin{proof}
The maximum defining the objective is continuous and at least
$|z|-\max_u|h_A(u)-h_K(u)|$. It is therefore coercive and has a minimizer
$z_*$. If zero were outside the convex hull of its active unit directions,
strict separation would provide $v$ and $a>0$ with $u\cdot v\le-a$
for every active direction. This inequality with $-a/2$ holds on a
neighborhood of the compact active set. On the compact complement, the
objective values at $z_*$ lie below the maximum by some $b>0$.
Replacing $z_*$ by $z_*+tv$ decreases every value in the neighborhood
by at least $ta/2$, while each remaining value increases by at most
$t|v|$. For $0<t<b/(2|v|)$ the new maximum is strictly smaller,
contradicting minimality.
Hence zero belongs to that convex hull. Carath\'eodory's theorem gives
balanced active directions of cardinality at most $d+1$.
Summing their active equalities cancels $z_*$ and gives the right side
of \eqref{eq:margin-dual}. Any balanced average is no larger than the
maximum at any translation, giving the opposite inequality.

Nonpositive margin is precisely support-function containment.
A negative margin gives a uniform positive distance from all supporting
hyperplanes. Conversely, a compact set in $\Int K$ has a positive
distance from the complement, giving a strictly negative maximum.
\end{proof}

\begin{corollary}[Fit factor and orientation elimination]\label{cor:fit}
Assume $rB_2^d\subset K$ with $r>0$. For $R\in\OO(d)$ define
\[
 s_K(A,R)=\inf\{t\ge0:\exists z,\ z+RA\subset tK\}.
\]
Then
\begin{equation}\label{eq:fit-factor}
 s_K(A,R)=\max_{\substack{2\le s\le d+1\\(u,\lambda)\text{ balanced}}}
 \frac{\sum_i\lambda_i h_A(R^Tu_i)}
      {\sum_i\lambda_i h_K(u_i)}.
\end{equation}
All denominators are at least $r$. The function is translation invariant
in $A$, homogeneous of degree one in $A$, and continuous in $(A,R)$.
Moreover, with $H$ the allowed orientation group,
\begin{equation}\label{eq:rot-escape}
 A\text{ has no strict placement in }K
 \quad\Longleftrightarrow\quad
 \min_{R\in H}s_K(A,R)\ge1.
\end{equation}
\end{corollary}
\begin{proof}
For $t>0$, apply \cref{thm:balance} to the full-dimensional body
$tK$. Its balanced inequalities say exactly that $t$ is at least every
ratio in \eqref{eq:fit-factor}. Taking the infimum over positive $t$
gives the same fitting infimum as allowing $t=0$: a set that fits in
$0K=\{0\}$ is a singleton and fits in $tK$ for every $t>0$. The numerator
of each ratio is nonnegative, since translating any point of $A$ to
zero makes each support nonnegative without changing the balanced sum.
The compact parameter set of balanced collections and the positive
denominator give attainment and continuity. If $z+RA\subset\Int K$,
there is $\eps>0$ with $h_{z+RA}\le h_K-\eps$ on the unit sphere.
Choose $0<t<1$ with $(1-t)\max_{|u|=1}h_K(u)<\eps$; then
$h_{z+RA}\le th_K$, so the same placement lies in $tK$.
Conversely, $tK\subset\Int K$ for $0\le t<1$. This proves
\eqref{eq:rot-escape}.
\end{proof}

This gives the exact convex Bellman formulation
\begin{equation}\label{eq:support-bellman}
 E_H(K)=\min_\gamma\left\{\len(\gamma):
         \min_{R\in H}s_K(\conv\Gamma_\gamma,R)\ge1\right\}.
\end{equation}
It eliminates translations without prescribing the temporal order of
support contacts. For a triangle there is one positive facet circuit;
for general polytopes several circuits compete.

\begin{proposition}[Polytope circuits]\label{prop:poly-circuit}
Suppose
\[
 K=\{x:a_i\cdot x\le b_i,\ 1\le i\le m\}
\]
is bounded with nonempty interior, with every $a_i\ne0$. Define
$C_K=\{\lambda\in\R_+^m:\sum_i\lambda_i a_i=0\}$.
At fixed orientation $R$, failure of a strict translating placement
of $RA$ is equivalent to
\begin{equation}\label{eq:poly-nonfit}
 \sum_i\lambda_i h_A(R^Ta_i)\ge\sum_i\lambda_i b_i
\end{equation}
for some extreme ray of $C_K$. Each such ray uses at most $d+1$ normals.
\end{proposition}
\begin{proof}
Minimizing the largest facet violation over translations is a coercive
finite maximum, because the facet normals positively span $\R^d$.
The active-set argument of \cref{thm:balance}, with normalization
$\sum_i\lambda_i=1$, gives the strict alternative
\eqref{eq:poly-nonfit}. Decompose a balanced vector into extreme rays.
If its violation is nonnegative, at least one constituent ray has
nonnegative violation. Finally, a ray supported on more than $d+1$
columns has a kernel of dimension at least two on that support.
Small positive and negative perturbations along an independent kernel
vector preserve balance and nonnegativity and contradict extremality.
\end{proof}

\begin{proposition}[Regular even-polygon circuit count]\label{prop:circuit-count}
The facet-normal configuration of a regular $2m$-gon, $m\ge2$,
has $m$ positive two-normal circuits and
$m(m-1)(m-2)/3$ positive three-normal circuits. Its total number
of extreme positive balanced rays is therefore
\[
 m+\frac{m(m-1)(m-2)}3.
\]
This counts circuits, not the additional temporal assignments in an
escape optimization.
\end{proposition}
\begin{proof}
The two-normal circuits are the $m$ antipodal pairs. For a positive
triple, its three positive cyclic gaps $a,b,c$, measured in normal
indices, sum to $2m$, and each is less than $m$; this is exactly
the condition that the origin lies in the triangle's interior.
There are $\binom{2m-1}{2}$ ordered positive gap triples before
this restriction. For a specified gap, $\binom m2$ violate it
by being at least $m$, and two gaps cannot simultaneously violate it.
Thus the admissible ordered gap count is
\[
 \binom{2m-1}{2}-3\binom m2=\frac{(m-1)(m-2)}2.
\]
Multiplication by $2m$ choices of initial normal and division by the
three choices of initial vertex of the same circuit give the result.
\end{proof}

The finite obstruction principle is classical \cite{Klee}.
The short non-fitting polygon reduction belongs to the shortest-billiard
theory of Bezdek--Bezdek and its asymmetric extensions
\cite{BezdekBezdek,AkopyanEtAl,KruppRudolf}; the following proof
records the exact open/closed conventions needed here.

\begin{theorem}[Finite-point obstruction and polygon reduction]\label{thm:helly}
If a compact set $A$ has no strict translate into a convex body $K$, a
subset of at most $d+1$ points already has no strict translate.
Consequently, for any finite convex positively homogeneous integrand
$g:\R^d\to[0,\infty)$ with $g(v)>0$ for $v\ne0$,
the shortest translatively non-fitting closed curve can be chosen to
be a closed polygon with at most $d+1$ vertices. The analogous open
curve can be chosen as a polygonal chain through at most $d+1$ points.
The resulting classes of polygons also suffice to test closed universal coverage.
\end{theorem}
\begin{proof}
By \cref{thm:balance}, non-fitting gives balanced directions with
$\sum_i\lambda_i(h_A(u_i)-h_K(u_i))\ge0$.
Choose $p_i\in A$ attaining $h_A(u_i)$. A strict translating placement
of these selected points would make every inequality
$u_i\cdot(z+p_i)<h_K(u_i)$ strict, contradicting the balanced sum.

For a curve, select times at these points in their actual traversal
order. Replace intervening subarcs by chords. Convexity and positive
homogeneity give
\[
 g(\gamma(b)-\gamma(a))\le\int_a^b g(\dot\gamma(t))\,dt.
\]
For a loop, include the closing cyclic interval; for an open curve,
discard the initial and final unused subarcs. The selected polygon or
chain is still non-fitting and is no longer than the original curve.
Existence follows by the same compactness argument as
\cref{thm:euclidean}, using $g(v)\ge c|v|$ for some $c>0$.

For closed universal coverage, if a trace has no closed translate,
its margin in \cref{thm:balance} is positive. The same construction
gives selected points, and hence a polygon, with positive margin.
Thus covering every such finite polygon covers every original curve.
\end{proof}
\begin{remark}
The finite-point principle is a form of the classical Helly--Klee
translation principle \cite{Klee}. Its closed-loop consequences are
central to the billiard literature \cite{BezdekConnelly,BezdekBezdek}.
The proof above is included to make the computational reductions
self-contained. For rotations, different orientations may require
different selected points. A single $d+1$-point certificate does not
follow for all orientations simultaneously.
\end{remark}

\subsection{Directed ledgers and optimal closed-path calibration}\label{sec:ledger}
Balanced support inequalities can be integrated or summed to certify
escape lower bounds, as in \cite{Deng2609}. The order of support
contacts along the path is essential. Their angular order on a convex
hull need not be their temporal order on an open path.

\begin{theorem}[Directed ledger in a normed space]\label{thm:normed-ledger}
Let $V$ be a real normed space, $a_1,\ldots,a_m\in V^*$ continuous
linear functionals with $\sum_j a_j=0$, and let $P_0,\ldots,P_N\in V$
with $N\ge1$.
Choose contact indices $k_j\in\{0,\ldots,N\}$ and put
$v_i=P_i-P_{i-1}$ and $R_i=\sum_{j:k_j\ge i}a_j\in V^*$. Then
\begin{equation}\label{eq:normed-ledger}
 \sum_j a_j(P_{k_j})=\sum_i R_i(v_i)
 \le \max_i\|R_i\|_*\sum_i\|v_i\|.
\end{equation}
The same bound holds for finitely many contacts on a rectifiable path,
with its length in place of the polygonal sum.
\end{theorem}
\begin{proof}
Use $P_{k_j}=P_0+\sum_{i\le k_j}v_i$. Balance cancels
$\sum_j a_j(P_0)$, and interchanging the remaining finite sums yields
$\sum_j\sum_{i\le k_j}a_j(v_i)=\sum_iR_i(v_i)$. By the definition
of the dual norm, $R_i(v_i)\le\|R_i\|_*\|v_i\|$ for every $i$;
summing proves \eqref{eq:normed-ledger}. For a rectifiable path, order
its finitely many contact times and use the corresponding chord
increments. Each chord norm is bounded by the intervening subarc length,
so their sum is at most the full length. Neither finite dimensionality
nor angular ordering of contacts is required.
\end{proof}

In finite dimension one can sharpen the constant for closed paths and
allow asymmetric gauges. Let a polygonal path have vertices
$P_0,\ldots,P_N$, $N\ge1$, and increments
$v_i=P_i-P_{i-1}$. Suppose vectors $a_1,\ldots,a_m$ are assigned to
contact indices $k_j\in\{0,\ldots,N\}$, with
$\sum_j a_j=0$. A weighted support sum has this form when
$a_j=\lambda_j u_j$ and $P_{k_j}$ maximizes the support in direction
$u_j$. Define the suffix states
\begin{equation}\label{eq:ledger-states}
 R_i=\sum_{j:k_j\ge i}a_j,\qquad 1\le i\le N.
\end{equation}

\begin{theorem}[Ledger identity and best closed-path constant]\label{thm:ledger}
For the gauge length $g=h_T$,
\begin{equation}\label{eq:ledger-identity}
 \sum_{j=1}^m a_j\cdot P_{k_j}
   =\sum_{i=1}^N R_i\cdot v_i
   \le\left(\max_i\mu_T(R_i)\right)\sum_{i=1}^N g(v_i).
\end{equation}
If the path is closed, the stronger estimate is
\begin{equation}\label{eq:centered-ledger}
 \sum_j a_j\cdot P_{k_j}
 \le C_T(R)\,\ell_T(P),\qquad
 C_T(R)=\min_{b\in\R^d}\max_i\mu_T(R_i-b).
\end{equation}
The minimum is attained, and $C_T(R)$ is the smallest nonnegative
constant valid for all closed increment sequences with these fixed states:
\begin{equation}\label{eq:ledger-optimal}
 C_T(R)=\max\left\{\sum_iR_i\cdot v_i:
         \sum_i v_i=0,\ \sum_i g(v_i)\le1\right\}.
\end{equation}
The contact-realization constraints of a particular geometric problem
can further reduce this maximum; they are not assumed in
\eqref{eq:ledger-optimal}.
\end{theorem}
\begin{proof}
Insert $P_{k_j}=P_0+\sum_{i\le k_j}v_i$ and interchange the finite
sums. Balance cancels the $P_0$ term. The polar inequality
$w\cdot v\le\mu_T(w)h_T(v)$ proves
\eqref{eq:ledger-identity}. If $\sum_i v_i=0$, replacing every
$R_i$ by $R_i-b$ leaves the sum unchanged, proving the upper
bound in \eqref{eq:centered-ledger}. Its minimizing objective is
continuous and coercive, so it attains a minimum.

Here is a finite-dimensional proof that the constant is optimal.
Let
\[
 D=\conv\bigcup_{i=1}^N
       \{(v,R_i\cdot v):v\in T^\circ\}\subset\R^{d+1}.
\]
This is compact and convex. Its points are precisely sums
$\sum_i(v_i,R_i\cdot v_i)$ with $\sum_i g(v_i)\le1$:
one direction follows by convexity of $g$; the other writes each
nonzero $v_i$ as $g(v_i)$ times a point of $T^\circ$, and uses
the remaining weight on zero. Let $M$ be the largest $t$ with
$(0,t)\in D$. Thus $M$ is the maximum in
\eqref{eq:ledger-optimal}. For every $\eps>0$, strictly separate
$(0,M+\eps)$ from $D$ by a functional $(p,q)$.
Comparison with $(0,M)\in D$ forces $q>0$. Dividing by $q$ gives
\[
 M+\eps>
 \max_i\max_{v\in T^\circ}(R_i+p/q)\cdot v
 =\max_i\mu_T(R_i+p/q)\ge C_T(R).
\]
Letting $\eps\downarrow0$ and combining with the already proved
upper bound yields $M=C_T(R)$. If $N=1$, closure forces $v_1=0$
and the formula gives $C_T(R)=0$; the nonnegative convention is
necessary here because every real constant would satisfy the zero
inequality. For $N\ge2$, a nonzero closed increment sequence exists,
and applying a proposed bound to it and its negative forces that
constant to be nonnegative.
\end{proof}

\begin{remark}[When centering is a strict improvement]
For a symmetric gauge and states $\{0,w\}$, the uncentered constant
is $\mu_T(w)$ whereas the centered constant is
$\mu_T(w)/2$, attained at $b=w/2$. The lower bound follows from
the triangle inequality. The improvement is a consequence of closure,
not an assumption that contacts occur in angular order. For a general
rectifiable curve, select the finitely many contact times and apply
the same identity to its intervening increments, whose total gauge
length is no larger than the original length.
\end{remark}

\begin{corollary}[Calibration lower bounds]\label{cor:calibration}
Suppose every feasible escape path admits a certificate of a type
$\alpha$ for which the support sum $S_\alpha$ and a constant
$C_\alpha>0$ satisfy
\[
 B_\alpha\le S_\alpha\le C_\alpha\,\len(\gamma),
 \qquad B_\alpha>0.
\]
Then $\len(\gamma)\ge B_\alpha/C_\alpha$ and the escape threshold
is at least $\inf_\alpha B_\alpha/C_\alpha$. One may use the dual
suffix bound from \cref{thm:normed-ledger}, or the centered constant
from \cref{thm:ledger} for closed gauge-length paths.
\end{corollary}
\begin{proof}
Transitivity of the two inequalities gives
$B_\alpha\le C_\alpha\len(\gamma)$. Division by the positive
constant gives the pathwise estimate, and taking the infimum over all
feasible paths gives the threshold bound. For balanced support
certificates, $B_\alpha$ is the appropriate weighted sum of the body's
support values. If a proposed certificate has $C_\alpha=0$ and
$B_\alpha>0$, the ledger makes that contact case infeasible rather
than authorizing division by zero.
\end{proof}

\subsection{Exact escape from centered support circuits}\label{sec:circuit-radius}
The centered ledger gives more than a lower-bound certificate. When the
contact points are variables, its sharp constant eliminates those points
altogether. For a nonempty compact set $S$ put
\[
 \operatorname{rad}_T(S)=\min_{z\in\R^d}\max_{x\in S}\mu_T(x-z).
\]
Thus $\operatorname{rad}_T(S)$ is the least factor of a translate of $T$
covering $S$. The orientation of $T$ is fixed.

Compared with the Minkowski-billiard characterization, the task here
is elimination: the non-translatable polygon and at-most-$d+1$-vertex
reduction are already known \cite{BezdekBezdek,AkopyanEtAl,KruppRudolf}.
We optimize the contact vertices in their prescribed temporal order
and retain only the balanced support vectors and a circumradius.
This is an elementary representation of the established escape/capacity
quantity, not a new symplectic invariant.

\begin{theorem}[Exact circuit--circumradius formula]\label{thm:circuit-radius}
Let $K,T$ be convex bodies with $0\in\Int T$. For each ordered list
$a=(a_1,\ldots,a_m)$ satisfying
\[
 2\le m\le d+1,\qquad \sum_{i=1}^m a_i=0,
 \qquad \sum_{i=1}^m|a_i|=1,
\]
define
\[
 R_i=-\sum_{j=1}^i a_j\quad(1\le i\le m),\qquad
 B_K(a)=\sum_i h_K(a_i),\qquad
 C_T(a)=\operatorname{rad}_T\{R_1,\ldots,R_m\}.
\]
Then $B_K(a),C_T(a)>0$ and
\begin{equation}\label{eq:circuit-radius}
 E_T^{\rm tr}(K)=\min_a\frac{B_K(a)}{C_T(a)}.
\end{equation}
The minimum is attained. In particular, $K$ is a closed universal
unit $T$-loop cover if and only if
\begin{equation}\label{eq:support-only-cover}
 \sum_i h_K(a_i)\ge C_T(a)\qquad\hbox{for every such ordered list }a.
\end{equation}
\end{theorem}
\begin{proof}
Translate $K$ to contain $rB_2^d$ and choose $M>0$ with
$T\subset MB_2^d$. Balance makes $B_K$ unchanged by this translation,
and $B_K(a)\ge r$. Consecutive states, including the cyclic pair
$R_m=0,R_1$, differ by the $a_i$. Since some $|a_i|\ge1/m$,
any translate of $cT$ containing all states has
$2cM\ge1/m$. Hence $C_T(a)\ge1/(2mM)>0$.
The parameter sets are compact, and the circumradius is continuous by
\cref{cor:fit}; thus the minimum on the right exists.

Let a non-fitting loop be given. By \cref{thm:balance}, select at
most $d+1$ support contacts and balanced weights that certify its
non-fitting. Put these contacts in their actual cyclic temporal order,
breaking ties arbitrarily, and normalize the
associated vectors so that $\sum|a_i|=1$. If their positions are
$q_0,\ldots,q_{m-1}$, then
\[
 \sum_{i=1}^m a_i\cdot q_{i-1}\ge B_K(a).
\]
The centered ledger in \cref{thm:ledger}, applied to the intervening
increments and the closing increment, bounds this sum by
$C_T(a)\ell_T(q)$. This proves the lower bound in
\eqref{eq:circuit-radius}.

Conversely, fix $a$. The optimality statement
\eqref{eq:ledger-optimal} gives increments $v_1,\ldots,v_m$ with
\[
 \sum_i v_i=0,\qquad \sum_i h_T(v_i)=1,
 \qquad \sum_i R_i\cdot v_i=C_T(a).
\]
The length budget is saturated because $C_T(a)>0$ and otherwise a
positive scaling would increase the objective. Define $q_0=0$ and
$q_i=\sum_{j\le i}v_j$, so $q_m=q_0$. Telescoping gives
\[
 \sum_{i=1}^m a_i\cdot q_{i-1}
 =\sum_{i=1}^m R_i\cdot v_i=C_T(a).
\]
Scale this polygon by $B_K(a)/C_T(a)$. The scaled support sum is
$B_K(a)$, so a strict translate into $K$ would imply
\[
 B_K(a)=\sum_i a_i\cdot(q_{i-1}+z)
       <\sum_i h_K(a_i)=B_K(a),
\]
a contradiction. Its length is the claimed ratio. The selected
vertices need not be prescribed maximizers of the individual supports:
the displayed inequality alone certifies non-fitting. This observation
is why eliminating the contact variables is exact here.

Finally, endpoint coverage and strictification give the equivalence
between universal unit coverage and $E_T^{\rm tr}(K)\ge1$.
Apply \eqref{eq:circuit-radius} to obtain
\eqref{eq:support-only-cover}.
\end{proof}

Haim-Kislev's combinatorial formula \cite[Theorem~1.1]{HaimKislevPolytopes}
optimizes ordered facet weights and symplectic pairings for general
polytopes in $\R^{2d}$. Specialized to a Lagrangian product, it
computes the same capacity as \eqref{eq:circuit-radius}, by the
published identity used in \cref{sec:symplectic}. Here the product
structure and the short-polygon reduction leave balanced circuits
in $\R^d$ and $d+1$ circumradius variables. This is a useful linear
programming representation, without a claim that the underlying
finite combinatorial capacity formula is new.

\begin{corollary}[Facet circuits and small circumradius programs]\label{cor:circuit-lp}
Let $K(b)=\{x:n_i\cdot x\le b_i,\ 1\le i\le N\}$ be a convex
body, with all $n_i\ne0$. Choose one nonzero generator $\lambda^I$
of every extreme ray of $\{\lambda\ge0:\sum_i\lambda_i n_i=0\}$,
and write $I=\supp\lambda^I$. For each ordering $\sigma$ of $I$,
form the states in \cref{thm:circuit-radius} from
$a_j=\lambda^I_{\sigma(j)}n_{\sigma(j)}$, without normalization, and put
\[
 c_I(T)=\max_\sigma C_T(a).
\]
Then
\begin{equation}\label{eq:facet-circumradius}
 E_T^{\rm tr}(K(b))=\min_I\frac{\lambda^I\cdot b}{c_I(T)}.
\end{equation}
If $T=\{y:f_l\cdot y\le d_l,\ 1\le l\le M\}$ with
$d_l>0$, each $C_T(a)$ is the minimum of the linear program
\begin{equation}\label{eq:small-radius-lp}
 \min_{z\in\R^d,\,c\ge0}c
 \quad\text{subject to}\quad
 f_l\cdot(R_j-z)\le c d_l\quad(1\le j\le |I|,\ 1\le l\le M).
\end{equation}
It has only $d+1$ real variables. If the facet data are rational,
choose rational generators $\lambda^I$. Then all displayed constants
and the escape value are rational and admit rational certificates.
\end{corollary}
\begin{proof}
Use \cref{prop:poly-circuit} to choose an extreme-ray certificate
for each non-fitting loop. Apply both constructions in the proof of
\cref{thm:circuit-radius} with those facet coefficients. Every
$\lambda^I\cdot b$ is positive, since an interior point of $K(b)$
gives strictly positive slacks and balance cancels that point.
For fixed $I$, minimizing its positive ratio over orders is the same
as maximizing its denominator. This proves \eqref{eq:facet-circumradius}.
The inequalities in \eqref{eq:small-radius-lp} express exactly
$R_j\in z+cT$. The optimum is attained and positive. Rational
linear-programming duality gives the stated certificates.
\end{proof}

Reciprocity is already implicit in the established capacity and
Minkowski-billiard dualities \cite{RudolfCapacity,KruppRudolf}.
The next argument derives it directly from the centered ledger.

\begin{corollary}[Elementary reciprocity of translative escape]\label{cor:escape-reciprocity}
For convex bodies $K,T$ containing the origin in their interiors,
\[
 E_T^{\rm tr}(K)=E_K^{\rm tr}(T).
\]
\end{corollary}
\begin{proof}
Traverse the state polygon $\{R_i\}$ in reverse order. Its edge
vectors are the $a_i$, up to cyclic reindexing, and its $K$-length
is $B_K(a)$. Conversely, every nonconstant closed polygon with at
most $d+1$ vertices is such a state polygon, up to translation,
positive scaling, and choice of traversal. Thus
\eqref{eq:circuit-radius} is the minimum of
$\ell_K(P)/\operatorname{rad}_T(P)$ over these polygons.
Scaling by the reciprocal circumradius makes a polygon non-fitting
at fitting scale one. Conversely, a shortest non-fitting polygon in
$T$ has fitting scale one, since a larger scale could be reduced.
\Cref{thm:minkowski} now identifies this minimum with
$E_K^{\rm tr}(T)$.
\end{proof}

In the plane, the following complete-family statement is
\cite[Theorem~3.3]{BMP}. That paper also observes that the approach
extends to higher dimensions. We provide the full dimensional proof
and retain the circuit counts needed for the explicit vertex bound
and the outer optimization below.

\begin{theorem}[A finite complete test family for a polytopal gauge]\label{thm:complete-test-family}
Let $T=\{y:n_i\cdot y\le b_i,\ 1\le i\le N\}$ be a
polytopal gauge with zero in its interior and $n_i\ne0$ for all $i$.
Consequently $b_i>0$. For every extreme positive facet circuit $\lambda^I$ and each cyclic ordering $\sigma$ of its
support, put
\[
 B_I=\sum_i\lambda_i^I b_i,
 \qquad
 Q_{I,\sigma}=B_I^{-1}
 \left\{-\sum_{j=1}^k\lambda^I_{\sigma(j)}n_{\sigma(j)}:
                          1\le k\le |I|\right\}.
\]
Use the polygonal loop through these vertices in reverse cyclic order
as the corresponding test object. Each such loop has $T$-length at most one, and exactly
one if all displayed facet offsets are supporting offsets.
A convex body $K$ covers every closed unit $T$-loop by translation
if and only if it covers every $Q_{I,\sigma}$ by translation; by
convexity, covering its vertices is equivalent to covering its loop.
Only $(|I|-1)!$ cyclic orders per circuit are needed, with reversals
retained. Every test polygon has at most $d+1$ vertices.
\end{theorem}
\begin{proof}
The reversed state polygon has edge vectors
$\lambda_i^I n_i/B_I$, so its length is
\[
 \frac{\sum_i\lambda_i^I h_T(n_i)}{B_I}\le1.
\]
Necessity follows. For sufficiency, translate $K$ so that zero lies
in its interior. By reciprocity and the facet-circuit formula applied
to the \emph{forest} $T$ with \emph{gauge} $K$,
\[
 E_T^{\rm tr}(K)=E_K^{\rm tr}(T)
   =\min_{I,\sigma}
       \frac{B_I}{\operatorname{rad}_K(B_IQ_{I,\sigma})}.
\]
Thus the escape value is at least one exactly when every displayed
test set has $K$-circumradius at most one, which is its closed
translative containment in $K$. Endpoint strictification finishes
the equivalence. A cyclic change of starting index translates the
state polygon, so only cyclic orders must be retained. Reversing
an order can give a different trace up to translation and is not
discarded for an asymmetric gauge.
\end{proof}

For a regular dodecagonal gauge, \cref{prop:circuit-count} gives six
two-normal circuits and forty three-normal circuits. Hence this
theorem produces a complete list of at most $6+2\cdot40=86$ test
loops, with at most $2\cdot6+6\cdot40=252$ vertices in total.
These are exact tests for that polytopal gauge. Approximation of the
Euclidean gauge remains governed by \cref{cor:gauge-rate}.

The reciprocity is consistent with the duality of Minkowski billiards
\cite{AkopyanEtAl,RudolfCapacity}. The contribution of
\eqref{eq:facet-circumradius} to the present framework is its explicit
elimination of path vertices by the sharp centered ledger. General
combinatorial capacity formulas for symplectic polytopes already exist
\cite{HaimKislevPolytopes}; no claim that finite capacity descriptions
originate here is intended.

\subsection{Orientation discretization and adaptive recovery}\label{sec:rotation}
The norm in this section is Euclidean. Let $S\subset H$ be a nonempty
finite orientation set. Define $e_S$ by the constraints
$s_K(\Gamma_\gamma,R)\ge1$ for $R\in S$ only.

\begin{lemma}[Centered orientation estimate]\label{lem:orientation-lip}
If a trace has a translate in $bB_2^d$ and $rB_2^d\subset K$, then
\begin{equation}\label{eq:orientation-lip}
 |s_K(A,R)-s_K(A,Q)|\le\frac b r\|R-Q\|_\op.
\end{equation}
For the trace of a curve of length $\ell$, one may take $b=\ell/2$.
\end{lemma}
\begin{proof}
Translate $A$ into $bB_2^d$. Its support function is $b$-Lipschitz on
unit directions. Each numerator in \eqref{eq:fit-factor} therefore
changes by at most $b\|R-Q\|_\op$, and each denominator is at least
$r$. Taking maxima preserves this estimate. For an open curve, its
arclength midpoint is within distance $\ell/2$ of every point; for a
closed curve any point has the same property by using the shorter
cyclic subarc. Translation invariance permits these choices of center.
\end{proof}

\begin{theorem}[Explicit rotational sampling rate]\label{thm:rotation-rate}
Suppose $rB_2^d\subset K$ and $S$ is an $\eta$-net of $H$ in operator
norm. The sampled infimum $e_S$ is attained and
\begin{equation}\label{eq:rot-rate}
 e_S\le E_H(K)\le
 \frac{e_S}{1-e_S\eta/(2r)}
 \qquad(e_S\eta<2r).
\end{equation}
In particular, if $B\ge E_H(K)$ and $B\eta<2r$, then
\begin{equation}\label{eq:rot-gap}
 0\le E_H(K)-e_S\le\frac{B^2\eta}{2r-B\eta}.
\end{equation}
One may take $B=\diam K$. Thus
\begin{equation}\label{eq:mesh-choice}
 \eta\le\frac{2r\eps}{B(B+\eps)}
\end{equation}
is sufficient for additive error at most $\eps$.
\end{theorem}
\begin{proof}
If $\eta=0$, the finite sample equals $H$ and the conclusion is exact.
Assume $\eta>0$. The finite sampled constraints are closed under uniform convergence:
support functions of uniformly convergent traces converge uniformly.
The arclength compactness argument gives a minimizing curve of length
$e=e_S$. The full problem has additional constraints, so $e\le E_H(K)$.
By \cref{lem:orientation-lip}, every orientation satisfies
\[
 s_K(\Gamma_\gamma,R)\ge1-\frac{e\eta}{2r}.
\]
Scale the curve by the reciprocal of the positive right-hand side.
Fit factors are homogeneous, so the scaled curve is continuum feasible
and has length $e/(1-e\eta/(2r))$. This proves \eqref{eq:rot-rate}.
Subtract $e$ and use that $x\mapsto x^2\eta/(2r-x\eta)$ is increasing
on $[0,2r/\eta)$ to obtain \eqref{eq:rot-gap}.
A diameter-length segment cannot have both endpoints in $\Int K$;
otherwise interior balls allow a longer segment in $K$. It is therefore
an escape path. Solving the last error bound for $\eta$ gives
\eqref{eq:mesh-choice}.
\end{proof}

For $H=\SO(2)$, $M$ equally spaced rotations have
$\eta=2\sin(\pi/(2M))$. For $\OO(2)$ both components must be sampled.
These statements do not require a path to visit supports in angular
order. That order is a discrete choice, not a consequence of using a
support function.

\begin{proposition}[Finite optimization and certified incumbents]\label{prop:finite-rot}
For $m$ sampled orientations, a minimizing sampled curve can be chosen
polygonal with at most $m(d+1)$ vertices. If $K$ is a polytope, $e_S$
is the minimum of finitely many second-order cone program values.
If a certified sampled lower bound is $l$ and a sampled feasible curve
has length at most $u$, then
\begin{equation}\label{eq:solver-gap}
 l\le E_H(K)\le \frac{u}{1-u\eta/(2r)}
 \qquad(u\eta<2r).
\end{equation}
\end{proposition}
\begin{proof}
For each orientation choose the at-most-$d+1$ support contact points
given by \cref{thm:helly}, applied to that rotated trace. Take their
union in traversal order and shortcut. Every orientation retains its
own certificate. The prefix before the first selected point can be
discarded, and translation makes this point the root zero.

For a polytope, enumerate one extreme facet circuit per sampled
orientation and assign its support directions to the ordered polygon
vertices $p_1,\ldots,p_N$, $N=m(d+1)$. Repeated vertices permit padding.
For each assignment impose
\[
 \sum_i\lambda_i a_i\cdot R p_{\sigma_R(i)}
       \ge\sum_i\lambda_i b_i.
\]
These are linear inequalities. Minimize
$\sum_{j=1}^{N-1}L_j$ subject to
$\|p_{j+1}-p_j\|_2\le L_j$ and $p_1=0$.
Every feasible branch is a valid sampled obstruction, and every
sampled minimizing polygon appears in some branch. This proves exactness.
The final inequality follows by applying the scaling repair to the
sampled feasible incumbent; the repair factor is increasing in its
length.
\end{proof}

\begin{corollary}[An effective finite-tolerance procedure]\label{cor:eps}
For a rational polytope $K$ with nonempty interior and a rotation group
admitting effectively constructible real-algebraic finite nets, the
open-path threshold is computable
to any positive additive tolerance, using exact real-algebraic
optimization of the finitely many branches and \eqref{eq:solver-gap}.
\end{corollary}
\begin{proof}
Choose computable bounds $r>0$ and $B\ge\diam K$, then a sufficiently
fine net using \eqref{eq:mesh-choice} with half the requested tolerance,
leaving a positive allowance for finite optimization error. In fixed dimension, the orthogonal
group is a compact real-algebraic set and admits effective algebraic
nets; in two dimensions the explicit angle formula suffices. Each
branch is a semialgebraic optimization problem, so quantifier
elimination provides enclosing intervals for its value. Refine those
intervals until the lower bound and the repaired upper bound differ
by at most the requested tolerance.
\end{proof}
This is an effectivity result, not a polynomial-time claim. The branch
enumeration can be very large, and numerical termination tolerances
without validated bounds are not certificates.

\begin{proposition}[Adaptive orientation exchange]\label{prop:adaptive}
Suppose exact sampled minimizers and exact worst-orientation separation
are available. At each iteration add an orientation minimizing
$s_K(\Gamma_\gamma,R)$ whenever its value is below one. Then sampled
values are nondecreasing and converge to $E_H(K)$. Unless a continuum
optimizer is found earlier, their maximal deficits
\[
 a_m=\max\{0,1-\min_{R\in H}s_K(\Gamma_{\gamma_m},R)\}
\]
tend to zero. Whenever $a_m<1$, the repaired curve
$\gamma_m/(1-a_m)$ supplies an upper bound, and the lower--upper gap
tends to zero.
\end{proposition}
\begin{proof}
Sampled lengths are at most $B=\diam K$. If a new orientation has deficit
at least $a>0$, \cref{lem:orientation-lip} puts it at distance at least
$2ra/B$ from every previously sampled orientation. Infinitely many such
insertions contradict compactness of $H$. Therefore $a_m\to0$.
Homogeneity repairs the exact worst-orientation deficit, giving
\[
 e_m\le E_H(K)\le e_m/(1-a_m).
\]
The boundedness of $e_m$ proves convergence and vanishing gap.
\end{proof}

\subsection{Exact polytopal linear programs}\label{sec:lp}
\begin{theorem}[Exact finite LP formula]\label{thm:lp}
Let $K=\{x:a_i\cdot x\le b_i,\ 1\le i\le m\}$ and
$T=\conv\{t_1,\ldots,t_s\}$ be polytopes, with every $a_i\ne0$
and $0\in\Int K\cap\Int T$.
For every $2\le k\le d+1$, extreme ray $\lambda$ of $C_K$, and map
$\sigma:\supp\lambda\to\{1,\ldots,k\}$, consider
\begin{equation}\label{eq:lp}
\begin{aligned}
 \min\quad &\sum_{j=1}^kL_j,\\
 \text{subject to}\quad
 &p_1=0,\quad p_{k+1}=p_1,\\
 &L_j\ge t_\ell\cdot(p_{j+1}-p_j)
            &&(1\le j\le k,\ 1\le\ell\le s),\\
 &\sum_{i\in\supp\lambda}\lambda_i a_i\cdot p_{\sigma(i)}
       \ge\sum_i\lambda_i b_i.
\end{aligned}
\end{equation}
Then $E_T^{\rm tr}(K)$ equals the minimum of these finitely many LP
values, with value $+\infty$ for an infeasible branch. For rational data
the threshold is rational and has a finite exact certificate.
\end{theorem}
\begin{proof}
The last inequality implies \eqref{eq:poly-nonfit}, since the actual
support is at least the chosen point evaluation. Thus each feasible
branch is a valid obstruction. The epigraph inequalities minimize
exactly the polygon's $T$-length. Conversely, a minimizing polygon
from \cref{thm:minkowski} has a nonnegative-violation extreme circuit;
assign every circuit normal to a support-attaining vertex. Its data
are feasible in the corresponding branch, proving equality.

For rational data the extreme rays can be rationally normalized.
Every feasible branch has a positive attained finite optimum:
bounded objective bounds all anchored vertices because
$h_T(v)\ge c|v|$. Rational linear programming gives rational optimal
solutions and dual certificates. Infeasible branches have rational
Farkas certificates. To certify a lower bound for the minimum, every
branch must be certified; one good branch alone is insufficient.
\end{proof}

\subsection{Diameter certificates and explicit meshes}\label{sec:diameter-cert}
\begin{proposition}[Diameter thickening certificate]\label{prop:diameter-hit}
Let $C$ be a nonempty compact metric space and let nonempty compact targets in a
finite-dimensional normed space of positive dimension satisfy
$\dH(N_c,N_{c'})\le L_Nd_C(c,c')$. Let $\Delta$ be the least diameter
of a compact set meeting every target and $\Delta_\eta$ the optimum on
a finite $\eta$-net. Both optima are attained, and
\begin{equation}\label{eq:diameter-net}
 \Delta_\eta\le\Delta\le\Delta_\eta+2L_N\eta.
\end{equation}
The statement remains valid with a prescribed anchor required in every
hitting set. For nested finite nets of meshes $\eta_m\to0$,
$\Delta_m\uparrow\Delta$; in particular $2L_N\eta_m\le\eps$ is an
explicit absolute-error certificate.
\end{proposition}
\begin{proof}
Hausdorff continuity on compact $C$ bounds the union of the targets,
so a sufficiently large closed ball is a feasible compact hitting set.
The finite problem is the minimum diameter of one selected point in
each sampled target, together with the anchor if prescribed. Diameter
is continuous on that finite compact product, so the minimum is attained.

For any nonempty compact $S$ and $a\ge0$, in the given norm,
\begin{equation}\label{eq:diameter-thickening-exact}
 \diam(S+aB)=\diam(S)+2a.
\end{equation}
The upper bound follows from the triangle inequality. If $x,y\in S$
realize positive diameter, let $v=(x-y)/\|x-y\|$. The two points
$x+av,y-av$ in $S+aB$ have distance $\diam(S)+2a$. If $S$ is a
singleton, two opposite radius-$a$ points prove the same formula.
Thicken a finite minimizer by $L_N\eta B$. For every $c\in C$, a
sampled $c'$ has $d_C(c,c')\le\eta$, and a selected point in $N_{c'}$
is within $L_N\eta$ of $N_c$. The thickened set meets $N_c$ and
retains the anchor. Equation \eqref{eq:diameter-thickening-exact}
proves the upper bound in \eqref{eq:diameter-net}; the lower bound
comes from omitting constraints.

For continuum attainment, take a sequence of feasible compact sets with
bounded diameter. Each meets a fixed compact target $N_{c_0}$, so all
these sets lie in a common compact ball. Selection in the hyperspace of
nonempty compact subsets gives a Hausdorff limit. For any fixed target,
selected intersection points have a convergent subsequence in that ball;
closedness of the target preserves the hit. The anchor is preserved as
well. Diameter is Hausdorff continuous, proving attainment. Finally,
nested samples give monotonicity, and \eqref{eq:diameter-net} gives
$0\le\Delta-\Delta_m\le2L_N\eta_m\to0$. If $L_N=0$, any nonempty
sample already has the exact value.
\end{proof}

\begin{corollary}[An explicit Euclidean configuration mesh]\label{cor:diameter-mesh}
Let $K\subset RB_2^d$ be a convex body and $H\le\OO(d)$ be compact.
Sample $x\in K$ with Euclidean mesh $\eta_x$ and $Q\in H$ with
operator-norm mesh $\eta_Q$. Let $D_{\rm sam}$ be the minimum diameter
of a compact set containing zero and meeting every sampled target
$Q^{-1}(\partial K-x)$. Then
\begin{equation}\label{eq:diameter-product-mesh}
 D_{\rm sam}\le D_H(K)\le
 D_{\rm sam}+2\eta_x+4R\eta_Q.
\end{equation}
Nested product nets with both meshes tending to zero converge monotonically.
\end{corollary}
\begin{proof}
For any $y\in\partial K$ and configurations $(x,Q),(x',Q')$,
\[
 \bigl|Q^{-1}(y-x)-Q'^{-1}(y-x')\bigr|
 \le |x-x'|+2R\|Q-Q'\|_{\op}.
\]
Matching the same $y$ in both directions proves this Hausdorff bound
for the targets. An anchored hitting set can be convexified without
increasing diameter by \cref{lem:convex-hull}. A rooted convex set is
connected: the segment from zero to a point failing interior containment
meets the boundary, by \cref{prop:motion}. Thus the full anchored hitting
optimum equals $D_H(K)$. Applying the thickening proof with radius
$\eta_x+2R\eta_Q$ gives \eqref{eq:diameter-product-mesh}. The same
constraint inclusion proves monotonicity.
\end{proof}

Anchoring supplies the premise for the boundary-contact argument.
Thickening controls diameter, but does not supply a comparable length
bound for a path tracing the thickened set.

\section{Optimization of Minkowski worms and universal covers}\label{sec:cover-optimization}
\subsection{Body and gauge stability under Hausdorff and affine changes}\label{sec:stability}
\begin{theorem}[One stability estimate for all homogeneous thresholds]
\label{thm:body-stability}
Let $J$ be a finite nonnegative, translation-invariant,
inclusion-monotone, degree-one homogeneous functional on convex bodies. If
$rB_2^d\subset K$, $\dH(K,L)\le\delta<r$, then
\begin{equation}\label{eq:body-stability}
 (1-\delta/r)J(K)\le J(L)\le(1+\delta/r)J(K).
\end{equation}
This applies to $E_H$, $D_H$, and the translative Minkowski threshold
defined in \cref{sec:minkowski}.
\end{theorem}
\begin{proof}
Support functions satisfy $|h_K(u)-h_L(u)|\le\delta$ for unit $u$ and
$h_K(u)\ge r$. Therefore
\[
 (1-\delta/r)h_K(u)\le h_L(u)\le(1+\delta/r)h_K(u),
\]
which is equivalent to
$(1-\delta/r)K\subset L\subset(1+\delta/r)K$.
Apply monotonicity and homogeneity. For escape thresholds, monotonicity
holds because enlarging the cover shrinks its non-fitting family, and
homogeneity follows by scaling objects and cover together.
\end{proof}

\begin{corollary}[Quotient continuity and polytope density]\label{cor:density}
If $J(K)>0$ and the hypotheses of \cref{thm:body-stability} hold, then
\begin{equation}\label{eq:quotient-perturbation}
 \left(\frac{1-\delta/r}{1+\delta/r}\right)^d
 \frac{\vol(K)}{J(K)^d}
 \le\frac{\vol(L)}{J(L)^d}\le
 \left(\frac{1+\delta/r}{1-\delta/r}\right)^d
 \frac{\vol(K)}{J(K)^d}.
\end{equation}
Consequently all corresponding quotient infima are unchanged if the
search is restricted to rational polytopes with nonempty interior.
\end{corollary}
\begin{proof}
The inclusions in the preceding proof bound volume by the $d$th powers
of $1\pm\delta/r$. Combine these bounds with \eqref{eq:body-stability}.
Rational polytopes are Hausdorff dense among convex bodies, and the
limiting threshold is positive, so their quotients converge to that of
each body.
\end{proof}

\begin{theorem}[Gauge comparison and affine covariance]\label{thm:gauge}
If $\alpha T\subset S\subset\beta T$, $0<\alpha\le\beta$, then
\begin{equation}\label{eq:gauge-comparison}
 \alpha E_T^{\rm tr}(K)\le E_S^{\rm tr}(K)\le\beta E_T^{\rm tr}(K),
 \qquad
 \beta^{-d}W_T\le W_S\le\alpha^{-d}W_T.
\end{equation}
For $A\in\GL(d)$,
\begin{equation}\label{eq:affine}
 E_{AT}^{\rm tr}(K)=E_T^{\rm tr}(A^TK),\qquad
 W_{AT}=|\det A|^{-1}W_T.
\end{equation}
Translation of $T$ leaves closed-loop lengths and $W_T$ unchanged,
provided the translated gauge is still represented with the origin
in its interior. Hence
$\Psi_d(T)=d!\vol(T)W_T$ is affine invariant, including translations.
\end{theorem}
\begin{proof}
Gauge inclusions imply $\alpha h_T\le h_S\le\beta h_T$.
Integrate on the common non-fitting loop class and take infima.
Apply \eqref{eq:minkowski-quotient} for the volume bounds.
The identity $h_{AT}(v)=h_T(A^Tv)$ gives
$\ell_{AT}(q)=\ell_T(A^Tq)$, while translation fitting in $K$
corresponds bijectively to fitting in $A^TK$. This proves the first
identity in \eqref{eq:affine}; changing variables in the quotient
proves the second. Finally
\[
 \ell_{T+b}(q)-\ell_T(q)
   =b\cdot\int_0^1\dot q(t)\,dt=0
\]
for a closed loop. The determinant factors in $\Psi_d$ cancel.
\end{proof}

\begin{corollary}[Explicit Hausdorff modulus in the gauge]\label{cor:gauge-hausdorff}
If $rB_2^d\subset T$ and $\dH(T,S)\le\delta<r$, then
\begin{align}
 (1-\delta/r)E_T^{\rm tr}(K)&\le E_S^{\rm tr}(K)
                  \le(1+\delta/r)E_T^{\rm tr}(K),\label{eq:gauge-hausdorff}\\
 (1+\delta/r)^{-d}W_T&\le W_S\le(1-\delta/r)^{-d}W_T.
                  \label{eq:gauge-w-hausdorff}
\end{align}
\end{corollary}
\begin{proof}
For unit $u$, $|h_S(u)-h_T(u)|\le\delta$ and $h_T(u)\ge r$.
Thus $(1-\delta/r)T\subset S\subset(1+\delta/r)T$. Apply
\cref{thm:gauge}. The estimates hold for every $K$ with the same
constants, making the gauge modulus uniform in the cover body.
\end{proof}

\begin{corollary}[Joint continuity and affine stability]\label{cor:psi}
The maps $(K,T)\mapsto E_T^{\rm tr}(K)$ and $T\mapsto W_T$ are
Hausdorff continuous at full-dimensional bodies. Under the inclusions
of \cref{thm:gauge},
\begin{equation}\label{eq:psi-stability}
 (\alpha/\beta)^d\Psi_d(T)\le\Psi_d(S)
       \le(\beta/\alpha)^d\Psi_d(T).
\end{equation}
\end{corollary}
\begin{proof}
The support sandwich used in \cref{thm:body-stability} also applies
to the gauge. Thus Hausdorff perturbations give inclusions with
$\alpha,\beta\to1$. Use \cref{thm:gauge} and then
\cref{thm:body-stability} for joint continuity. Multiplying the
volume bounds $\alpha^d\vol(T)\le\vol(S)\le\beta^d\vol(T)$ by
the corresponding bounds for $W$ proves \eqref{eq:psi-stability}.
\end{proof}
Here $0<W_T<\infty$: a sufficiently small loop through
$0,\eps e_1,\ldots,\eps e_d,0$ has $T$-length at most one and
forces every cover to contain a translated simplex of positive volume.
Conversely, $h_T(v)\ge c|v|$ bounds the Euclidean length of every unit
$T$-loop by $1/c$, so a sufficiently large fixed ball covers all of them.
Thus $\Psi_d(T)>0$, and in particular
\[
 |\log\Psi_d(S)-\log\Psi_d(T)|\le d\log(\beta/\alpha).
\]
Affine invariance allows the inclusions to be optimized over affine
normalizations; this gives the corresponding Banach--Mazur stability
estimate without choosing a preferred coordinate system.

The inclusion comparison can be sharpened by exploiting the cancellation
of translations in the length of a closed curve.

\begin{theorem}[Optimal translation-invariant gauge comparison]\label{thm:translated-gauges}
For convex gauges $S,T$ containing zero in their interiors, set
$D(S,T)=\operatorname{rad}_T(S)$. Then
\begin{equation}\label{eq:optimal-gauge-ratio}
 D(S,T)=\max_{\substack{q\text{ nonconstant, rectifiable}\\q\text{ closed}}}
              \frac{\ell_S(q)}{\ell_T(q)}.
\end{equation}
A maximizing polygon has at most $d+1$ edges. The constant is unchanged
by independent translations of $S$ and $T$ within the admissible gauge
class, and
\begin{align}
 D(T,S)^{-1}E_T^{\rm tr}(K)
   &\le E_S^{\rm tr}(K)\le D(S,T)E_T^{\rm tr}(K),
       \label{eq:translated-escape}\\
 D(S,T)^{-d}W_T&\le W_S\le D(T,S)^dW_T.
       \label{eq:translated-worm}
\end{align}
The two length-comparison constants are optimal over all closed curves;
optimality is not asserted separately for every fixed escape body.
\end{theorem}
\begin{proof}
A placement $S+z\subset D(S,T)T$ gives
$h_S(v)+z\cdot v\le D(S,T)h_T(v)$. Integrate and use
$\int\dot q=0$ to prove the upper bound in
\eqref{eq:optimal-gauge-ratio}. The balanced fitting formula
\eqref{eq:fit-factor} gives a maximizing collection
$(\lambda_i,u_i)_{i=1}^m$, $m\le d+1$, with
\[
 D(S,T)=\frac{\sum_i\lambda_i h_S(u_i)}
                  {\sum_i\lambda_i h_T(u_i)}.
\]
The vectors $v_i=\lambda_i u_i$ close to a polygon and realize
exactly this ratio of lengths. The denominator is positive because
the balanced weights sum to one and $T$ contains an inball.
Translations add linear terms that cancel both in the balanced ratios
and in closed-curve lengths. Apply the resulting length inequalities
to the common class of non-fitting loops, and then use
\eqref{eq:minkowski-quotient}, to obtain the two displayed comparisons.
\end{proof}

For instance, translating a gauge within the class of bodies containing
zero may substantially worsen its origin-centered inclusion factors.
Nevertheless $D(T+z,T)=D(T,T+z)=1$, and all closed-loop escape and
worm values remain exactly equal. This removes a coordinate artifact
from the earlier multiplicative bounds.

\begin{proposition}[Reversal, asymmetry, and central symmetrization]\label{prop:gauge-symmetrization}
Let $\sigma(T)=D(-T,T)$ and $T_s=(T-T)/2$. Then
$1\le\sigma(T)\le d$, with equality $\sigma(T)=1$ precisely for
centrally symmetric bodies up to translation. Moreover
\begin{align}
 E_{-T}^{\rm tr}(K)&=E_T^{\rm tr}(K),\qquad W_{-T}=W_T,
       \label{eq:gauge-reversal}\\
 E_T^{\rm tr}(K)&\le E_{T_s}^{\rm tr}(K)
                  \le\frac{1+\sigma(T)}2 E_T^{\rm tr}(K),
       \label{eq:gauge-sym-escape}\\
 W_{T_s}&\le W_T\le
                  \left(\frac{1+\sigma(T)}2\right)^dW_{T_s}.
       \label{eq:gauge-sym-worm}
\end{align}
The length-comparison factor $(1+\sigma(T))/2$ is optimal.
\end{proposition}
\begin{proof}
If $q^-(t)=q(1-t)$, then $\ell_T(q^-)=\ell_{-T}(q)$, while
the traces and their fitting properties coincide. This proves
\eqref{eq:gauge-reversal}. Also
\[
 \ell_{T_s}(q)=\frac{\ell_T(q)+\ell_T(q^-)}2.
\]
Both terms on the right are at least $E_T^{\rm tr}(K)$ for a
non-fitting loop. Conversely, choose a shortest $T$-loop and use
\cref{thm:translated-gauges} to bound the length of its reversal by
$\sigma(T)$ times its original length. This proves
\eqref{eq:gauge-sym-escape}, and the quotient gives
\eqref{eq:gauge-sym-worm}. The same averaging identity and the
attained maximum in \eqref{eq:optimal-gauge-ratio} show that
$D(T_s,T)=(1+\sigma(T))/2$.

Volume monotonicity gives $\sigma(T)\ge1$. At equality, an attained
translate of $-T$ in $T$ has equal volume and hence equals $T$,
which is central symmetry. To prove the upper bound, inscribe a
maximum-volume simplex $\Delta\subset T$ and let $z$ be its
centroid. Replacing one vertex by $x\in T$ shows that each barycentric
coordinate $\beta_i(x)$ has absolute value at most one. In particular
$\beta_i(x)\le1$. The barycentric coordinates of
$z-(x-z)/d$ are $(1-\beta_i(x))/d$, so this point belongs to
$\Delta$. Consequently
$-(T-z)\subset d(\Delta-z)\subset d(T-z)$, proving
$\sigma(T)\le d$.
\end{proof}

The orientation distinction is visible in exact rational data. For
\[
 T=\conv\{(1,0),(0,1),(-1,-1)\},\qquad
 q:0\longrightarrow(1,1)\longrightarrow(-1,2)\longrightarrow0,
\]
the three forward edge supports are all one, whereas their negatives
have support two. Thus
\begin{equation}\label{eq:asymmetric-example}
 \ell_T(q)=3,\qquad \ell_T(q^-)=6,\qquad
 \ell_{T_s}(q)=9/2,
 \qquad\sigma(T)=2.
\end{equation}
The last equality follows from the ratio and $\sigma(T)\le d$.
This is why the object must retain its traversal even though the
placement constraint uses only its trace.

\subsection{Attainment, finite test families, and recovery of covers}\label{sec:outer}
The inner escape problem minimizes the size of a non-fitting object. The
outer covering problem minimizes the volume of a body accommodating an
entire object family. The following result supplies the compactness needed
for the latter problem and an explicit finite-family error estimate.

\begin{theorem}[Outer attainment and multiplicative net bound]\label{thm:outer-net}
Let $\mathscr A$ be a compact family of nonempty compact subsets of
$\R^d$, in Hausdorff distance, and let $H\subset\OO(d)$ be compact.
Placements have the form $z+RA$, with $R\in H$. Suppose some
$A_0\in\mathscr A$ has a convex hull containing a ball of radius
$\rho>0$, and some convex body covers $\mathscr A$ with volume at
most $V_+<\infty$. Let $W$ be the infimum of covering volumes.
Then $W>0$ and it is attained.

Let $\mathscr A_\delta\subset\mathscr A$ be a finite Hausdorff
$\delta$-net which includes $A_0$, and let $W_\delta$ be its
minimum covering volume. Then $W_\delta$ is attained and
\begin{equation}\label{eq:outer-net}
 W_\delta\le W\le(1+\delta/\rho)^d W_\delta.
\end{equation}
Consequently a nested sequence of finite subfamilies containing $A_0$
and with dense union has minima increasing to $W$. Normalized minimizing
covers have Hausdorff-convergent subsequences, and every such limit along
vanishing meshes is a minimum cover of $\mathscr A$.
\end{theorem}
\begin{proof}
A placement of $A_0$ in any feasible convex body places a ball of
radius $\rho$ in that body. Translate this ball's center to zero.
If $x\in K$, the cone with apex $x$ and base
$\rho B_2^d\cap x^\perp$ lies in $K$. For $d\ge2$ its volume is
$\kappa_{d-1}\rho^{d-1}|x|/d$, where $\kappa_j$ is the volume
of the $j$-dimensional Euclidean unit ball. Thus
\begin{equation}\label{eq:outer-radius}
 \rho B_2^d\subset K\subset RB_2^d,
 \qquad R=\frac{dV}{\kappa_{d-1}\rho^{d-1}}
 \quad\text{if }\vol(K)\le V.
\end{equation}
For $d=1$ the corresponding boundedness follows from the length of an
interval. A minimizing sequence has a common bound with $V=V_++1$.
Blaschke selection gives a Hausdorff-convergent subsequence of its
normalized bodies. Volume is continuous for convex bodies.

For each fixed $A\in\mathscr A$, choose placements
$z_n+R_nA\subset K_n$. A fixed point of $A$, boundedness of $K_n$,
and compactness of $H$ bound the translations. Passing to a subsequence
for this particular $A$ yields a placement in the limit body. A common
subsequence for all $A$ is unnecessary: feasibility is an intersection
of closed conditions. The limit covers $\mathscr A$, proving attainment.
The same proof applies to a finite subfamily containing $A_0$.
The inball gives $W\ge\kappa_d\rho^d>0$.

Let $K_\delta$ be a finite-family minimizer, normalized to contain
$\rho B_2^d$. For every $A\in\mathscr A$, choose
$A'\in\mathscr A_\delta$ with $\dH(A,A')\le\delta$.
A placement of $A'$ in $K_\delta$ places $A$ in
$K_\delta+\delta B_2^d$, because orthogonal maps preserve distance.
The inclusion
\[
 K_\delta+\delta B_2^d\subset(1+\delta/\rho)K_\delta
\]
proves \eqref{eq:outer-net}. A nested dense sequence is eventually a
$\delta$-net for every $\delta>0$, by compactness. Finally, the
radius bound and the same closedness argument, with nearby sampled
objects, prove the stated recovery of minimizers.
\end{proof}

\begin{proposition}[Stability of a finite covering optimum]\label{prop:finite-stability}
Let $F(A_0,\ldots,A_m)$ denote the minimum convex covering volume
under the motions in \cref{thm:outer-net}. Retain the same witness
$A_0$, whose convex hull contains a radius-$\rho$ ball. If
$\dH(A_i,A_i')\le\delta$ for all $i\ge1$, then
\[
 (1+\delta/\rho)^{-d}F(A_0,A_1,\ldots,A_m)
 \le F(A_0,A_1',\ldots,A_m')
 \le(1+\delta/\rho)^dF(A_0,A_1,\ldots,A_m).
\]
\end{proposition}
\begin{proof}
Normalize a minimizing cover of the first family to contain the ball
provided by its placement of $A_0$. Parallel enlargement by $\delta$
covers the second family and multiplies volume by at most
$(1+\delta/\rho)^d$. Interchange the two families for the other
inequality. Their minima exist by \cref{thm:outer-net}.
\end{proof}

The witness itself can also be perturbed. More precisely, if
$\dH(A_i,A_i')\le\delta<\rho$ for all $0\le i\le m$, and $F,F'$
are the respective optimum volumes, then
\begin{equation}\label{eq:all-witness-stability}
 \left(1+\frac{\delta}{\rho-\delta}\right)^{-d}F
 \le F'\le\left(1+\frac{\delta}{\rho}\right)^dF.
\end{equation}
Indeed convexification does not increase Hausdorff distance. After
centering the original inball, the support inequality
$h_{\conv A_0'}\ge h_{\conv A_0}-\delta\ge\rho-\delta$ shows that
the perturbed witness forces an inball of radius $\rho-\delta$.
Apply the same enlargement argument first from the original family to
the perturbed one, using $\rho$, and then in reverse using
$\rho-\delta$. This proves continuity even when every test object is
allowed to vary.

\begin{corollary}[Additive parallel-body certificate]\label{cor:steiner-error}
Under \cref{thm:outer-net}, for $d\ge2$,
\begin{equation}\label{eq:steiner-error}
 0\le W-W_\delta\le\kappa_d\bigl((R+\delta)^d-R^d\bigr),
 \qquad R=\frac{dV_+}{\kappa_{d-1}\rho^{d-1}}.
\end{equation}
\end{corollary}
\begin{proof}
Normalize the finite-family minimizer as in \eqref{eq:outer-radius}.
Since its volume $W_\delta\le V_+$, it lies in $RB_2^d$. The body
$K_\delta+\delta B_2^d$ covers the continuum family. Steiner's
formula writes the volume of $K+tB_2^d$ as a polynomial in $t$ whose
positive-degree coefficients are nonnegative multiples of the intrinsic
volumes of $K$; these coefficients are monotone under inclusion
\cite{Schneider}. Comparing the coefficients for $K_\delta$ and
$RB_2^d$ gives
\[
 \vol(K_\delta+\delta B_2^d)-\vol(K_\delta)
 \le\vol((R+\delta)B_2^d)-\vol(RB_2^d).
\]
The constant terms cancel before this comparison. The enlarged body is
feasible, so this difference bounds $W-W_\delta$, proving the claim.
\end{proof}

\begin{corollary}[Applicability to the principal Euclidean families]\label{cor:outer-families}
Minimum-volume convex covers exist for open Euclidean curves of length
at most one under translations and any compact orientation group, for
all compact sets of diameter at most one under the same motions, and
for closed curves of $T$-length at most one under translations.
Each problem admits the bounds and recovery statement of
\cref{thm:outer-net}.
\end{corollary}
\begin{proof}
Every Euclidean curve of length at most one admits a rooted
parametrization on $[0,1]$ with Lipschitz constant at most one. These
maps form a compact family by Arzel\`a--Ascoli, and their traces are
exactly the normalized admissible traces. For the Minkowski family,
choose $c>0$ with $h_T(v)\ge c|v|$. Use rooted maps of Lipschitz
constant at most $1/c$ that are closed and have $T$-length at most
one. Closedness and the length constraint persist under uniform limits,
by the polygonal-sum argument in \cref{thm:minkowski}, so this is again
a compact family representing all admissible traces. Taking traces is
continuous from uniform to Hausdorff convergence. Diameter objects containing zero form a closed
family of compact subsets of the unit ball and hence are compact.
Each family contains a sufficiently small full-dimensional simplex
trace or a path visiting its vertices. A sufficiently large ball is
a universal cover. All hypotheses of \cref{thm:outer-net} follow.
\end{proof}

\begin{proposition}[Minkowski superadditivity and extremality]\label{prop:minkowski-extreme}
For fixed $T$, the unit translation-cover class is convex under
Minkowski averaging, and
\[
 E_T^{\rm tr}(K_0+K_1)\ge E_T^{\rm tr}(K_0)+E_T^{\rm tr}(K_1).
\]
If a minimum-volume cover is a nontrivial Minkowski average of two
feasible covers, those two covers are translates of each other.
\end{proposition}
\begin{proof}
If a loop trace $Q$ has placements $Q+z_i\subset K_i$, then
\[
 Q+(1-t)z_0+tz_1\subset(1-t)K_0+tK_1.
\]
This proves convexity of the unit-cover class. Write
$e_i=E_T^{\rm tr}(K_i)>0$ and use endpoint coverage of $K_i/e_i$.
Their Minkowski average with weights $e_i/(e_0+e_1)$ is a unit cover.
Strictification yields the displayed superadditivity.

If $K_*=(1-t)K_0+tK_1$, $0<t<1$, minimizes volume, each
$\vol(K_i)\ge\vol(K_*)$. Brunn--Minkowski gives
\[
 \vol(K_*)^{1/d}\ge(1-t)\vol(K_0)^{1/d}
       +t\vol(K_1)^{1/d}\ge\vol(K_*)^{1/d}.
\]
Equality forces homothety; equal volumes force the homothety ratio
to be one, so the bodies differ by translation. Existence of a
minimizer is proved in \cref{thm:outer-net}.
\end{proof}

\subsection{Finite polygon witnesses and two-sided computation}\label{sec:hierarchies}
\begin{proposition}[Finite polygon witnesses and algebraic certification]\label{prop:finite-witness}
For closed translation covers with a fixed gauge $T$, it suffices to
test polygonal loops with at most $d+1$ vertices and $T$-length at
most one. If $T$ is a rational polytope, the constant $W_T$ is a
computable real in the following concrete sense: there are finite
procedures producing rational lower and upper bounds with arbitrarily
small prescribed gap. Upper bounds can be certified by the LPs in
\cref{thm:lp}; lower bounds can be certified by finite families of
algebraic polygons and real-algebraic decision procedures.
This is a decidability and convergence statement, not a polynomial-time
complexity bound.
\end{proposition}
\begin{proof}
If a loop has no closed translate in $K$, the signed support margin
in \cref{thm:balance} is positive. Choosing its at most $d+1$
support points and shortcutting in cyclic temporal order preserves
that positive margin and does not increase length. Thus testing
these polygons is sufficient for closed coverage. Necessity is immediate.
Normalize their first vertex to zero. The resulting bounded family
of at most $d+1$ vertices, with repetitions allowed, is compact.

For a rational polytopal gauge this vertex family is a compact
semialgebraic set defined by rational data. An explicit finite
$\delta$-net can be constructed by partitioning a known bounding
box into small rational cubes, deciding which cubes meet the family,
and selecting a real-algebraic point in each nonempty intersection.
Choose the cubes so that corresponding vertices are within $\delta$
in Euclidean distance. Linear interpolation along corresponding edges
then puts the polygonal traces within Hausdorff distance $\delta$.
Include a fixed rational simplex loop of sufficiently small positive
length; it provides an explicitly positive inradius $\rho$.
The test family is now finite and has a quantitative error estimate
from \cref{thm:outer-net}.

For finitely many polytopes a minimum convex cover is the convex hull
of their placed vertices. Translations can be bounded by the radius
argument, after normalizing the retained simplex witness. There are
only finitely many labeled hull and triangulation types. On each
type, containment, orientation signs, and volume inequalities are
semialgebraic, the volume being a sum of simplex determinants with
prescribed signs. Quantifier elimination therefore brackets the
minimum volume to any rational accuracy. This gives lower bounds,
and \eqref{eq:outer-net} gives matching upper bounds.

Alternatively, rational polytopal bodies are Hausdorff dense, so
\cref{cor:density,thm:minkowski} show that their quotients approach
$W_T$. Each such quotient is exactly computable by the finite rational
LP formula. It is an upper bound for $W_T$; the corresponding scaled
polytope is an actual universal cover. Enumerating these candidates
and the finite-witness lower bounds yields arbitrarily close rational
brackets. One must certify all LP branches when claiming a lower
bound on an escape threshold; one favorable branch alone is insufficient.
\end{proof}

\begin{theorem}[Nested two-sided certificate hierarchy]\label{thm:hierarchy}
Fix a rational polytopal gauge $T$. There are nested finite families
$\mathscr A_m$ of algebraic polygonal test loops and nested nonempty
finite families $\mathscr P_m$ of rational polytopal bodies such that
\begin{align}
 L_m&:=\min\{\vol K:K\text{ covers }\mathscr A_m\}\uparrow W_T,
       \label{eq:hierarchy-lower}\\
 U_m&:=\min_{P\in\mathscr P_m}
       \frac{\vol P}{E_T^{\rm tr}(P)^d}\downarrow W_T.
       \label{eq:hierarchy-upper}
\end{align}
One can compute nondecreasing rational numbers $l_m\le L_m$ with
$l_m\to W_T$, and every $U_m$ is rational and exactly computable.
For every rational $\eps>0$, the stopping condition
$U_m-l_m\le\eps$ is eventually met and certifies
$l_m\le W_T\le U_m$.
\end{theorem}
\begin{proof}
Use the finite semialgebraic nets from \cref{prop:finite-witness}
with meshes tending to zero. Let $\mathscr A_m$ be their cumulative
union through step $m$, always including one fixed full-dimensional
simplex loop. Each family is finite; the sequence is nested and dense
in the compact normalized polygon family. The finite-polygon reduction
and \cref{thm:outer-net} prove \eqref{eq:hierarchy-lower}, including
attainment of every finite minimum.

Enumerate rational polytopes with nonempty interior containing zero in
their interiors, and let $\mathscr P_m$ contain the first $m$. Their
quotients are defined and rational by \cref{thm:lp}; rational polytopes
have rational volume. Each quotient is the volume of a universal
cover obtained by dividing $P$ by its escape threshold. Hence
$U_m\ge W_T$. The minima decrease, and quotient continuity and density
in \cref{cor:density} give \eqref{eq:hierarchy-upper}.

By the real-algebraic bracketing procedure in \cref{prop:finite-witness},
compute rational $a_m$ with $a_m\le L_m<a_m+2^{-m}$. Put
$l_m=\max\{0,a_1,\ldots,a_m\}$. Nested constraints ensure
$l_m\le L_m\le W_T$, and $a_m\to W_T$ implies $l_m\to W_T$.
Thus $U_m-l_m\to0$. Checking the rational stopping inequality is
finite at every step and terminates for each positive tolerance. The
argument gives certification and completeness, without an efficiency
claim for enumeration or quantifier elimination.
\end{proof}

\begin{corollary}[Polytopal minimizers and algebraic worm constants]\label{cor:polytopal-minimizer}
For a polytopal gauge $T$, a minimum universal unit-loop cover can be
chosen as a polytope with at most
\begin{equation}\label{eq:cover-vertex-bound}
 N_T=\sum_I |I|!
\end{equation}
vertices, where $I$ runs over its extreme positive facet circuits.
If the facet data of $T$ are real algebraic, then $W_T$ is real
algebraic and a minimizing cover has real-algebraic vertices.
In particular this holds for rational gauges and regular polygonal
gauges. This is an exact finite-family conclusion, not only a limit
of increasingly fine polygonal test nets.
\end{corollary}
\begin{proof}
First the hull of arbitrary placements of the complete family is
full-dimensional. The facet normals positively span $\R^d$ because
the gauge is bounded. Express each $-n_i$ as a nonnegative combination
of the normals and sum the resulting balances. This constructs a
strictly positive vector in the balanced cone. Decomposing it into
extreme rays shows that every normal occurs in some positive circuit.
An edge of the corresponding test polygon is a nonzero multiple of
that normal. Translation leaves these edge differences unchanged;
hence the differences of placed vertices span $\R^d$.

Now choose a minimum cover by \cref{cor:outer-families} and place every
member of the finite complete test family inside it. The convex hull
of their placed vertices is a convex body by the preceding argument,
and is still a universal cover by \cref{thm:complete-test-family}.
It therefore has the same minimum volume and at most the total number
of listed vertices, $\sum_I |I|(|I|-1)!=N_T$.

For algebraic facet data the circuits and test vertices are algebraic.
The inball witness and radius bound in \cref{thm:outer-net} allow
one minimizing configuration to be translated so that all placement
vectors lie in a fixed box with rational endpoints. A root of each
test polygon is zero, so each placement vector belongs to the cover
and is bounded by that radius. The volume of the hull of finitely
many translated algebraic vertices is a continuous semialgebraic
function of the placement vectors: finitely many hull and triangulation
types express it through determinant sums and orientation signs.
Let $F(x)$ denote this function on the compact rational box $B$.
Its graph is definable over the real algebraic field
$\R_{\rm alg}=\overline{\mathbb Q}\cap\R$. The sentence
\[
 \exists x\in B\ \exists v\;
   \bigl(v=F(x)\ \wedge\ \forall y\in B\ F(y)\ge v\bigr)
\]
is true over $\R$, by continuity and compactness. Quantifier
elimination for real closed fields implies that
$\R$ is an elementary extension of $\R_{\rm alg}$: a formula
with algebraic parameters has the same truth value in both fields
\cite[Chapter~2]{BasuPollackRoy}; see also
\cite[Theorem~2.1]{BasuSurvey}. Thus there are algebraic $x,v$
satisfying this sentence over $\R_{\rm alg}$. Transfer back to
$\R$ shows that this $x$ minimizes over \emph{all real} placements
in $B$, not just the algebraic ones. This proves both algebraicity
of the value and existence of an algebraic minimizing configuration. This proves the last assertions. The rational hierarchy
in \cref{cor:facet-hierarchy} supplies a different route to certified
brackets without using quantifier elimination.
\end{proof}

\begin{lemma}[Positive concave functions on pointed polyhedra]
\label{lem:positive-concave}
Let $P$ be a nonempty pointed polyhedron in a finite-dimensional
real vector space. If $F:P\to(0,\infty)$ is finite and concave,
then $F$ is nondecreasing on every recession ray, and
\[
 \min_{x\in P}F(x)=\min_{v\in\operatorname{Vert}P}F(v).
\]
In particular the minimum is attained at a vertex, even if $P$
is unbounded.
\end{lemma}
\begin{proof}
For $w$ in the recession cone and $0\le s_1<s_2$, a strict
decrease $F(x+s_2w)<F(x+s_1w)$ would, by concavity, give
\[
 F(x+sw)\le F(x+s_1w)+\frac{s-s_1}{s_2-s_1}
                  \bigl(F(x+s_2w)-F(x+s_1w)\bigr),\quad s>s_2.
\]
The right side becomes negative, a contradiction. The finite
polyhedral decomposition theorem \cite{Rockafellar} writes every
$x\in P$ as $\sum_j t_jv_j+w$, where $t_j\ge0$, $\sum_jt_j=1$,
the $v_j$ are vertices, and $w$ is a recession vector. A nonempty
pointed polyhedron has a nonempty finite vertex set. Therefore
\[
 F(x)\ge F\Bigl(\sum_jt_jv_j\Bigr)
       \ge\sum_jt_jF(v_j)\ge\min_jF(v_j).
\]
The finite minimum on the right is attained and belongs to $P$.
\end{proof}

\begin{theorem}[Exact minimum covers with prescribed facet directions]\label{thm:fixed-normals}
Fix a rational polytopal gauge $T$ and nonzero rational vectors
$n_1,\ldots,n_N$ positively spanning $\R^d$. Choose rational
generators $\lambda^I$ and let $c_I(T)$ be as in \cref{cor:circuit-lp}. The complete set of offsets
for universal unit $T$-loop covers with these facet directions is
the rational polyhedron
\begin{equation}\label{eq:offset-polyhedron}
 \mathcal D_T=\{b\in\R^N:\lambda^I\cdot b\ge c_I(T)
                        \text{ for every }I\}.
\end{equation}
Every $b\in\mathcal D_T$ defines a convex body $K(b)$. Remove
translation by fixing $b_{i_1}=\cdots=b_{i_d}=0$ for any $d$ linearly
independent normals. The resulting polyhedron $\mathcal D_T^0$ is
nonempty and pointed, and
\begin{equation}\label{eq:fixed-normal-min}
 V_T(n_1,\ldots,n_N):=\min_{b\in\mathcal D_T}\vol K(b)
   =\min_{b\in\operatorname{Vert}(\mathcal D_T^0)}\vol K(b).
\end{equation}
Consequently the restricted optimum is rational and can be computed
exactly by finitely many rational linear programs, vertex enumerations,
and determinant evaluations. Redundant facet inequalities are allowed.
\end{theorem}
\begin{proof}
For every nonzero nonnegative balanced vector $\lambda$, its
extreme-ray decomposition and \eqref{eq:offset-polyhedron} give
$\lambda\cdot b>0$. The strict finite alternative used in
\cref{prop:poly-circuit} therefore gives a point satisfying
$n_i\cdot x<b_i$ for all $i$. Positive spanning makes the intersection
bounded. Hence $K(b)$ is a convex body. By
\cref{cor:circuit-lp} and endpoint strictification, its universal-cover
property is exactly \eqref{eq:offset-polyhedron}. A sufficiently
large homothet of $K(1,\ldots,1)$ is feasible, so the polyhedron
is nonempty.

Write $A$ for the matrix with rows $n_i^T$. Its nonnegative balanced
cone spans $\ker A^T$. Indeed it contains a strictly positive vector:
express each $-n_i$ as a nonnegative combination of the normals and
sum the resulting balances. A small perturbation of this positive
vector in either direction of any vector in $\ker A^T$ stays
nonnegative. It follows that the lineality space of $\mathcal D_T$
is exactly $(\ker A^T)^\perp=\operatorname{im}A$. This is precisely
the change in offsets under translation of $K(b)$. The stated
independent-coordinate section intersects each translation class once
and removes all lineality.

On $\mathcal D_T^0$ the function $F(b)=\vol K(b)^{1/d}$ is positive
and concave. In fact
\[
 (1-t)K(b)+tK(b')\subset K((1-t)b+tb'),
\]
and Brunn--Minkowski proves concavity.
\Cref{lem:positive-concave} proves \eqref{eq:fixed-normal-min},
including attainment on the unbounded pointed polyhedron.
Rationality of the circumradius constants follows from
\cref{cor:circuit-lp}. The two polyhedral vertex computations are
rational, and triangulating each rational $K(b^{(j)})$ expresses its
volume as a finite sum of absolute rational determinants divided by
$d!$. This proves the effective assertion.
\end{proof}

\begin{proposition}[Exact facet-offset benchmarks]\label{prop:facet-benchmarks}
Put
\begin{align*}
 \Delta&=\conv\{(1,0),(0,1),(-1,-1)\},\\
 Q&=\conv\{(1,0),(0,1),(-1,0),(0,-2)\},\\
 H&=\conv\{(1,0),(1,1),(0,1),(-1,0),(-1,-1),(0,-1)\}.
\end{align*}
Let $J(x,y)=(-y,x)$. The following covers have unit escape threshold
and minimize area both for the listed facet directions and without
that restriction:
\begin{center}\small
\begin{tabular}{@{}llll@{}}
\toprule
Gauge $T$ & Forest normals & Minimum cover $K$ & $W_T$\\
\midrule
$\Delta$ & $\pm e_1,\pm e_2,\pm(e_1+e_2)$
& $\conv(J\Delta\cup-J\Delta)/3$ & $1/3$\\
$[-1,1]^2$ & $\pm(1,1),\pm(1,-1)$
& $\{|x+y|,|x-y|\le1/4\}$ & $1/8$\\
$Q$ & $\pm e_1,\pm e_2$
& $[-1/4,1/4]\times[-1/6,1/6]$ & $1/6$\\
$H$ & $(-1,-2),(2,1),(-1,1)$
& $J\Delta/3$ & $1/6$\\
\bottomrule
\end{tabular}
\end{center}
All four values are known sharp cases: the simplex and parallelotope
results are in \cite[Theorems~5.1--5.2]{Balitskiy2020}; the planar
quadrilateral and affinely regular hexagon bounds are in
\cite[Theorem~1.2]{BMP}. The purpose of this proposition is exact
validation of the facet-offset pipeline.
\end{proposition}
\begin{proof}
Use the explicit normal order specified below for the triangular gauge;
for the other rows, order the positive normals before their negatives.
For $\Delta$, its own three facet normals are
$(1,1),(-2,1),(1,-2)$ with offsets one. The balanced-support formula
gives, for every finite set $S$,
\[
 \operatorname{rad}_\Delta(S)
 =\tfrac13\bigl(h_S(1,1)+h_S(-2,1)+h_S(1,-2)\bigr).
\]
Its use on the forest state pairs gives $2/3$; on each of the two
three-normal circuits, both cyclic orders give one. Thus the
offset inequalities are
\[
 b_1+b_3\ge2/3,\quad b_2+b_4\ge2/3,\quad b_5+b_6\ge2/3,
 \quad b_1+b_2+b_6\ge1,\quad b_3+b_4+b_5\ge1,
\]
where the normal order is $e_1,e_2,-e_1,-e_2,e_1+e_2,-e_1-e_2$.
After $b_1=b_2=0$, these five boundary equations are consistent
and have rank four. Every independent four-row subsystem therefore
has the same solution
$(b_3,b_4,b_5,b_6)=(2/3,2/3,-1/3,1)$.
It is feasible and is the only vertex. Translating gives all six
offsets $1/3$, hence the stated hexagon of area $1/3$.

For the square gauge, the two opposite-normal state pairs have
circumradius $1/2$. The opposite-facet sums $w_1,w_2$ must therefore
be at least $1/2$. The determinant of the two independent forest
normals has absolute value two, giving area $w_1w_2/2\ge1/8$.
For $Q$, the horizontal and vertical state pairs have circumradii
$1/2$ and $1/3$: their projection widths force these bounds, and
the horizontal and vertical central chords of $Q$ attain them.
The rectangle constraints are therefore $w_1\ge1/2$, $w_2\ge1/3$,
with minimum area $1/6$.

For $H$, the single circuit has weights $(1,1,1)$. Its state
triangles, up to translation, are
\[
 \{0,(1,2),(-1,1)\},\qquad \{0,(1,2),(2,1)\}.
\]
Their centers $(0,1)$ and $(1,1)$ put them in $H$; their vertical
width two forces circumradius at least one. Thus the only offset
condition is $b_1+b_2+b_3\ge1$. The area is
$(b_1+b_2+b_3)^2/6$, as follows by intersecting the three lines
and taking their determinant. At offsets $1/3$ this is $J\Delta/3$.

The facet formula proves escape one for each displayed cover.
The restricted minima follow either from these calculations or
\cref{thm:fixed-normals}. Finally the cited sharp product inequalities,
the published capacity identity, and the quotient formula give the
unrestricted lower bounds $1/(2\area T)$, which equal the displayed
areas because $\area\Delta=3/2$, $\area([-1,1]^2)=4$, and
$\area Q=\area H=3$. This last step uses established global bounds;
restricted facet optimization alone would not prove global optimality.
\end{proof}

The supplied exact verifier enumerates every positive circuit,
cyclic order, and independent active offset subsystem in these
four cases. The respective counts are $5,2,2,1$ circuits,
$7,2,2,2$ circumradius programs, and one feasible offset vertex
per case. Primal and dual rational radius certificates are replayed
independently. These examples demonstrate the full finite computation
on small instances, rather than a practical complexity bound.
Leipold--Vallentin \cite{LeipoldVallentin} prove NP-hardness of
EHZ-capacity computation for general rational polytopes, even
simplices in the ambient symplectic space. Their result is not
asserted here for the special subclass of Lagrangian products.

\begin{corollary}[Two-sided rational approximation from facet meshes]\label{cor:facet-hierarchy}
Let $d\ge2$ and $T$ be a rational polytopal gauge. Suppose a unit
$T$-loop witness has convex hull containing a radius-$\rho$ ball,
and a known universal cover has volume at most $V_+$. Put
\[
 R=\frac{dV_+}{\kappa_{d-1}\rho^{d-1}},
 \qquad \kappa_j=\vol_j(B_2^j).
\]
For a rational unit-normal $\eta$-net $U\subset\Sph^{d-1}$,
$0<\eta\le1$, which positively spans $\R^d$, let $V_T(U)$ be the
exact optimum in \cref{thm:fixed-normals}. Then
\begin{equation}\label{eq:facet-bracket}
 \frac{V_T(U)}{(1+3R\eta/\rho)^d}
       \le W_T\le V_T(U).
\end{equation}
Nested such normal sets with mesh tending to zero have exact optima
decreasing to $W_T$. In particular, two-sided rational certificates
for $W_T$ can be obtained using linear programming and determinants,
without real-algebraic quantifier elimination.
\end{corollary}
\begin{proof}
Normalize a minimum cover $K_*$ to contain $\rho B_2^d$.
The radius estimate \eqref{eq:outer-radius} gives $K_*\subset RB_2^d$.
Circumscribe the polytope
$P=\{x:u\cdot x\le h_{K_*}(u),\ u\in U\}$.
For $x\in P\setminus\{0\}$ choose $u\in U$ within $\eta$ of
$x/|x|$. Since both directions are unit,
\[
 |x|(1-\eta^2/2)\le u\cdot x\le R,
 \qquad |x|\le2R.
\]
For an arbitrary unit $v$, choose $u\in U$ with $|u-v|\le\eta$.
The support function of $K_*$ is $R$-Lipschitz, so
\[
 v\cdot x-h_{K_*}(v)
 \le |v-u|\,|x|+h_{K_*}(u)-h_{K_*}(v)
 \le3R\eta.
\]
Thus
\[
 K_*\subset P\subset K_*+3R\eta B_2^d
                    \subset(1+3R\eta/\rho)K_*.
\]
The polytope $P$ is a feasible competitor for the fixed-normal problem,
even though its displayed offsets need not be rational. This gives
$W_T\le V_T(U)\le(1+3R\eta/\rho)^dW_T$.
Adding normals preserves every old feasible polytope by allowing
redundant inequalities, proving monotonicity.

All auxiliary bounds can be made rational. For example, take the loop
through $0,\eps e_1,\ldots,\eps e_d,0$, with rational $\eps>0$
small enough that its rational $T$-length is at most one. Its hull
contains the ball centered at $(\rho,\ldots,\rho)$ with
$\rho=\eps/(2d)$. Indeed its distances to the coordinate facets are
$\rho$, and its distance to $\sum_i x_i=\eps$ is
$(\eps-d\rho)/\sqrt d=\eps/(2\sqrt d)\ge\rho$.
Choose rational $r>0$
with $rB_2^d\subset T$; this follows, for example, by inscribing
a rational cube using the facet inequalities of $T$. A unit $T$-loop
has Euclidean length at most $1/r$; after translating its root to
zero it lies in the cube $[-1/(2r),1/(2r)]^d$. Hence one may take
$V_+=r^{-d}$. Finally replace $\kappa_{d-1}$ in $R$ by the rational
lower bound $(2/(d-1))^{d-1}$, the volume of a cube contained in the unit ball.
This enlarges $R$ and preserves \eqref{eq:facet-bracket}.
For rational certificate endpoints choose rational $\eta>0$ as well.
An explicit rational mesh is obtained by choosing an integer $M\ge d/\eta$
and taking the images of $z\in M^{-1}\mathbb Z^{d-1}\cap[-1,1]^{d-1}$ under
\[
 F_\pm(z)=\left(\frac{2z}{1+|z|^2},
          \ \pm\frac{1-|z|^2}{1+|z|^2}\right).
\]
Every direction in the respective hemisphere has a preimage with $|z|\le1$.
Direct expansion gives the chord-distance identity
\[
 |F_\pm(z)-F_\pm(w)|
 =\frac{2|z-w|}{\sqrt{(1+|z|^2)(1+|w|^2)}}\le2|z-w|.
\]
Rounding each preimage coordinate to the nearest grid coordinate
changes $z$ by at most $\sqrt{d-1}/(2M)$ and stays in the displayed
cube. The image error is therefore at most $\sqrt{d-1}/M\le\eta$.
The set includes the coordinate directions and their negatives, so
it positively spans the space. With the rational enlarged $R$, these
choices make both endpoints in \eqref{eq:facet-bracket} rational.
Their gap tends to zero; taking successive unions makes the meshes nested.
\end{proof}

This hierarchy strengthens the earlier finite-witness construction:
its upper step optimizes all offsets for the chosen normals, and the
same step supplies a lower certificate through the explicit mesh bound.
It is a finite exact procedure, with no polynomial-time claim; the
numbers of circuits, orders, and polyhedral vertices may be large.

\subsection{Regular-polygon gauge approximation}\label{sec:regular}
\begin{corollary}[Euclidean gauge approximation]\label{cor:gauge-rate}
In the plane let $Q_N$ be the regular circumscribed $N$-gon about
the unit disk, $N\ge3$. Then
\begin{equation}\label{eq:regular-gauge}
 W_{Q_N}\le\Wet\le\sec^2(\pi/N)W_{Q_N},\qquad
 0\le\frac{\Wet-W_{Q_N}}{W_{Q_N}}\le\tan^2(\pi/N).
\end{equation}
Thus the gauge discretization has a certified relative error
$O(N^{-2})$, independently of the subsequent body optimization.
\end{corollary}
\begin{proof}
The inclusions $B_2^2\subset Q_N\subset\sec(\pi/N)B_2^2$
and \cref{thm:gauge} give the result. Here $\Wet=W_{B_2^2}$.
\end{proof}

\begin{remark}[A conditional dodecagon benchmark]\label{rem:dodecagon}
\Cref{prop:circuit-count} gives 46 facet circuits for a regular
dodecagon. For $Q_{12}$ circumscribed about the unit disk,
\[
 \area(Q_{12})=12(2-\sqrt3),\qquad
 \frac1{2\area(Q_{12})}=\frac{2+\sqrt3}{24}.
\]
If the fixed-gauge inequality
$W_{Q_{12}}\ge(2+\sqrt3)/24$ were established, then
\eqref{eq:regular-gauge} would yield the same lower bound for
$\Wet$. This implication is conditional on that fixed-gauge inequality. Counting circuits does not prove the required outer
volume inequality. The unconditional improvement proved next uses
a different geometric argument.
\end{remark}

\section{An explicit strengthening of Wetzel's lower bound}\label{sec:wetzel}
Throughout this section $K\subset\R^2$ is a convex translation cover
of every closed curve of Euclidean length at most one, and
$\Wet$ denotes the minimum possible area. Degenerate doubled segments
are among the allowed closed curves. If one starts with only curves
of length exactly one, padding by an out-and-back excursion shows
that the same convex-cover problem results.

Wetzel's comparison function is
\begin{equation}\label{eq:wetzel-f}
 f(r)=\pi r^2+3r\sqrt{\tfrac14-r}
       -3r^2\arccos\!\left(\frac{r}{\frac12-r}\right),
 \qquad \frac16\le r\le\frac14.
\end{equation}
We give the comparison construction and all additional estimates,
so that the strict quantitative improvement does not depend on a
diagrammatic equality-case assertion.

\subsection{Geometric restrictions and exclusion of the classical lower bound}
\begin{lemma}[Width, inradius, and balanced contacts]\label{lem:wetzel-inradius}
Every such cover satisfies
\begin{equation}\label{eq:wetzel-width}
 w(K):=\min_{u\in\Sph^1}\bigl(h_K(u)+h_K(-u)\bigr)\ge\frac12,
 \qquad r(K)\ge r_0:=\frac{\sqrt3}{9}.
\end{equation}
If a largest inscribed disk has radius $r<1/4$ and is centered at
zero, it has three supporting unit normals $u_1,u_2,u_3$ that
positively balance, with each consecutive angular gap less than $\pi$.
\end{lemma}
\begin{proof}
For each direction, the doubled segment of length $1/2$ has total
curve length one; its translation into $K$ proves the width inequality.
For a largest inball $rB_2^2$, let
$U=\{u:h_K(u)=r\}$ be the set of active normals. If
$0\notin\conv U$, a strict separating direction permits a small
translation of the center that increases every active clearance.
The separating inequality remains uniformly positive near the compact
active set; on the complementary closed set the original clearances
exceed $r$ by a positive amount. A sufficiently small center translation
therefore increases the minimum clearance, contradicting maximality. Hence $0\in\conv U$.
Carath\'eodory supplies a minimal balancing set with two or three
members. Two members are antipodal and give width $2r$, so
$r\ge1/4$ in that case.

In the three-member case the balance has positive coefficients and
each angular gap is less than $\pi$. The triangle with vertices
$ru_i$ cannot be translated into $\Int K$: a translating vector $z$
would satisfy $z\cdot u_i<0$ for all three normals, contradicting
their positive balance. Its perimeter is therefore at least one,
by \cref{lem:strictification}. If the consecutive angular gaps are
$\theta_1,\theta_2,\theta_3$, their sum is $2\pi$, and its
perimeter is
\[
 2r\sum_{i=1}^3\sin(\theta_i/2)
 \le6r\sin(\pi/3)=3\sqrt3\,r
\]
by concavity of sine on $[0,\pi]$. Thus $r\ge r_0$.
If $r<1/4$, the two-member alternative is impossible.
\end{proof}

\begin{lemma}[Area, diameter, and width]\label{lem:area-diameter}
Every planar convex body obeys
$\area(K)\ge\tfrac12\diam(K)w(K)$.
\end{lemma}
\begin{proof}
Take diameter endpoints as a common base and support points
on the two parallel supporting lines on opposite sides of that base.
The two resulting triangles lie in $K$, have disjoint interiors, and
have combined area equal to one-half the base length times the width
in its perpendicular direction. In particular a Wetzel cover has
$\diam(K)\le4\area(K)$, an additional compactness bound.
\end{proof}

\begin{lemma}[Three disjoint caps and their area]\label{lem:wetzel-caps}
Suppose $r_0\le r<1/4$ and $rB_2^2$ is a largest inball of $K$.
Put
\[
 a=\tfrac12-r,\qquad s=\sqrt{\tfrac14-r},\qquad
 \alpha=\arccos(r/a).
\]
There exist three points $y_1,y_2,y_3\in K$, each at distance $a$
from zero, for which the open angular intervals
\[
 I_i=\{u\in\Sph^1:y_i\cdot u>r\}
\]
are pairwise disjoint. Consequently
\begin{equation}\label{eq:phi-area}
 \Phi:=\conv(rB_2^2,y_1,y_2,y_3)\subset K,
 \qquad \area(\Phi)=f(r).
\end{equation}
Moreover, $f$ is strictly increasing on $[1/6,1/4]$.
\end{lemma}
\begin{proof}
Choose the balanced contact normals from
\cref{lem:wetzel-inradius}, in cyclic order. The width bound gives
points $x_i\in K$ with $u_i\cdot x_i\le-a$. In particular
$|x_i|\ge a$. Set $y_i=ax_i/|x_i|$, which belongs to $K$ by
convexity and $0\in K$.

We first locate its direction. The direction of $x_i$ cannot lie in
either gap adjacent to $u_i$. Indeed, rotate one such gap to
$[0,g]$, with $u_i$ at zero, the adjacent normal at $g<\pi$,
and $x_i$ at angle $\theta\in[0,g]$. The negative projection onto
$u_i$ forces $\theta>\pi/2$. Since
$0\le g-\theta<\pi-\theta<\pi/2$,
\[
 u_{\rm adjacent}\cdot x_i
   =|x_i|\cos(g-\theta)
   >-|x_i|\cos\theta\ge a>r,
\]
contradicting its supporting inequality for $K$. The other adjacent
gap is treated by reversing angular orientation. Thus $y_i$ points
into the gap between the other two normals which does not contain
$u_i$.

Every contact normal has $u_j\cdot y_i\le r$. The connected open
interval $I_i$, centered on the direction of $y_i$, therefore contains
none of the three normals and is confined to that opposite gap.
The three intervals are consequently disjoint. Each has length
$2\alpha$. Since
$h_\Phi(u)=\max\{r,y_1\cdot u,y_2\cdot u,y_3\cdot u\}$,
each open interval $I_i$ is exactly the set where its point $y_i$
contributes the support above the disk. On the complement the disk
supplies the support. Thus the boundary consists of the retained disk
arcs and the two tangents from each $y_i$ to the disk, with coincident
tangency endpoints allowed. No additional region between distinct
caps is introduced, and their interiors are disjoint. The area added by one cap is
\[
 r\sqrt{a^2-r^2}-r^2\alpha=rs-r^2\alpha,
\]
as follows by subtracting the disk sector from the two right
triangles at its tangencies. Adding the disk proves
\eqref{eq:phi-area}.

For $1/6\le r<1/4$, differentiation gives
\begin{equation}\label{eq:f-prime}
 f'(r)=2r(\pi-3\alpha)+3s\frac{a-r}{a}>0.
\end{equation}
Here $\alpha\le\pi/3$, $s>0$, and $a-r>0$; at the endpoint
$r=1/4$, continuity completes strict monotonicity. For example,
the derivative identity follows from $s'=-1/(2s)$ and
$\alpha'=-1/(2as)$.
\end{proof}

\begin{lemma}[Uniform width deficit of the comparison bodies]\label{lem:wetzel-defect}
For every comparison body in \cref{lem:wetzel-caps} with
$r_0\le r\le1/5$,
\begin{equation}\label{eq:phi-width}
 w(\Phi)\le b_0:=\frac32r_0+
       \frac{\sqrt3}{2}\sqrt{\frac14-r_0}<\frac12.
\end{equation}
More precisely, there is a unit vector $u$ with
$h_\Phi(u)=r$ and $h_\Phi(u)+h_\Phi(-u)\le b_0$.
This estimate is independent of the three point directions.
\end{lemma}
\begin{proof}
Let $g$ be the largest angular gap between consecutive directions of
the three points. Then $g\ge2\pi/3$. Between the corresponding cap
intervals there is a closed interval $J$ of disk-support directions
of length
\[
 L=g-2\alpha\ge2\pi/3-2\alpha>0,
 \qquad h_\Phi(u)=r\quad(u\in J).
\]
If the antipodal interval $-J$ contains a direction outside all
open cap intervals, some opposite pair has width $2r$. Otherwise
the connected interval $-J$ lies in a single cap interval. One of
its endpoints is at angular distance at least $L/2$ from the center
of that cap. At the opposite pair through that endpoint the width
is at most
\[
 r+a\cos(L/2)
 \le r+a\cos(\pi/3-\alpha)
 =\frac32r+\frac{\sqrt3}{2}\sqrt{\frac14-r}=:B(r).
\]
In the second case $L/2\le\alpha<\pi/2$, so the stated
monotonicity of cosine applies. We have proved
$w(\Phi)\le\max\{2r,B(r)\}$.

The derivative
\[
 B'(r)=\frac32-\frac{\sqrt3}{4\sqrt{\frac14-r}}
\]
is negative for $r>1/6$, and $B(1/6)=1/2$. Hence
$B(r)\le B(r_0)=b_0<1/2$. Moreover
$b_0\ge B(1/5)=3/10+\sqrt{15}/20>49/100>2/5$;
the strict radical inequality follows from $15>(19/5)^2$.
Thus $2r\le2/5<b_0$ on the specified interval, proving \eqref{eq:phi-width}.
In both cases the selected opposite pair contains a direction in
$J$, where $h_\Phi=r$; orienting that direction as $u$ proves
the asserted tangency as well.
\end{proof}

\begin{proposition}[Attainment and strict exclusion of Wetzel's lower bound]
\label{prop:wetzel-strict}
The minimum $\Wet$ is attained, and $\Wet>f(r_0)$.
\end{proposition}
\begin{proof}
Attainment follows from \cref{cor:outer-families}. It can also be
seen directly from \cref{lem:area-diameter}: bounded area and the
width lower bound control diameter, so translated minimizing covers
have a Blaschke-convergent subsequence. For each fixed curve, bounded
translating vectors preserve its placement in the limit; volume is
continuous. The resulting body is a minimum cover.

Suppose some cover $K$ had area $f(r_0)$. Its inradius $r$ is at
least $r_0$. If $r\ge1/4$, its area is at least
$\pi/16=f(1/4)>f(r_0)$, a contradiction. Otherwise the comparison
body $\Phi\subset K$ of \cref{lem:wetzel-caps} has area $f(r)$.
Strict monotonicity of $f$ forces $r=r_0$ and equal areas of $K$
and $\Phi$. A proper inclusion of full-dimensional convex bodies
strictly increases volume: a point outside the smaller body, strictly
separated from it, and an interior ball generate a positive-volume
cone outside that body. Thus $K=\Phi$. But
\cref{lem:wetzel-defect} gives $w(K)\le b_0<1/2$, contradicting
\cref{lem:wetzel-inradius}. Every cover consequently has area greater
than $f(r_0)$. Attainment then implies the strict inequality for the
minimum itself, rather than only for each individual competitor.
\end{proof}

\subsection{Quantifying the width deficit by additional area}
\begin{lemma}[Optimized disk-chord area bound]\label{lem:cap-area}
Let $\Phi,K$ be convex bodies with $rB_2^2\subset\Phi\subset K$.
For a unit vector $u$,
write $H=h_\Phi(u)\ge r>0$. If $h_K(u)\ge H+\delta$ with
$\delta>0$, then
\begin{equation}\label{eq:extra-cap}
 \area(K)\ge\area(\Phi)+
       \frac{r\delta^2}{\sqrt{(H+\delta)^2-r^2}}.
\end{equation}
\end{lemma}
\begin{proof}
Choose $x\in K$ with $u\cdot x=H+\delta$, by contracting a support
point toward zero if necessary. Put $A=H+\delta$ and
$t=r^2/A<r\le H$. For a unit vector $v\perp u$, the disk chord
with endpoints $tu\pm\sqrt{r^2-t^2}\,v$ has half-length
$\sqrt{r^2-t^2}$. Its triangle with apex $x$ lies in $K$.
The part above $u\cdot y=H$ has height $\delta$ and area
\[
 \frac{\sqrt{r^2-t^2}\,\delta^2}{A-t}
 =\frac{r\delta^2}{\sqrt{A^2-r^2}}.
\]
It lies outside the interior of $\Phi$. The choice $t=r^2/A$
maximizes $\sqrt{r^2-t^2}/(A-t)$ over $-r\le t\le r$;
the displayed identity alone suffices for the bound. The tangential
component of $x$ is unrestricted.
\end{proof}

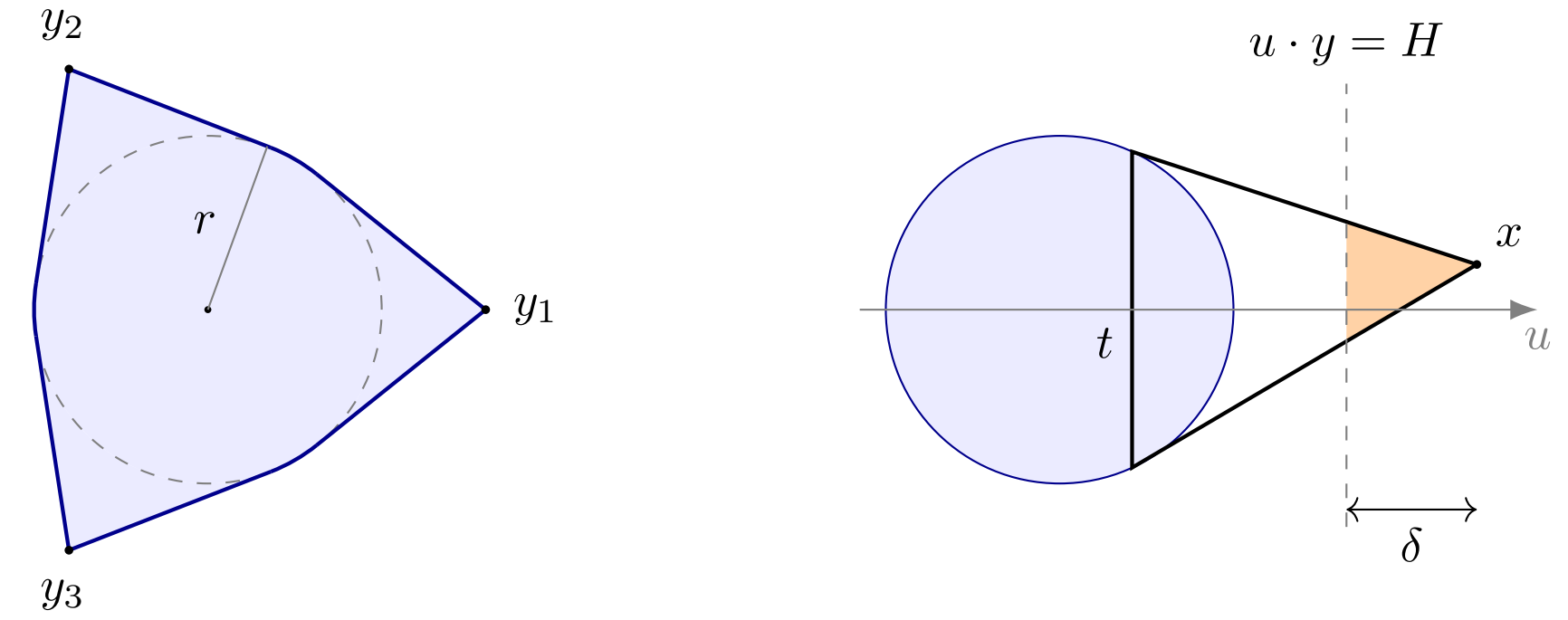
\begin{figure}[htbp]
\centering
\begin{tikzpicture}[scale=1.27,font=\small]
\begin{scope}[shift={(-2.6,0)}]
\def\aa{1.598076211}
\pgfmathsetmacro{\al}{acos(1/\aa)}
\fill[blue!8] (0,0) circle (1);
\foreach \ang in {0,120,240}{
  \fill[blue!8] (\ang:\aa) -- ({\ang-\al}:1)
       arc[start angle=\ang-\al,end angle=\ang+\al,radius=1] -- cycle;
}
\draw[gray,dashed] (0,0) circle (1);
\foreach \ang in {0,120,240}{
  \draw[blue!55!black,thick] ({\ang-\al}:1) -- (\ang:\aa)
        -- ({\ang+\al}:1);
  \draw[blue!55!black,thick] ({\ang+\al}:1)
        arc[start angle=\ang+\al,end angle=\ang+120-\al,radius=1];
  \fill (\ang:\aa) circle (0.025);
}
\fill (0,0) circle (0.02);
\draw[gray] (0,0) -- (70:1);
\node at (92:0.5) {$r$};
\node[anchor=west] at (1.65,0) {$y_1$};
\node[anchor=south] at (120:1.67) {$y_2$};
\node[anchor=north] at (240:1.67) {$y_3$};
\node at (0,-2.0) {Three-cap comparison body $\Phi$};
\end{scope}
\begin{scope}[shift={(2.3,0)}]
\def\tt{0.416666667}
\def\bb{0.909059343}
\def\hh{1.65}
\def\xx{2.4}
\def\yy{0.26}
\pgfmathsetmacro{\yu}{\yy+(\bb-\yy)*(\xx-\hh)/(\xx-\tt)}
\pgfmathsetmacro{\yl}{\yy+(-\bb-\yy)*(\xx-\hh)/(\xx-\tt)}
\fill[blue!8] (0,0) circle (1);
\fill[orange!35] (\hh,\yl) -- (\xx,\yy) -- (\hh,\yu) -- cycle;
\draw[blue!55!black] (0,0) circle (1);
\draw[thick] (\tt,-\bb) -- (\xx,\yy) -- (\tt,\bb) -- cycle;
\draw[gray,dashed] (\hh,-1.25) -- (\hh,1.3);
\draw[gray,-{Latex[length=2mm]}] (-1.15,0) -- (2.75,0)
      node[below] {$u$};
\fill (\xx,\yy) circle (0.025) node[above right] {$x$};
\node[below left] at (\tt,0) {$t$};
\node[above] at (\hh,1.3) {$u\cdot y=H$};
\draw[<->] (\hh,-1.15) -- (\xx,-1.15)
      node[midway,below] {$\delta$};
\node at (0.8,-2.0) {Area beyond a supporting line};
\end{scope}
\end{tikzpicture}
\caption{Left: the comparison body consists of an inscribed disk and
three disjoint caps. Right: a disk chord at $t=r^2/(H+\delta)$
and a point of excess support generate the shaded triangular cap.
The right panel illustrates the general lemma and uses schematic
values, not the small numerical excess in the final bound.}
\label{fig:wetzel}
\end{figure}

\begin{lemma}[Combining opposite support deficits]\label{lem:two-cap}
Let $\Phi,K$ be convex bodies with $rB_2^2\subset\Phi\subset K$,
$0<r<1/4$, and suppose $K$
has width at least $1/2$ in a direction in which $\Phi$ has
width $w<1/2$. Then
\begin{equation}\label{eq:two-cap}
 \area(K)-\area(\Phi)\ge
       \frac{r(1/2-w)^2}{2\sqrt{1/16-r^2}}.
\end{equation}
\end{lemma}
\begin{proof}
Let the two supports of $\Phi$ be $H_+,H_-\ge r$, with
$H_++H_-=w$. The two corresponding support excesses of $K$
are nonnegative and have sum at least $1/2-w$. Choose numbers
$d_+,d_-$ no larger than those excesses such that
$d_++d_-=1/2-w$, and write $A_\pm=H_\pm+d_\pm$.
Thus $A_++A_-=1/2$.

Apply \cref{lem:cap-area} on both sides. The two added triangular
caps lie in opposite open half-planes outside $\Phi$ and are
disjoint. Put $w_\pm=\sqrt{A_\pm^2-r^2}$. If $w_\pm=0$,
then $A_\pm=H_\pm=r$ and $d_\pm=0$, so that side contributes
zero area and is omitted from the sum. At least one $w_\pm$ is
positive, because $A_++A_-=1/2>2r$. Cauchy's inequality on the
positive-weight terms gives
\[
 \area(K)-\area(\Phi)
 \ge r\sum_{\pm:\,w_\pm>0}\frac{d_\pm^2}{w_\pm}
 \ge\frac{r(d_++d_-)^2}{w_++w_-}.
\]
The function
$A\mapsto\sqrt{A^2-r^2}$ is concave on $[r,\infty)$, since
its second derivative is $-r^2/(A^2-r^2)^{3/2}$ for $A>r$.
Its denominator sum is therefore at most
$2\sqrt{(1/4)^2-r^2}$. This proves \eqref{eq:two-cap}.
\end{proof}

\begin{lemma}[Area forced by a tangent deficient support]
\label{lem:tangent-deficit}
Suppose $r_0\le r\le1/5$, $rB_2^2\subset\Phi\subset K$,
every width of $K$ is at least $1/2$, and a unit vector $u$ satisfies
$h_\Phi(u)=r$ and $h_\Phi(u)+h_\Phi(-u)\le b_0$.
Put
\[
 D=1/2-b_0,\qquad A_0=\frac{r_0}{\sqrt{1/4-r_0}},\qquad
 B_0=\frac{r_0}{\sqrt{2r_0+D}}.
\]
Then
\begin{equation}\label{eq:tangent-deficit}
 \area(K)-\area(\Phi)
 \ge\min_{0\le x\le D}\{A_0(D-x)^2+B_0x^{3/2}\}.
\end{equation}
\end{lemma}
\begin{proof}
Write $H=h_\Phi(-u)$. The two nonnegative support excesses
$h_K(u)-r$ and $h_K(-u)-H$ sum to at least $D$.
Choose $x,y\ge0$ not exceeding these excesses and with $x+y=D$.
For $x>0$, \cref{lem:cap-area} gives an added cap on the $u$ side
of area at least
\[
 \frac{rx^2}{\sqrt{(r+x)^2-r^2}}
 =\frac{rx^{3/2}}{\sqrt{2r+x}}
 \ge\frac{r}{\sqrt{2r+D}}x^{3/2}.
\]
For $y>0$, the opposite cap has area at least
$ry^2/\sqrt{(H+y)^2-r^2}$. Since
$H+y\le b_0-r+D-x\le1/2-r$, this is at least
$ry^2/\sqrt{1/4-r}$. The two caps lie in the disjoint half-planes
$u\cdot z>r$ and $u\cdot z<-H$, so their areas add.
If an excess is zero its contribution is zero; the displayed
positive-excess formulas have positive denominators.

Both coefficients increase with $r$: their derivatives are
\[
 \frac{d}{dr}\frac r{\sqrt{1/4-r}}
 =\frac{1/4-r/2}{(1/4-r)^{3/2}}>0,\qquad
 \frac{d}{dr}\frac r{\sqrt{2r+D}}
 =\frac{r+D}{(2r+D)^{3/2}}>0.
\]
Replacing them by their values at $r_0$, and then minimizing over
$x\in[0,D]$, proves \eqref{eq:tangent-deficit}.
\end{proof}

\begin{theorem}[Explicit improvement of Wetzel's lower bound]\label{thm:wetzel-improvement}
Let $r_0=\sqrt3/9$, let $b_0$ be as in \cref{lem:wetzel-defect},
and use $D,A_0,B_0$ from \cref{lem:tangent-deficit}. Define
\begin{align}
 \delta_0&=D/2,\qquad c_0=8\delta_0^2/\sqrt{11},\label{eq:c0}\\
 t_0&=\frac{8A_0D}{\sqrt{9B_0^2+64A_0^2D}+3B_0},\qquad
 c_1=A_0(D-t_0^2)^2+B_0t_0^3.\label{eq:c1}
\end{align}
Every convex translation cover of closed unit-length planar curves
satisfies
\begin{equation}\label{eq:wetzel-improved}
 \area(K)\ge f(r_0)+c_1>0.15545785378555.
\end{equation}
In particular the same lower bound holds for $\Wet$. Exact decimal
enclosures are
\begin{align*}
 0.000009944944227372766220899
 &<c_1<0.000009944944227372766220900,\\
 0.1554578537855501720883
 &<f(r_0)+c_1<0.1554578537855501720884.
\end{align*}
Moreover $1.2946<c_1/c_0<1.2947$.
\end{theorem}
\begin{proof}
The function $G(x)=A_0(D-x)^2+B_0x^{3/2}$ is strictly convex
on $[0,D]$. Its one-sided derivative at zero is $-2A_0D<0$,
and $G'(D)>0$. Its unique minimum is therefore interior.
Writing $x=t^2$, its critical equation is
\[
 4A_0t^2+3B_0t-4A_0D=0.
\]
The positive root, after rationalization, is exactly $t_0$ in
\eqref{eq:c1}. Thus $0<t_0^2<D$ and $\min G=c_1>0$.

Center a largest inball of $K$ at zero. Its radius $r$ is at least
$r_0$ by \cref{lem:wetzel-inradius}. For $r\le1/5$, the comparison
body has area $f(r)$ and a tangent deficient support by
\cref{lem:wetzel-caps,lem:wetzel-defect}. Apply
\cref{lem:tangent-deficit} and the increase of $f$ to obtain
$\area K\ge f(r)+c_1\ge f(r_0)+c_1$.

For $1/5<r<1/4$, the comparison gives $\area K\ge f(1/5)$;
for $r\ge1/4$ the disk gives $\area K\ge\pi/16=f(1/4)>f(1/5)$.
These branches also exceed the claimed bound by the following exact
inequalities. Since $r_0<39/200$, $1/5-r_0>1/200$.
On $[r_0,1/5]$, $s>1/5$, $a-r\ge1/10$, and $a<1/3$,
so \eqref{eq:f-prime} gives $f'(r)>9/50$. Hence
$f(1/5)-f(r_0)>9/10000$.
Also $b_0>49/100$, so $D<1/100$. Since $r_0<1/5$,
$A_0<1$ (indeed $r_0^2+r_0<1/25+1/5<1/4$).
Consequently $c_1\le G(0)=A_0D^2<1/10000$, which suffices.

For comparison, discarding tangency and using \cref{lem:two-cap}
in the low-inradius branch gives the weaker increment
\[
 \frac{r_0D^2}{2\sqrt{1/16-r_0^2}}
 =\frac{2D^2}{\sqrt{11}}=c_0.
\]
The rational intervals in \cref{app:arithmetic} verify all displayed
decimals and $1.2946<c_1/c_0<1.2947$ using exact fractions.
They are used for numerical reporting and comparison of the two
increments; the geometric inequality itself precedes them.
\end{proof}

\begin{remark}[Comparison with the literature]
Wetzel's original lower bound is
$f(r_0)=0.1554479088413227993220\ldots$ \cite{Wetzel1973}.
The increment in \cref{thm:wetzel-improvement} is strictly positive
and explicit; mere exclusion of equality would not by itself produce
this number. The theorem concerns translation covers of closed curves,
not congruence covers of open unit arcs. It should therefore not be
compared numerically with the approximately $0.26$ areas appearing
in the latter problem.

Wetzel's 1973 asserted upper bound $0.15900$ is not a valid benchmark;
Rudolf discusses the error \cite{Rudolf}. Bezdek--Connelly and
Mall\'ee give later constructions \cite{BezdekConnelly,Mallee};
the recent survey of the planar situation in \cite{BMP} identifies
Mall\'ee's parabolic triangle as the best known upper construction.
None of those upper constructions is needed in the proof above.
\end{remark}

\subsection{Finite witnesses and effective approximation errors}
\begin{proposition}[A circle-witness finite-net bound]\label{prop:circle-net}
Let a finite translation test family of closed curves, each of length
at most one, contain the circle of radius $1/(2\pi)$ and a $\delta$-net
of normalized triangles of perimeter at most one. If its minimum cover area is $W_\delta^\circ$, then
\begin{equation}\label{eq:circle-net}
 0\le\Wet-W_\delta^\circ
 \le 4\pi W_\delta^\circ\delta+\pi\delta^2.
\end{equation}
In particular, for any certified $V_+\ge\Wet$, both
$4\pi V_+\delta+\pi\delta^2$ and
$4\pi^2 V_+\delta+\pi\delta^2$ are valid explicit error bounds.
\end{proposition}
\begin{proof}
The finite family consists of allowable unit-length tests, so
$W_\delta^\circ\le\Wet$. Its circle witness forces an inball of
radius $\rho=1/(2\pi)$ in every convex cover and ensures attainment
by \cref{thm:outer-net}. Let $K_\delta$ be a minimizing cover with
that inball centered at zero. The parallel body
$K_\delta+\delta B_2^2$ covers every triangle, hence every closed
unit-length curve by \cref{prop:finite-witness}.

For any polygon containing $\rho B_2^2$, its area is one-half the
sum over edges of edge length times support distance from zero.
It is therefore at least $\rho\per(K)/2$. Approximation by convex
polygons, with convergence of area and perimeter, gives the same
inequality for a convex body. Thus
$\per(K_\delta)\le2W_\delta^\circ/\rho=4\pi W_\delta^\circ$.
Planar Steiner's formula now gives \eqref{eq:circle-net}.

For completeness, the coarser coefficient follows independently from
\eqref{eq:outer-radius}: $K_\delta\subset(2\pi V_+)B_2^2$.
Monotonicity of perimeter gives
$\per(K_\delta)\le2\pi(2\pi V_+)=4\pi^2V_+$; applying Steiner's
formula gives the second stated error. The inradius estimate below
allows a finite triangle family without the circle witness.
\end{proof}

\begin{corollary}[Improved finite-triangle approximation error]\label{cor:wetzel-net}
Take a finite Hausdorff $\delta$-net of normalized closed triangles
of perimeter at most one, allowing degenerate triangles, and retain
one nondegenerate triangle. Let $W_\delta$ be its minimum convex
translation-cover area. For $0<\delta<r_0$,
\begin{align}
 W_\delta&\le\Wet\le
       \left(\frac{r_0}{r_0-\delta}\right)^2W_\delta,
       \label{eq:wetzel-net-multiplicative}\\
 0\le\Wet-W_\delta
   &\le\frac{2W_\delta\delta}{r_0-\delta}+\pi\delta^2.
       \label{eq:wetzel-net-additive}
\end{align}
An independently certified upper bound $V_+$ for $\Wet$ may replace
$W_\delta$ on the right of \eqref{eq:wetzel-net-additive}.
\end{corollary}
\begin{proof}
By \cref{prop:finite-witness}, triangles and doubled segments are a
sufficient closed-curve test family. If $K_\delta$ minimizes the
finite family, $K_\delta+\delta B_2^2$ covers every triangle and
therefore every closed unit-length curve. Its inradius is at least
$r_0$ by \cref{lem:wetzel-inradius}.

For a convex body, the inradius has the support representation
\[
 r(K)=\max_z\min_{u\in\Sph^1}(h_K(u)-z\cdot u).
\]
Adding a disk adds its radius to every expression inside the minimum;
thus $r(K+\delta B_2^2)=r(K)+\delta$ and
$r(K_\delta)\ge r_0-\delta$.
Normalize the corresponding inball at zero. The inclusion
\[
 K_\delta+\delta B_2^2\subset
 \left(1+\frac{\delta}{r_0-\delta}\right)K_\delta
\]
proves \eqref{eq:wetzel-net-multiplicative}.
The support-area identity
$2\area(K)=\int h_K\,dS_K$ and the inball imply
$\per(K_\delta)\le2W_\delta/(r_0-\delta)$.
Steiner's planar formula
$\area(K+\delta B)=\area(K)+\delta\per(K)+\pi\delta^2$
then proves \eqref{eq:wetzel-net-additive}.
\end{proof}

\begin{corollary}[An explicit finite rational witness: an existence certificate]\label{cor:explicit-witness}
Put $L_0=f(r_0)+c_1$ and
\[
 \delta_* =r_0\left(1-\sqrt{f(r_0)/L_0}\right)
            >0.00000615580355820.
\]
Every finite triangle net of mesh $\delta<\delta_*$ as in
\cref{cor:wetzel-net} has minimum covering area strictly greater
than Wetzel's classical lower bound $f(r_0)$.

In particular, set $h=1/750000$ and let $\mathscr T_h$ consist
of all triangles with vertices $0,p,q$, where
\[
 p,q\in h\mathbb Z^2\cap[-1/2,1/2]^2,
 \qquad |p|+|q-p|+|q|\le1.
\]
This explicitly specified finite family of rational triangles has
minimum translation-cover area greater than $f(r_0)$.
\end{corollary}
\begin{proof}
Equation \eqref{eq:wetzel-net-multiplicative} and
\cref{thm:wetzel-improvement} give
\[
 W_\delta\ge(1-\delta/r_0)^2L_0>f(r_0)
 \quad\text{if }\delta<\delta_*.
\]
It remains to verify the asserted mesh. A triangle of perimeter at most
one rooted at zero has $|p|,|q|\le1/2$. Contract its two other vertices
by $1-4h$, then round both coordinates of each to the nearest
multiple of $h$. Each rounding error is at most $h/\sqrt2<h$.
The rounded vertices remain in $[-1/2,1/2]^2$. Perimeter increases
by at most twice the sum of the two rounding errors, which is at
most $4h$; hence the rounded triangle still has perimeter at most
$(1-4h)+4h=1$. Each vertex moves by at most $2h+h=3h$ from
its original value. Corresponding line segments show that the
triangle traces, and also their hulls, are within Hausdorff distance
$3h=1/250000$.

The family includes the nondegenerate triangle with vertices
$0,(1/4,0),(0,1/4)$. The rational verifier checks
$1/250000<\delta_*$ and, more specifically,
\[
 W_{1/250000}\ge(1-1/(250000r_0))^2L_0
             >0.15545139159107>f(r_0).
\]
This proves the finite-family claim. Before the perimeter filter,
the specification ranges over $(750001)^4>3\cdot10^{23}$ ordered
grid-point pairs. This is an existence and computability certificate;
the family has not been enumerated and is not a practical numerical
scheme. For rigorous optimization of a small selected family,
the geometric box-search strategy of Grechuk--Som-am
\cite{GrechukSomam} is a relevant computational comparison. Their
closed-curve problem allows congruent placement, so its area bound
is not a competing value of the translation-only $\Wet$.
\end{proof}

\begin{corollary}[Algebraic witnesses of exactly unit perimeter]
\label{cor:unit-algebraic}
There is an explicitly specified finite family of triangles with
real-algebraic vertices, allowing degeneracies, and perimeter exactly
one whose minimum convex
translation-cover area exceeds the rational number
\[
 q=\frac{15545139159107}{10^{14}}>f(r_0).
\]
\end{corollary}
\begin{proof}
Take the rational triangles $A=\conv\{0,p,q'\}$ in
\cref{cor:explicit-witness}. Omit the singleton, and write
$s_A=|p|+|q'-p|+|q'|\in(0,1]$, counting the boundary as the
three-edge closed path also in degenerate cases. Replace $A$ by
$\widehat A=s_A^{-1}A$. Each $s_A$ is real algebraic, so the new
vertices are real algebraic and the associated closed polygon has
perimeter exactly one. Since $0\in A$ and $s_A\le1$, convexity gives
$A\subset\widehat A$. Any convex body covering the new triangle
boundaries also covers their hulls and hence every original triangle.
Its area is therefore greater than $q$ by
\cref{cor:explicit-witness}. The omitted singleton places no additional
constraint on a nonempty covering body. The arithmetic certificate in
\cref{app:arithmetic} also verifies $q>f(r_0)$.
\end{proof}

\section{Minkowski bounds, symplectic ratios, and further consequences}\label{sec:consequences}
\subsection{Convex hulls, realization cost, and path structure}\label{sec:hull}
Support constraints see a path's convex hull, whereas the objective
counts the length needed to generate that hull. Keeping these two
roles separate prevents an invalid substitution of boundary length
for open traversal length.

\begin{proposition}[Support variation and exact area calculus]\label{prop:shape-calculus}
Let a strictly convex planar body have a $C^2$ support function $h$
with $h+h''>0$, and let $h_t=h+t\varphi$ for a $C^2$, $2\pi$-periodic
function $\varphi$. For sufficiently small $|t|$ these are support
functions of strictly convex bodies $K_t$, and
\begin{align}
 \area(K_t)&=\frac12\int_0^{2\pi}(h_t^2-h_t'^2)\,d\theta,
       \label{eq:shape-area}\\
 \left.\frac{d}{dt}\area(K_t)\right|_{t=0}
   &=\int_0^{2\pi}\varphi(h+h'')\,d\theta,\qquad
 \frac{d^2}{dt^2}\area(K_t)
   =\int_0^{2\pi}(\varphi^2-\varphi'^2)\,d\theta.
       \label{eq:shape-area-variation}
\end{align}
For a fixed compact object $A$, let $\Omega_*(A,K)$ denote all
maximizing balanced collections in \eqref{eq:margin-dual}. The
one-sided derivatives of its translation margin are
\begin{align}
 \left.\frac{d}{dt^+}m_{K_t}(A)\right|_{t=0}
  &=\max_{(u,\lambda)\in\Omega_*(A,K)}
                    \left(-\sum_i\lambda_i\varphi(u_i)\right),
       \label{eq:margin-variation}\\
 \left.\frac{d}{dt^-}m_{K_t}(A)\right|_{t=0}
  &=\min_{(u,\lambda)\in\Omega_*(A,K)}
                    \left(-\sum_i\lambda_i\varphi(u_i)\right).
       \nonumber
\end{align}
Here $\varphi(u_\theta)$ means $\varphi(\theta)$.
\end{proposition}
\begin{proof}
The positive minimum of $h+h''$ ensures positivity of
$h_t+h_t''$ for small $|t|$. Set
$x_t(\theta)=h_tu_\theta+h_t'u_\theta^\perp$ and
$\rho_t=h_t+h_t''>0$, so $x_t'=\rho_tu_\theta^\perp$.
For each fixed direction $u_\psi$,
\[
 \frac{d}{d\theta}\bigl(x_t(\theta)\cdot u_\psi\bigr)
       =\rho_t(\theta)\sin(\psi-\theta).
\]
This projection strictly decreases on $(\psi,\psi+\pi)$ and
strictly increases on $(\psi+\pi,\psi+2\pi)$. Its unique maximum
is therefore $h_t(\psi)$ at $\theta=\psi$ modulo $2\pi$.
The convex hull of the curve has support $h_t$, with exactly one
point in every supporting face. It is the required strictly convex
body, and $x_t$ parametrizes its boundary once counterclockwise.
The signed area formula gives
$\area(K_t)=\tfrac12\int h_t(h_t+h_t'')$, and integration by
parts gives \eqref{eq:shape-area}. Differentiation of this quadratic
expression proves \eqref{eq:shape-area-variation}.

For the margin, \cref{thm:balance} expresses $m_{K_t}(A)$ as a
maximum over a fixed compact parameter set of affine functions
\[
 G(u,\lambda)-t\sum_i\lambda_i\varphi(u_i).
\]
For $t>0$, each maximizer at zero gives a lower bound on the
difference quotient, so its limit inferior is at least the maximum
in \eqref{eq:margin-variation}. Choose maximizing collections at
$t\downarrow0$. Compactness yields cluster points, all maximizing
at zero, while optimality at zero bounds the corresponding difference
quotients above by $-\sum_i\lambda_i\varphi(u_i)$.
This gives the matching limit superior. Replacing $\varphi$ by
$-\varphi$ proves the left derivative formula.
\end{proof}

These formulas justify support-based shape variations and expose the
active balanced collections. They do not differentiate the globally
optimized escape threshold without an additional envelope argument
over paths and orientations. For example, whenever $E_H(K_t)$ is
actually differentiable and nonzero, the ordinary chain rule gives
\[
 \left.\frac{d}{dt}\frac{\area(K_t)}{E_H(K_t)^2}\right|_{t=0}
 =\frac{\area(K)}{E_H(K)^2}
   \left(\frac{\area'(K_t)|_0}{\area(K)}
                  -2\frac{E_H'(K_t)|_0}{E_H(K)}\right).
\]
Differentiability of this optimized threshold is an explicit premise.

For a compact convex set $A$ and a root $o\in A$, define
\[
 \tau_o(A)=\inf\{\len(\gamma):\gamma(0)=o,
                         \conv\gamma([0,1])=A\}.
\]
Set $\tau_o(A)=\infty$ if no such rectifiable path exists.

\begin{proposition}[Exact hull-realization formulas]\label{prop:hull}
Every compact path generating $A$ contains all extreme points of $A$.
The following consequences hold.
\begin{enumerate}[label=(\alph*)]
\item If $A$ is a polytope with vertices $v_1,\ldots,v_n$, then
\[
 \tau_o(A)=\min_{\sigma\in\mathfrak S_n}
       \left(|o-v_{\sigma(1)}|
         +\sum_{i=1}^{n-1}|v_{\sigma(i+1)}-v_{\sigma(i)}|\right).
\]
Thus realization is an exact rooted Hamiltonian-path problem.
\item If $A\subset\R^2$ is a strictly convex body, then
\[
 \tau_o(A)=\per(A)+\dist(o,\partial A).
\]
\item If $d\ge3$ and $A\subset\R^d$ is a strictly convex body,
then $\tau_o(A)=\infty$.
\end{enumerate}
\end{proposition}
\begin{proof}
If an extreme point $x$ belongs to the convex hull of a compact trace,
Carath\'eodory expresses it as a finite convex combination of trace
points. Extremality forces every point with positive weight to be $x$.
Thus $x$ lies on the trace.

For a polytope, order the first visits to its vertices and replace
intervening subpaths by straight segments. The triangle inequality
proves the lower bound in (a). Conversely each listed polygonal path
lies in $A$, visits every vertex, and generates $A$.

In (b) every boundary point is extreme. Up to its first boundary
contact, the path has length at least $\dist(o,\partial A)$.
Its remaining trace contains the whole boundary, since before that
contact the trace lies in the interior. Its length is therefore at
least $\HH^1(\partial A)=\per(A)$. A shortest segment to the
boundary followed by one full boundary traversal attains this bound.

In (c) the trace would contain $\partial A$. This boundary has
infinite one-dimensional Hausdorff measure: its orthogonal projection
onto a suitable $(d-1)$-dimensional hyperplane contains a ball with
positive $(d-1)$-dimensional volume, whose $\HH^1$ measure is
infinite when $d-1>1$. Lipschitz projections do not increase
$\HH^1$. A rectifiable path has finite $\HH^1$ measure, a
contradiction.
\end{proof}

\begin{corollary}[Exact shape--traversal separation]\label{cor:shape-traversal}
For the open rotational escape problem,
\[
 E_H(K)=\inf\left\{\tau_o(A):o\in A,
                   \min_{R\in H}s_K(A,R)\ge1\right\}.
\]
Here $A$ ranges over compact convex sets and $s_K$ is the fitting
scale of \cref{cor:fit}. This is an exact reformulation, not an
assertion that all admissible hulls have finite realization cost.
\end{corollary}
\begin{proof}
Convexity of $K$ makes fitting depend only on the convex hull.
Every path contributes its hull with no greater realization infimum.
Conversely, paths approaching $\tau_o(A)$ have the same non-fitting
hull and hence are feasible escape paths. Taking the two infima proves
the identity.
\end{proof}

\subsection{Restricted symplectic ratios}\label{sec:symplectic}
The symplectic connection uses an external theorem, not a consequence
of compactness alone. With the usual Ekeland--Hofer--Zehnder capacity
normalization, the polygonal non-fitting characterization gives
\begin{equation}\label{eq:external-capacity}
 c_{\rm EHZ}(K\times T)=E_T^{\rm tr}(K)
\end{equation}
for convex bodies $K,T\subset\R^d$ with $0\in\Int T$; this is
\cite[Theorem~1]{RudolfCapacity}, combined with the finite-polygon
reduction in \cref{thm:helly}. See also \cite{Rudolf,AkopyanEtAl}.
The cited theorem uses closed polygons that cannot be translated into
$\Int K$ and applies without smoothness or strict convexity assumptions. The symplectic form pairs the two Euclidean
factors. For a body $T$ whose interior does not contain zero,
translate it before applying the identity: a product translation is
symplectic, and translations of an admissible gauge do not change
closed-loop lengths. Thus all subsequent product ratios can be
normalized with both origins in the interiors. Equation
\eqref{eq:external-capacity} is the only capacity theorem used in
the new deductions below.

The reciprocal worm/capacity relation is established in
\cite{Rudolf}. We state it with the compactness argument for
attainment used here. In the plane, the known bounds are
$(3+\sqrt5)/5\le\kappa_2^{\rm Lag}<\sqrt2$ \cite{BMP}; attainment
does not determine this constant.

\begin{theorem}[Fixed-gauge and relaxed product ratios]\label{thm:capacity-ratio}
Define
\[
 \Sigma_d(T)=\sup_K
  \frac{c_{\rm EHZ}(K\times T)^d}{d!\vol(K)\vol(T)},
 \qquad \Psi_d(T)=d!\vol(T)W_T.
\]
Then
\begin{equation}\label{eq:capacity-ratio}
 \Sigma_d(T)=\frac1{\Psi_d(T)},\qquad
 \kappa_d^{\rm Lag}:=\sup_{K,T}
  \frac{c_{\rm EHZ}(K\times T)^d}{d!\vol(K)\vol(T)}
       =\frac1{\inf_T\Psi_d(T)}.
\end{equation}
For every fixed dimension these suprema are finite and attained.
In particular the restricted inequality
$c_{\rm EHZ}(K\times T)^d\le d!\vol(K)\vol(T)$ for all $K$
with fixed $T$ is equivalent to
\[
 W_T\ge\frac1{d!\vol(T)}.
\]
\end{theorem}
\begin{proof}
Substitute \eqref{eq:external-capacity} and take reciprocals in
\eqref{eq:minkowski-quotient}. Positivity permits this operation.
Taking the supremum in $T$ proves the second identity.
For fixed $T$, let $K_*$ be a minimum unit cover given by
\cref{cor:outer-families}. It has $E_T^{\rm tr}(K_*)=1$:
strictification gives at least one, and a larger value would allow
scaling down to a smaller universal cover. Thus it attains the
fixed-gauge ratio.

John's ellipsoid theorem \cite{John,Schneider} gives, after an affine
map and a translation, the normalization
$B_2^d\subset T\subset dB_2^d$. Such a transformation leaves
$\Psi_d$ unchanged by \cref{thm:gauge}. This normalized class is
Hausdorff compact, and $\Psi_d$ is continuous by \cref{cor:psi}.
Furthermore
\[
 \Psi_d(T)\ge d!\kappa_d\,d^{-d}W_{B_2^d}>0.
\]
Hence its positive infimum is attained by some $T_*$. Choosing a
minimum cover for $T_*$ attains the joint ratio and proves finiteness.
The equivalence in the last statement follows directly from
\eqref{eq:capacity-ratio}.
\end{proof}

\begin{lemma}[Volume forced by translated independent segments]\label{lem:independent-segments}
If a convex body $K\subset\R^d$ contains a translate of each segment
$[0,v_i]$, where $v_1,\ldots,v_d$ are linearly independent, then
\begin{equation}\label{eq:independent-segments}
 \vol(K)\ge\frac{|\det(v_1,\ldots,v_d)|}{d!}.
\end{equation}
Both the simplex $\conv(0,v_1,\ldots,v_d)$ and the centrally
symmetric crosspolytope $\conv\{\pm v_i/2:1\le i\le d\}$ attain
this bound and contain the required segment translates.
\end{lemma}
\begin{proof}
An invertible linear change of variables reduces to $v_i=e_i$.
Perform Steiner symmetrization successively along the coordinate
axes: replace each fiber parallel to an axis by the centered interval
of the same length. Convexity is preserved because fiber lengths are
concave on the projection; Fubini's theorem preserves volume.
Symmetrization along $e_i$ preserves a unit chord parallel to $e_i$.
It also preserves the existence of every unit chord parallel to
$e_j$, $j\ne i$: the projection of that chord lies in the projected
body, and the centers of the corresponding nonempty fibers form
the same unit chord in the hyperplane $x_i=0$.

The operation commutes with reflections in the other coordinate
hyperplanes. After all $d$ symmetrizations the body is invariant under
every coordinate reflection and still has a unit chord in each axis
direction. Averaging the reflections in the other coordinates moves
such a chord onto its coordinate axis; convexity and reflection in
the hyperplane $x_i=0$ then give $[-e_i/2,e_i/2]$ in the body. Their convex hull
has volume $2^d(1/2)^d/d!=1/d!$. Volume preservation and the
determinant change of variables prove \eqref{eq:independent-segments}.
The volumes and containment assertions for the two displayed examples
follow by triangulating into simplices.
\end{proof}

The next sharp value and product inequality are known in every
dimension by \cite[Theorem~5.2]{Balitskiy2020}. The simplex case
appears in its Theorem~5.1; see also the planar trapezoid results
of \cite{RudolfR4}. We retain an independent support-and-volume
proof as a model application of the complete test family.

\begin{theorem}[Exact worms and product inequalities for parallelotope gauges]\label{thm:parallelotope}
For any parallelotope gauge $T\subset\R^d$,
\begin{equation}\label{eq:parallelotope-worm}
 W_T=\frac1{d!\vol(T)},\qquad \Psi_d(T)=\Sigma_d(T)=1.
\end{equation}
Consequently every convex body $K$ satisfies
\begin{equation}\label{eq:parallelotope-viterbo}
 c_{\rm EHZ}(K\times T)^d\le d!\vol(K)\vol(T).
\end{equation}
For $T=[-1,1]^d$, both
\[
 \conv\{0,e_1/2,\ldots,e_d/2\}
 \quad\text{and}\quad
 \conv\{\pm e_i/4:1\le i\le d\}
\]
are minimum covers.
\end{theorem}
\begin{proof}
First take $T=[-1,1]^d$. Its only extreme positive facet circuits
are the $d$ antipodal pairs $(e_i,-e_i)$, each with offsets one.
\Cref{thm:complete-test-family} says that universal coverage is
equivalent to containing a translate of each segment $[0,e_i/2]$.
The lemma gives $\vol(K)\ge1/(2^d d!)$, and its two equality
examples give the displayed minimum covers. Since $\vol(T)=2^d$,
this proves \eqref{eq:parallelotope-worm} for the cube.
For $T=A[-1,1]^d+z$, affine gauge covariance gives
$W_T=|\det A|^{-1}/(2^d d!)=1/(d!\vol(T))$.
\Cref{thm:capacity-ratio} yields the product inequality and its sharpness.

\end{proof}

These sharp model cases are recovered here from finite complete tests
and a volume argument. They provide unconditional applications of the
Minkowski framework in all dimensions. They do not extend the product
inequality to every gauge. Broader established classes and the general symplectic
connection are discussed in \cite{Rudolf,ArtsteinKarasevOstrover,BMP}.

\begin{corollary}[Reference-gauge transfer and an ellipsoidal product bound]\label{cor:reference-capacity}
Suppose $W_S\ge L_S>0$ for a reference gauge $S$. For every gauge $T$,
\begin{equation}\label{eq:reference-capacity}
 W_T\ge\frac{L_S}{D(T,S)^d},\qquad
 \Sigma_d(T)\le\frac{D(T,S)^d}{d!\vol(T)L_S}.
\end{equation}
In dimension two, write $L_0=f(r_0)+c_1$ as in
\cref{thm:wetzel-improvement}. For every ellipse $T$ and every convex
body $K\subset\R^2$,
\begin{equation}\label{eq:ellipse-viterbo}
 \frac{c_{\rm EHZ}(K\times T)^2}{2\area(K)\area(T)}
 \le\frac1{2\pi L_0}
 <1.023781939711102.
\end{equation}
More generally, if a translate of $T$ is contained in an ellipse $E$,
\begin{equation}\label{eq:ellipsoid-transfer}
 W_T\ge\frac{\pi L_0}{\area(E)},\qquad
 \Sigma_2(T)\le\frac{\area(E)/\area(T)}{2\pi L_0}.
\end{equation}
One may minimize the last right-hand side over circumscribed ellipses.
\end{corollary}
\begin{proof}
Apply \eqref{eq:translated-worm} with the gauges interchanged and
then \eqref{eq:capacity-ratio}. If $T=AB_2^2+z$, affine covariance gives
$W_T=|\det A|^{-1}\Wet$, while $\area(T)=\pi|\det A|$.
Thus $\Psi_2(T)=2\pi\Wet\ge2\pi L_0$, proving
\eqref{eq:ellipse-viterbo}. Finally, $T\subset E$ up to translation
implies $W_T\ge W_E$ by the gauge comparison. The same affine formula
for $E$ proves \eqref{eq:ellipsoid-transfer}.
\end{proof}

For comparison, using the classical lower bound alone would give the
weaker factor
\[
 \frac1{2\pi f(r_0)}=1.0238474372425079\ldots.
\]
The decrease exceeds $0.00006549753140655$. Both factors and their
difference are enclosed by the exact interval verifier in
Appendix~\ref{app:arithmetic}. Thus the quantitative Wetzel improvement
also strengthens an unconditional volume--capacity estimate for every
ellipsoidal planar gauge. It does not replace its factor by the
conjectural value one.

For a fixed pair $(K,T)$, \cref{thm:circuit-radius} gives an additional
exact criterion: the product inequality with a factor $A>0$,
\[
 E_T^{\rm tr}(K)^d\le A\,d!\vol(K)\vol(T),
\]
holds if and only if some ordered balanced list satisfies
\begin{equation}\label{eq:viterbo-circuit-target}
 C_T(a)\ge
 \frac{B_K(a)}{(A\,d!\vol(K)\vol(T))^{1/d}}.
\end{equation}
For polytopal $K$, one can restrict to its finitely many facet circuits
and orders. This is a finite geometric target for a particular product;
the universal assertion with $A=1$ cannot be inferred from the
existence of the formulation.

The symmetric diagonal identity is
due to Artstein-Avidan--Karasev--Ostrover
\cite{ArtsteinKarasevOstrover}; see also \cite{AkopyanEtAl}.
The elementary centroid proof below follows the closed-curve
covering argument of Chakerian--Klamkin \cite{ChakerianKlamkin},
whose relation to the billiard proof is noted in \cite{BMP}.

\begin{proposition}[Known symmetric diagonal identity: centroid proof]\label{prop:diagonal}
If $K=-K$ is a convex body, then
\begin{equation}\label{eq:diagonal}
 E_{K^\circ}^{\rm tr}(K)=4.
\end{equation}
\end{proposition}
\begin{proof}
Here $h_{K^\circ}=\mu_K$ is a norm. Let a closed curve of norm
length $L>0$ be parametrized by norm arclength $s\in[0,L]$, and
let $c=L^{-1}\int_0^Lq(s)\,ds$. For fixed $t$, the shorter of
the two curve arcs bounds the distance:
\[
 \|q(t)-c\|_K
 \le\frac1L\int_0^L\|q(t)-q(s)\|_K\,ds
 \le\frac1L\int_0^L\min\{|s-t|,L-|s-t|\}\,ds
 =\frac L4.
\]
Thus every loop of length less than four translates into $\Int K$.
The doubled diameter $-u\to u\to-u$, $u\in\partial K$, has
length four and cannot translate into $\Int K$, since its width in
a supporting direction at $u$ equals the width of $K$. This proves
\eqref{eq:diagonal} without capacity theory.

The quotient identity gives
$W_{K^\circ}\le\vol(K)/4^d$.
\end{proof}

\begin{remark}[Logical status]
The unrestricted Viterbo conjecture is false \cite{HKO}; a theorem
about all convex symplectic domains must not be inferred from
\eqref{eq:capacity-ratio}. The restricted product constants remain meaningful.
\cref{thm:capacity-ratio} proves existence of an optimal relaxed
ratio but does not calculate that ratio. Recent results for planar
polygonal gauges in \cite{BMP} are compatible with this framework
and should be used with their stated gauge and equality hypotheses.
\end{remark}

The bound $2+2/d$ and its simplex sharpness are established in
\cite{AkopyanEtAl}. The next calculation keeps the origin-dependent
asymmetry parameter in the centroid argument
\cite{ChakerianKlamkin}. We use it as a quantitative refinement
within that established method, without a priority claim.

\begin{proposition}[An asymmetric centroid bound and sharp simplex calculation]\label{prop:asymmetric-centroid}
Let $U$ be a convex body with $0\in\Int U$, and let
$s_o(U)=\min\{s\ge1:-U\subset sU\}$, with this specified origin.
Every closed curve of $\mu_U$-length $L$ has a translate contained in
\begin{equation}\label{eq:asymmetric-centroid}
 \frac{s_o(U)}{2(1+s_o(U))}\,L U.
\end{equation}
Consequently,
\begin{equation}\label{eq:asymmetric-diagonal}
 E_{U^\circ}^{\rm tr}(U)\ge2+\frac2{s_o(U)},\qquad
 W_{U^\circ}\le
 \left(\frac{s_o(U)}{2(1+s_o(U))}\right)^d\vol(U).
\end{equation}
For a simplex whose centroid is zero, $s_o(U)=d$ and
\begin{equation}\label{eq:simplex-diagonal}
 E_{U^\circ}^{\rm tr}(U)=2+2/d.
\end{equation}
\end{proposition}
\begin{proof}
The assertion for $L=0$ is immediate. Otherwise use oriented
$\mu_U$-arclength on $[0,L]$ and put
$c=L^{-1}\int_0^Lq(s)\,ds$. Since
$\mu_U(-v)\le s_o(U)\mu_U(v)$, the two oriented arcs from $q(s)$
to $q(t)$ give
\[
 \mu_U(q(t)-q(s))\le\min\{r,s_o(U)(L-r)\},
 \qquad r=(t-s)\pmod L.
\]
Convexity and integration over a full period imply
\[
 \mu_U(q(t)-c)
 \le\frac1L\int_0^L\min\{r,s_o(U)(L-r)\}\,dr
 =\frac{s_o(U)L}{2(1+s_o(U))}.
\]
Splitting the integral at $r=s_o(U)L/(1+s_o(U))$ proves the equality.
This gives the containment and both consequences in
\eqref{eq:asymmetric-diagonal}.

For the exact simplex value, linear covariance reduces the calculation
to $U=\Delta^\circ$, where
$\Delta=\conv\{e_1,\ldots,e_d,-\mathbf1\}$.
Its facet normals are $n_i=e_i$ and $n_{d+1}=-\mathbf1$, all with
offset one, and their equal-weight circuit has $B=d+1$.
For the order $e_1,\ldots,e_d,-\mathbf1$, the state set is
\[
 S=\{0,-e_1,-e_1-e_2,\ldots,-\mathbf1\}.
\]
The facets of $\Delta$ have offset one and normals
$f_0=\mathbf1$, $f_i=\mathbf1-(d+1)e_i$ for $1\le i\le d$.
Their unique positive circuit has equal weights. The fitting formula
therefore yields
\[
 \operatorname{rad}_{\Delta}(S)
 =\frac{\sum_{i=0}^d h_S(f_i)}{d+1}
 =\frac{0+\sum_{i=1}^d(d+1-i)}{d+1}=\frac d2.
\]
The center and the matching certificate are explicit. Put
$z_i=-1+i/(d+1)$ for $1\le i\le d$. Then
\[
 \sum_i z_i=-d/2,\qquad
 f_i\cdot z=h_S(f_i)-d/2\quad(0\le i\le d).
\]
Consequently $S\subset z+(d/2)\Delta$. The equal dual weights
$1/(d+1)$ on these supporting facets have zero weighted normal and
unit weighted offset, and their weighted support sum is $d/2$.
This proves the same radius by matching primal and dual certificates.

One can also display an attaining escape polygon directly. Set
\[
 q_k=\frac{2k}{d}\mathbf1-\frac{2(d+1)}d\sum_{i=1}^k e_i
       \quad(0\le k\le d),\qquad q_{d+1}=q_0=0.
\]
Its successive edges are $2f_1/d,\ldots,2f_d/d,2f_0/d$.
Since $h_\Delta(f_i)=1$, its $\Delta$-length is $2(d+1)/d$.
The balanced contacts with the facets of $U$ give
\[
 \sum_{i=1}^d e_i\cdot q_{i-1}-\mathbf1\cdot q_d
       =\sum_{i=1}^d\frac{2(i-1)}d+2=d+1.
\]
A strict translate into $U$ would make this sum smaller than the
sum $d+1$ of the facet offsets, a contradiction. Thus the polygon
attains the upper bound $2+2/d$ without any existence inference.
For a centroid-centered simplex, each opposite vertex has signed
facet height $-d$ when that facet has height one; hence $s_o(U)=d$.
The already proved lower bound matches this explicit upper bound.
\end{proof}

The sharp dimension-only diagonal estimate $2+2/d$ is known for
arbitrary convex bodies \cite[Theorem~4.1]{AkopyanEtAl};
\eqref{eq:simplex-diagonal} is included as an explicit check of the
circuit formula, not as a new simplex theorem. The centroid estimate
\eqref{eq:asymmetric-diagonal} retains the actual origin-dependent
asymmetry. When that origin is a Minkowski center, meaning that
$s_o(U)=\sigma(U)$, it interpolates
between the symmetric value four and the simplex value $2+2/d$,
and is stronger than the dimension-only lower bound if $s_o(U)<d$.
Changing this origin changes $U^\circ$ and must not be done silently.

\begin{corollary}[Finite certificates for a strict ratio improvement]\label{cor:ratio-certificate}
If a pair of convex bodies has a volume--capacity ratio strictly
larger than a specified number, then a pair of rational polytopes
has the same strict inequality. Its translative escape threshold
admits finite rational LP certificates. Conversely, failure to find
such a certificate in a finite search is not an upper bound on the
continuum supremum.
\end{corollary}
\begin{proof}
First translate both bodies so that their interiors contain zero;
volume and the product capacity are unchanged. Approximate the
translated bodies by rational polytopes whose interiors still contain
zero. Joint continuity from \cref{cor:psi} and volume
continuity preserve a strict inequality for all sufficiently close
approximants. Apply \cref{thm:lp} to the selected pair.
\end{proof}

\subsection{Homogeneous spaces and curvature rescaling}\label{sec:variants}
The metric arguments extend beyond Euclidean convexity. The quotient
arguments do not automatically do so: they need a scaling operation
within the same problem class and an endpoint strictification theorem.

\begin{remark}[Restricted motions and degeneracy under unrestricted scaling]
If $H_0\subset H_1\le\OO(d)$, then
$E_{H_0}(K)\le E_{H_1}(K)$ and $D_{H_0}(K)\le D_{H_1}(K)$:
non-fitting against more placements is a stronger requirement.
The corresponding universal-cover constants have the reverse
inequality, because more motions make covering easier. These facts
follow directly from their defining feasible families.
If the allowed placements instead include arbitrary contractions,
every compact object fits inside every body with nonempty interior;
the escape threshold becomes infinite and the unnormalized covering
volume infimum becomes zero. General affine and similarity variants
therefore require an explicit size or determinant normalization.
\end{remark}

\begin{proposition}[Homogeneous metric escape]\label{prop:homogeneous}
Let a locally compact metrizable group act continuously, properly,
transitively, and isometrically on a proper geodesic metric space $X$.
Let $F$ be compact, and suppose the open-path escape threshold is
positive and finite. Then it is attained, its critical paths have
closed placements in $F$, and the corresponding inverse-boundary
TSPN problems on nested dense finite configuration sets converge with
recovery of minimizing paths.
\end{proposition}
\begin{proof}
Transitivity normalizes the initial point of every candidate path
to $o$. Paths in a minimizing sequence have uniformly bounded length,
stay in a compact ball, and admit a common Lipschitz bound after
arclength normalization. Uniform limits preserve non-fitting:
a compact trace with a strict placement in $F$ has positive distance
from $X\setminus\Int F$, and close traces share that placement.
Length is lower semicontinuous, so a minimizer exists. Initial
truncations of a critical open path approach its trace, giving closed
endpoint placement by \cref{thm:endpoint}. Properness makes the
configuration set compact, and inverse boundaries have a closed graph.
The hitting interpretation \eqref{eq:reversal} and
\cref{thm:metric} give convergence and recovery.
\end{proof}

For a fixed class of compact forests in a Riemannian homogeneous space
and $L>0$, define $W^{\rm cl}(L)$ using closed placements of all open paths of
length at most $L$. Define the robust interior value by
\[
 W^{\rm rob}(L)=\inf\{\vol(F):
   \text{every path of length }<L\text{ fits in }\Int F\}.
\]
Under the hypotheses providing endpoint placement,
\begin{equation}\label{eq:curved-robust}
 W^{\rm rob}(L)=\inf_{E(F)\ge L}\vol(F),
 \qquad W^{\rm cl}(L)\le W^{\rm rob}(L).
\end{equation}
The equality is the threshold identity; the inequality follows from
endpoint placement. Equality of the two covering values requires
an additional strictification argument for the chosen forest class.

\begin{proposition}[Constant-curvature scaling]\label{prop:curvature}
Let $d\ge2$, and let $W_{\kappa,d}(L)$ denote either covering value above in the
complete simply connected $d$-dimensional space form of sectional
curvature $\kappa$, with a forest class preserved by metric rescaling.
Then, for $L>0$,
\begin{equation}\label{eq:curvature-scaling}
 W_{\kappa,d}(L)=L^dW_{\kappa L^2,d}(1).
\end{equation}
The identity also holds with infinite values. It does not assert
homogeneity in $L$ while keeping nonzero curvature fixed.
\end{proposition}
\begin{proof}
Replace the metric tensor $g$ by $L^{-2}g$. Lengths are divided by
$L$, volumes by $L^d$, and sectional curvature is multiplied by
$L^2$. Isometries, containment, and the chosen forest class correspond
under this common metric rescaling. Unit covers in the rescaled
problem are therefore exactly length-$L$ covers in the original
problem, with the displayed volume factor.
\end{proof}
In particular, $\kappa=0$ gives the Euclidean law
$W_{0,d}(L)=L^dW_{0,d}(1)$, since $0\cdot L^2=0$.

\subsection{Cooperating paths and weighted arrival costs}\label{sec:search-variants}
\begin{proposition}[A finite number of cooperating paths]\label{prop:multi}
Under the target hypotheses of \cref{thm:metric}, fix roots
$o_1,\ldots,o_k$, with $1\le k<\infty$. Require that the union of $k$
paths meet every target, and minimize either the sum or the maximum
of their lengths. If a finite feasible tuple exists, the optimum is
attained. Nested dense finite-target optima increase to it, and
normalized minimizing tuples are precompact with optimal cluster
points. A finite-target problem is exactly a finite minimum over
target assignments to paths, visitation orders, and points in the
assigned targets, with geodesic distances as edge costs.
\end{proposition}
\begin{proof}
A bound on either objective bounds each path length. Normalize the
paths separately on $[0,1]$ and use the product of the compactness
arguments in \cref{thm:metric}. For any converging sequence of target
parameters, select a hitting path index and a hitting time. A
subsequence has constant index and convergent times; the closed graph
passes the hit to the limit tuple. This proves closed feasibility,
attainment, and the dense-limit argument. The objective is lower
semicontinuous for both sum and maximum. For finitely many targets,
assign each target to one visiting path and order the chosen visits.
Shortcutting by geodesics cannot increase either objective, and
concatenating geodesics realizes every such discrete candidate.
\end{proof}

The min-mean Bellman's lost-in-a-forest problem asks for expected escape time;
early work includes Gluss \cite{Gluss} and Finch--Shonder
\cite{FinchShonder}. Related average disk-inspection optimization
is treated by Georgiou \cite{Georgiou}; its information model
differs from an unknown initial position in a disk. The following
standard finite-target ordering argument gives a discrete skeleton.
It does not solve the continuum expected-exit problem.

\begin{proposition}[Finite weighted first-arrival search]\label{prop:latency}
In a proper geodesic space, let $N_1,\ldots,N_m$ be nonempty closed
targets and $p_i>0$ with $\sum_i p_i=1$. A path begins at $o$ and
travels at speed at most one (that is, it is $1$-Lipschitz in time);
let $T_i$ be its first arrival time at $N_i$. Then
\begin{equation}\label{eq:latency}
 \inf_\gamma\sum_i p_iT_i
 =\min_{\substack{\sigma\in\mathfrak S_m\\x_j\in N_{\sigma(j)}}}
   \sum_{j=1}^m p_{\sigma(j)}
          \sum_{l=1}^j d(x_{l-1},x_l),\qquad x_0=o.
\end{equation}
Both sides attain their values. An incidental earlier visit along a
geodesic can only improve the realized first-arrival objective.
\end{proposition}
\begin{proof}
Order a feasible path's actual first arrivals, breaking ties
arbitrarily, and choose the first-hit points. Distances between
successive points do not exceed the intervening travel times, so
the expression on the right for this choice is at most the path's
weighted arrival cost. Conversely, concatenate geodesics through
any right-hand candidate. The actual arrival at each target occurs
no later than its scheduled visit, giving the opposite inequality
between infima.

For existence, fix a finite feasible tour and its finite cost $C$.
Every edge length on a candidate of cost at most $C$ is included in
the arrival time of the final target, whose weight is at least
$p_{\min}=\min_i p_i>0$. Hence the sum of edge lengths is at most
$C/p_{\min}$. The selected points lie in a common compact ball;
the targets are closed and the cost is continuous. There are finitely
many permutations, so the right side attains a minimum. Its geodesic
realization attains the common infimum on the left.
\end{proof}

\begin{remark}[What requires additional hypotheses]
The weighted finite theorem does not itself justify quadrature
convergence for a continuum probability law: first-arrival times
can be discontinuous. Likewise, a control-system extension needs a
specified proper geodesic travel metric, or a replacement compactness
theorem for admissible trajectories. Moving forests and feedback
policies are problems on trajectories or information states, not
automatically the same static-placement problem. These distinctions
identify additional assumptions needed for a theorem, rather than
replacing them with a formal analogy.
\end{remark}

\section{Consequences}\label{sec:outlook}
The facet-offset method computes the best cover for each prescribed
normal set and bounds the unrestricted fixed-gauge optimum with an
explicit mesh error. Its four small exact benchmarks recover known
sharp cases. The centered ledger and gauge comparisons provide the
support calculations needed by this method; their finite formulations
alone do not determine any unresolved extremal value.

The Wetzel improvement uses additional geometric information: a
uniform width deficit with one tangent support, followed by a
quantified area bound outside the comparison body. The finite-net
estimate uses the same inradius information to control outer
approximation error.

\appendix
\section{Dependency map of the main theorems}\label{app:dependencies}
\begin{figure}[htbp]
\centering
\begin{tikzpicture}[
  box/.style={draw=blue!45!black,rounded corners=2pt,
    fill=blue!3,align=center,text width=3.65cm,
    minimum height=1.0cm,inner sep=4pt,font=\small},
  arr/.style={-{Latex[length=2mm]},draw=black!65,semithick}]
\node[box] (threshold) at (0,0)
  {Threshold and endpoint\\\cref{thm:threshold,thm:endpoint}};
\node[box] (supports) at (4.8,0)
  {Balanced supports\\\cref{thm:balance}};
\node[box] (outer) at (9.6,0)
  {Outer compactness\\\cref{thm:outer-net}};
\node[box] (quotient) at (0,-1.75)
  {Homogeneous quotients\\\cref{thm:quotient,thm:minkowski}};
\node[box] (polygon) at (4.8,-1.75)
  {Finite polygon\\reduction\\\cref{thm:helly}};
\node[box] (witness) at (9.6,-1.75)
  {Finite witness\\reduction\\\cref{prop:finite-witness}};
\node[box] (stability) at (0,-3.5)
  {Body and gauge stability\\\cref{thm:body-stability,thm:gauge}};
\node[box] (lp) at (4.8,-3.5)
  {Exact rational LPs\\\cref{thm:lp}};
\node[box] (hierarchy) at (9.6,-3.5)
  {Two-sided computation\\\cref{thm:hierarchy}};
\draw[arr] (threshold) -- (quotient);
\draw[arr] (supports) -- (polygon);
\draw[arr] (polygon) -- (lp);
\draw[arr] (outer) -- (witness);
\draw[arr] (polygon) -- (witness);
\draw[arr] (witness) -- (hierarchy);
\draw[arr] (lp) -- (hierarchy);
\draw[arr] (stability.south) -- ++(0,-0.45) -| (hierarchy.south);
\draw[arr] (quotient.west) -- ++(-0.55,0) |- (stability.west);
\end{tikzpicture}
\caption{The main duality and finite-certification dependencies.
The quotient input to the stability box concerns the universal-cover
constant $W_T$; body stability itself follows directly from inclusion
and homogeneity. The hypotheses and further branches are detailed below.}
\label{fig:dependencies}
\end{figure}
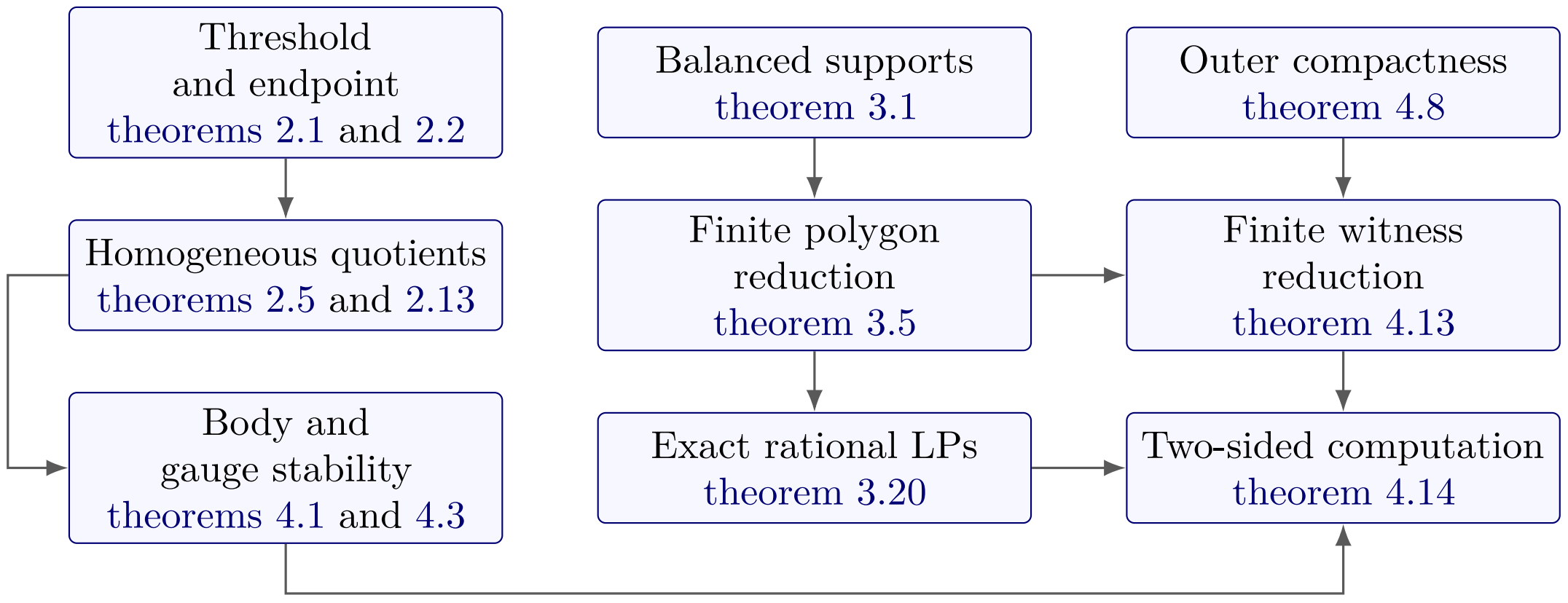

The arrows in \cref{fig:dependencies} indicate the principal inputs to
the displayed conclusions. The table gives the additional hypotheses
and covers the extensions omitted from the compact diagram. Classical
tools such as separation, Blaschke selection, and Steiner's formula
are identified explicitly. The geometric hypotheses are part of the
statements and must be checked before each theorem is applied.

\small
\begin{longtable}{@{}>{\raggedright\arraybackslash}p{0.26\textwidth}
  >{\raggedright\arraybackslash}p{0.34\textwidth}
  >{\raggedright\arraybackslash}p{0.34\textwidth}@{}}
\toprule
Conclusion & Inputs established in this article & Additional hypotheses or tools\\
\midrule
\endfirsthead
\toprule
Conclusion & Inputs established in this article & Additional hypotheses or tools\\
\midrule
\endhead
Threshold and cofinality
(\cref{thm:threshold,lem:cofinal}) &
Definitions of non-fitting size and greatest lower bound &
Order completeness; no convexity, topology, or compactness.\\[5pt]
Closed endpoint placement
(\cref{thm:endpoint,cor:critical}) &
Subcritical fitting from \cref{thm:threshold} &
Proper action, compact forest, and approximants within the object class.\\[5pt]
Metric TSPN and recovery
(\cref{thm:metric,cor:gamma}) &
Motion reversal and compact states in \cref{prop:motion} for forest targets &
Proper geodesic space, closed target graph, compact parameter space,
and one finite feasible path.\\[5pt]
Euclidean and Minkowski quotients
(\cref{thm:euclidean,thm:minkowski}) &
\Cref{thm:quotient,lem:strictification,lem:convex-hull};
\cref{thm:helly} for the Minkowski polygon assertion &
Homogeneity, positive finite thresholds, compact motions,
and convex covers.\\[5pt]
Balanced and finite support certificates
(\cref{thm:balance,prop:poly-circuit,thm:helly}) &
Signed support margin and its active-normal balance &
Finite-dimensional separation and Carath\'eodory; Farkas extreme rays
for polytopes.\\[5pt]
Ledger and calibration
(\cref{thm:normed-ledger,thm:ledger,cor:calibration}) &
Temporal telescoping and balanced contact functionals &
Dual-norm inequality; closure and finite-dimensional separation for
the optimal centered constant.\\[5pt]
Rotational certificates
(\cref{thm:rotation-rate,prop:finite-rot,prop:adaptive}) &
\Cref{thm:balance,lem:orientation-lip,thm:helly} &
Positive inradius, compact orientation group, and certified global
finite optimization.\\[5pt]
Diameter mesh recovery
(\cref{prop:diameter-hit,cor:diameter-mesh}) &
\Cref{prop:motion,lem:convex-hull}; exact thickening identity
\eqref{eq:diameter-thickening-exact} &
Compact Hausdorff-Lipschitz targets; an anchor for the forest interpretation.\\[5pt]
Exact translative LPs
(\cref{thm:lp}) &
\Cref{prop:poly-circuit,thm:helly,thm:minkowski} &
Polytopal body and gauge; rational data for rational primal, dual,
and infeasibility certificates.\\[5pt]
Stability and affine invariance
(\cref{thm:body-stability,thm:gauge,cor:psi}) &
Support-function inclusions; \cref{thm:minkowski} for $W_T$ &
Full-dimensional bodies and gauges; determinant change of variables.\\[5pt]
Outer attainment and approximation
(\cref{thm:outer-net,cor:steiner-error}) &
Inball witness, cone radius bound, and closed placement conditions &
Compact test family and motions; Blaschke selection, volume continuity,
and Steiner's formula for the additive estimate.\\[5pt]
Two-sided hierarchy
(\cref{thm:hierarchy}) &
\Cref{thm:outer-net,prop:finite-witness,cor:density,thm:lp} &
Rational polytopal gauge, exhaustive candidate families, and
real-algebraic decision procedures.\\[5pt]
Strict Wetzel lower bound
(\cref{prop:wetzel-strict}) &
\Cref{cor:outer-families,lem:wetzel-inradius,lem:wetzel-caps,lem:wetzel-defect} &
Closed planar loops under translations; strict volume monotonicity
for convex bodies.\\[5pt]
Explicit Wetzel increment
(\cref{thm:wetzel-improvement}) &
\Cref{lem:wetzel-inradius,lem:wetzel-caps,lem:wetzel-defect,lem:cap-area,lem:tangent-deficit} &
Tangency of one deficient support; the geometric inequality is analytic. Appendix~\ref{app:arithmetic}
certifies only its displayed decimals. Attainment is not needed for this
pointwise area bound.\\[5pt]
Finite strict Wetzel witnesses
(\cref{cor:explicit-witness,cor:unit-algebraic}) &
\Cref{thm:wetzel-improvement,cor:wetzel-net}; explicit grid and
algebraic perimeter normalization &
Finite triangle reduction and exact rational bounds for the numerical
separation from the classical lower bound.\\[5pt]
Hull realization
(\cref{prop:hull,cor:shape-traversal}) &
Extreme points in the trace; planar perimeter monotonicity &
Compact traces; rectifiable curves and traversal conventions.\\[5pt]
Attained capacity ratios
(\cref{thm:capacity-ratio}) &
\Cref{thm:minkowski,thm:gauge,cor:psi,cor:outer-families} &
The published capacity identity \eqref{eq:external-capacity} and
John's ellipsoid theorem.\\[5pt]
Exact centered circuit formula
(\cref{thm:circuit-radius,cor:circuit-lp}) &
\Cref{thm:balance,prop:poly-circuit,thm:ledger} &
Ordered balanced contacts; closure; circumradius in the fixed gauge.
Only $d+1$ variables per polytopal circumradius LP.\\[5pt]
Complete tests and polytopal minimizers
(\cref{thm:complete-test-family,cor:polytopal-minimizer}) &
\Cref{cor:circuit-lp,cor:escape-reciprocity,cor:outer-families} &
Extreme positive facet circuits; reverse traversal; convex hulls of
placed tests. Algebraicity additionally uses real-algebraic elimination.\\[5pt]
Exact facet-offset minimum and rational brackets
(\cref{thm:fixed-normals,cor:facet-hierarchy}) &
\Cref{cor:circuit-lp,lem:positive-concave,thm:outer-net}; the radius bound
\eqref{eq:outer-radius} &
Rational polytopal gauge, positive spanning normals, Brunn--Minkowski,
and finite polyhedral decomposition after removing translations.\\[5pt]
Exact facet benchmarks
(\cref{prop:facet-benchmarks}) &
\Cref{cor:circuit-lp,thm:fixed-normals}; exact determinants &
Known global model inequalities \cite{Balitskiy2020,BMP} for the final
unrestricted optimality assertion.\\[5pt]
Optimal gauge comparison and symmetrization
(\cref{thm:translated-gauges,prop:gauge-symmetrization}) &
\Cref{cor:fit,thm:minkowski}; balanced closed edge sequences &
Parametrized closed curves; reversal retained explicitly.\\[5pt]
Asymmetric centroid and simplex value
(\cref{prop:asymmetric-centroid}) &
Convexity, oriented arclength, and \cref{cor:circuit-lp} &
Origin-dependent asymmetry; a centroid-centered simplex for equality.\\[5pt]
Sharp parallelotope worms and products
(\cref{lem:independent-segments,thm:parallelotope}) &
\Cref{thm:complete-test-family,thm:gauge,thm:capacity-ratio} &
Antipodal cube circuits, coordinate Steiner symmetrization,
affine covariance, and the published capacity identity.\\[5pt]
Ellipsoidal capacity
(\cref{cor:reference-capacity}) &
\Cref{thm:wetzel-improvement,thm:gauge,thm:capacity-ratio} &
Published capacity identity.\\[5pt]
Curved and cooperative variants
(\cref{prop:homogeneous,prop:curvature,prop:multi,prop:latency}) &
\Cref{prop:motion,thm:endpoint,thm:metric} for hitting and recovery &
Transitivity for normalization; metric rescaling for curvature;
finitely many paths or finitely many positive target weights.\\
\bottomrule
\end{longtable}
\normalsize
\begin{figure}[htbp]
\centering
\begin{tikzpicture}[
 box/.style={draw=blue!45!black,rounded corners=2pt,fill=blue!3,
   align=center,text width=4.1cm,minimum height=1.0cm,
   inner sep=4pt,font=\small},
 arr/.style={-{Latex[length=2mm]},draw=black!65,semithick}]
\node[box] (support) at (-3.1,0)
 {Balanced supports\\\cref{thm:balance}};
\node[box] (ledger) at (3.1,0)
 {Sharp centered ledger\\\cref{thm:ledger}};
\node[box] (circuit) at (0,-1.75)
 {Exact circuit--circumradius formula\\\cref{thm:circuit-radius}};
\node[box] (offset) at (-3.1,-3.65)
 {Facet-offset polyhedron\\\cref{thm:fixed-normals}};
\node[box] (tests) at (3.1,-3.65)
 {Finite complete polygon tests\\\cref{thm:complete-test-family}};
\node[box] (bracket) at (-3.1,-5.45)
 {Rational worm brackets\\\cref{cor:facet-hierarchy}};
\node[box] (capacity) at (3.1,-5.45)
 {Sharp classes and ratio criteria\\\cref{thm:parallelotope}};
\draw[arr] (support) -- (circuit);
\draw[arr] (ledger) -- (circuit);
\draw[arr] (circuit) -- (offset);
\draw[arr] (circuit) -- (tests);
\draw[arr] (offset) -- (bracket);
\draw[arr] (tests) -- (capacity);
\draw[arr] (bracket) -- (capacity);
\end{tikzpicture}
\caption{The Minkowski proof chain. The last box also uses the
escape quotient and the published capacity identity. The parallelotope
case uses the independent-segments lemma.}
\label{fig:minkowski-dependencies}
\end{figure}
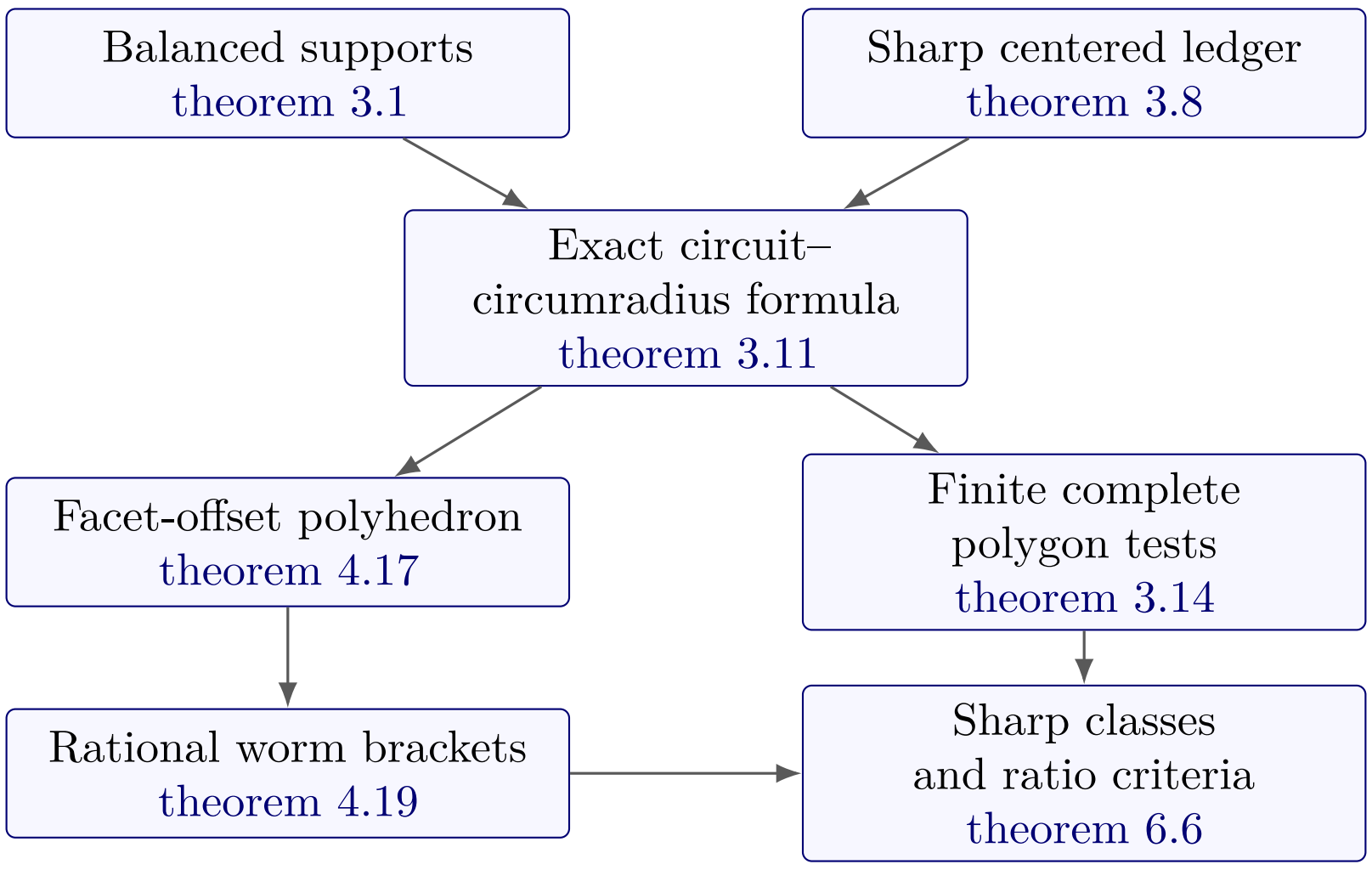

\section{Exact arithmetic verification and reproducibility}\label{app:arithmetic}
The accompanying standard-library Python file
\texttt{verify\_wetzel\_bound.py} checks the decimal inequalities
in \cref{thm:wetzel-improvement,cor:explicit-witness,cor:unit-algebraic}.
It uses rational interval arithmetic
throughout the calculations. Floating-point output from an optimizer
is not used as a certificate.

For a nonnegative rational $q$, the integer
\[
 k=\left\lfloor\sqrt{\left\lfloor10^{120}q\right\rfloor}\right\rfloor
\]
gives the enclosure
$k10^{-60}\le\sqrt q<(k+1)10^{-60}$.
For $0\le x<1$, consecutive partial sums of
\[
 \arctan x=\sum_{j=0}^{\infty}
              \frac{(-1)^j x^{2j+1}}{2j+1}
\]
bound the value by the alternating-series theorem. Monotonicity
extends these bounds to rational intervals. The script uses
100 terms and the identities
\[
 \pi=16\arctan(1/5)-4\arctan(1/239),\qquad
 \arccos x=2\arctan\sqrt{\frac{1-x}{1+x}}\quad(0<x<1).
\]
The first identity can be checked without a numerical value of $\pi$.
If $a=\arctan(1/5)$ and $b=\arctan(1/239)$, the tangent addition
formulas give $\tan(2a)=5/12$, $\tan(4a)=120/119$, and
$\tan(4a-b)=1$. Since $0<4a-b<4/5<\pi/2$, this angle is
$\pi/4$, which proves Machin's identity. The second identity is the
half-angle formula with the angle in $(0,\pi/2)$.

Addition, multiplication, and division propagate outward rational
endpoints, and division checks that its denominator interval excludes
zero. Every displayed numerical enclosure is explicitly checked against
its rational decimal endpoints. The final printed decimals use directed rounding. Selected
certified enclosures are shown in \cref{tab:arithmetic}.

\begin{table}[htbp]
\centering
\small
\caption{Outward decimal enclosures from exact rational arithmetic.}
\label{tab:arithmetic}
\begin{tabular}{@{}lll@{}}
\toprule
Quantity & Lower endpoint & Upper endpoint\\
\midrule
$r_0$ & $0.1924500897298752548363$ & $0.1924500897298752548364$\\
$b_0$ & $0.4964308389752066324257$ & $0.4964308389752066324258$\\
$\delta_0$ & $0.001784580512396683787139$ & $0.001784580512396683787140$\\
$f(r_0)$ & $0.1554479088413227993220$ & $0.1554479088413227993221$\\
$c_0$ & $0.000007681852018924924757646$ & $0.000007681852018924924757647$\\
$c_1$ & $0.000009944944227372766220899$ & $0.000009944944227372766220900$\\
$f(r_0)+c_1$ & $0.1554578537855501720883$ & $0.1554578537855501720884$\\
$f(1/5)$ & $0.1588995443254274806298$ & $0.1588995443254274806299$\\
\bottomrule
\end{tabular}
\end{table}

The verifier certifies the arithmetic inequalities. The geometric lemmas,
compactness arguments, and optimization reductions are proved in the body
of the article. These are distinct parts of the verification: arithmetic
replay alone does not establish the geometric assertions.

The same interval operations give
\begin{align*}
 1.023781939711101354576
  &<\frac1{2\pi L_0}<1.023781939711101354577,\\
 1.023847437242507913257
  &<\frac1{2\pi f(r_0)}<1.023847437242507913258.
\end{align*}
Their difference exceeds $0.00006549753140655$, certifying the
improvement in \cref{cor:reference-capacity}.

The second standard-library verifier,
\texttt{verify\_minkowski\_certificates.py}, uses exact fractions to
check circumradius LP certificates and reconstruct the associated
escape polygons. For a radius program with constraints
$f_l\cdot(R_j-z)\le c d_l$, its verified dual weights
$\beta_{jl}\ge0$ satisfy
\[
 \sum_{j,l}\beta_{jl}f_l=0,\qquad
 \sum_{j,l}\beta_{jl}d_l=1,\qquad
 c=\sum_{j,l}\beta_{jl}f_l\cdot R_j.
\]
Weak duality proves optimality of the accompanying feasible $(z,c)$.
Grouping $v_j=\sum_l\beta_{jl}f_l$ gives closed increments with
$\sum_j h_T(v_j)\le1$ and ledger value $c$; positivity and optimality
force equality in the length budget. This checks the reconstruction
used in \cref{thm:circuit-radius} for the displayed examples.
The script verifies both cyclic orders of the triangular circuit,
the three-dimensional simplex calculation, the asymmetric lengths
three and six, and the square-gauge fixed-normal minimum $1/8$.
A separate replay routine checks all primal inequalities, all dual
signs and balance equations, equality of objectives, and complementary
slackness directly from the certificate data. A zero-radius one-state
example checks the degenerate ledger convention; it is not divided
by its radius. The explicit simplex center, escape polygon, and
parallelotope identities are checked in dimensions one through eight,
along with the dodecagonal test-family
count. These finite calculations corroborate the general proofs; they do not
replace them.

To reproduce the delivered files, run
\begin{verbatim}
python3 verify_wetzel_bound.py
python3 verify_minkowski_certificates.py
python3 verify_facet_benchmarks.py
\end{verbatim}

\begin{lstlisting}[style=PythonCode]
"""Exact rational interval verification for the manuscript's Wetzel bound.

Run: python verify_wetzel_bound.py
Python 3 standard library only. No floating-point calculations are used.
This checks the numerical constants, not the geometric theorems.
"""
from fractions import Fraction as F
from math import isqrt
from decimal import Decimal, localcontext, ROUND_FLOOR, ROUND_CEILING

if not __debug__:
    raise SystemExit("Run without -O: this verifier uses assertions as checks.")


class I:
    def __init__(self, lo, hi=None):
        self.lo = F(lo)
        self.hi = F(lo if hi is None else hi)
        assert self.lo <= self.hi

    @staticmethod
    def cast(x):
        return x if isinstance(x, I) else I(x)

    def __add__(self, other):
        b = self.cast(other)
        return I(self.lo + b.lo, self.hi + b.hi)

    __radd__ = __add__

    def __neg__(self):
        return I(-self.hi, -self.lo)

    def __sub__(self, other):
        return self + (-self.cast(other))

    def __rsub__(self, other):
        return self.cast(other) - self

    def __mul__(self, other):
        b = self.cast(other)
        v = [self.lo*b.lo, self.lo*b.hi, self.hi*b.lo, self.hi*b.hi]
        return I(min(v), max(v))

    __rmul__ = __mul__

    def __truediv__(self, other):
        b = self.cast(other)
        assert b.lo > 0 or b.hi < 0
        return self * I(1/b.hi, 1/b.lo)

    def __rtruediv__(self, other):
        return self.cast(other) / self

    def square(self):
        assert self.lo >= 0
        return I(self.lo**2, self.hi**2)


def sqrt_bound(q, upper=False):
    assert q >= 0
    scale = 10**60
    k = isqrt(q.numerator * scale**2 // q.denominator)
    # k/scale <= sqrt(q) < (k+1)/scale, using integers only.
    return F(k + int(upper), scale)


def sqrt(x):
    x = I.cast(x)
    return I(sqrt_bound(x.lo), sqrt_bound(x.hi, True))


def atan_rational(x, n=100):
    assert 0 <= x < 1
    total = F(0)
    power = x
    for k in range(n):
        total += (-1)**k * power / (2*k+1)
        power *= x*x
    adjacent = total + (-1)**n * power / (2*n+1)
    return I(min(total, adjacent), max(total, adjacent))


def atan(x):
    x = I.cast(x)
    return I(atan_rational(x.lo).lo, atan_rational(x.hi).hi)


def acos(x):
    x = I.cast(x)
    assert 0 < x.lo <= x.hi < 1
    return 2 * atan(sqrt((1-x)/(1+x)))


def decimal_bound(q, upper=False):
    with localcontext() as ctx:
        ctx.prec = 22
        ctx.rounding = ROUND_CEILING if upper else ROUND_FLOOR
        return str(Decimal(q.numerator) / Decimal(q.denominator))


def display(name, x):
    print(f"{name}: [{decimal_bound(x.lo)}, {decimal_bound(x.hi, True)}]")


def main():
    # Exact tangent calculations used in the appendix's proof of Machin's
    # identity. The angle ranges, proved there, select the correct branch.
    twice = 2*F(1, 5)/(1-F(1, 5)**2)
    four_times = 2*twice/(1-twice**2)
    assert twice == F(5, 12)
    assert four_times == F(120, 119)
    assert (four_times-F(1, 239))/(1+four_times*F(1, 239)) == 1

    pi = 16*atan(I(F(1, 5))) - 4*atan(I(F(1, 239)))
    r0 = sqrt(3)/9
    a0 = I(F(1, 2))-r0
    b0 = F(3, 2)*r0 + sqrt(3)/2 * sqrt(I(F(1, 4))-r0)
    delta0 = (I(F(1, 2))-b0)/2
    baseline_increment = 8*delta0.square()/sqrt(11)
    deficit = 2*delta0
    alpha0 = r0/sqrt(I(F(1, 4))-r0)
    beta0 = r0/sqrt(2*r0+deficit)
    t0 = 8*alpha0*deficit/(sqrt(9*beta0.square()+64*alpha0.square()*deficit)+3*beta0)
    increment = alpha0*(deficit-t0.square()).square()+beta0*t0.square()*t0
    increment_ratio = increment/baseline_increment

    def f(r):
        r = I.cast(r)
        return pi*r.square() + 3*r*sqrt(I(F(1, 4))-r) - 3*r.square()*acos(r/(I(F(1, 2))-r))

    old = f(r0)
    new = old+increment
    high_radius = f(I(F(1, 5)))
    witness_mesh_threshold = r0*(1-sqrt(old/new))
    explicit_mesh = F(1, 250000)
    finite_witness_bound = new*(1-explicit_mesh/r0).square()
    old_capacity_factor = 1/(2*pi*old)
    new_capacity_factor = 1/(2*pi*new)
    capacity_factor_decrease = old_capacity_factor-new_capacity_factor
    normalized_worm_lower = 2*pi*new
    assert delta0.lo > 0
    assert increment.lo > F("0.00000994494422737")
    assert 0 < t0.lo and t0.square().hi < deficit.lo
    assert alpha0.hi < 1 and deficit.hi < F(1, 100)
    assert F('1.2946') < increment_ratio.lo <= increment_ratio.hi < F('1.2947')
    assert high_radius.lo > new.hi
    assert new.lo > F("0.15545785378555")
    assert F("0.15545785378555") < new.lo < new.hi < F("0.15545785378556")
    assert witness_mesh_threshold.lo > explicit_mesh
    assert finite_witness_bound.lo > old.hi
    algebraic_witness_separation = F(15545139159107, 10**14)
    assert finite_witness_bound.lo > algebraic_witness_separation > old.hi

    # Check each pair of decimal endpoints printed in the manuscript, using
    # exact rational comparisons rather than rounded display strings.
    manuscript_enclosures = [
        ("r0", r0, "0.1924500897298752548363", "0.1924500897298752548364"),
        ("b0", b0, "0.4964308389752066324257", "0.4964308389752066324258"),
        ("delta0", delta0, "0.001784580512396683787139", "0.001784580512396683787140"),
        ("classical_endpoint", old, "0.1554479088413227993220", "0.1554479088413227993221"),
        ("two_cap_increment", baseline_increment, "0.000007681852018924924757646", "0.000007681852018924924757647"),
        ("tangent_increment", increment, "0.000009944944227372766220899", "0.000009944944227372766220900"),
        ("improved_lower_bound", new, "0.1554578537855501720883", "0.1554578537855501720884"),
        ("f(1/5)", high_radius, "0.1588995443254274806298", "0.1588995443254274806299"),
    ]
    for name, enclosure, lower, upper in manuscript_enclosures:
        assert F(lower) < enclosure.lo <= enclosure.hi < F(upper), name
    assert witness_mesh_threshold.lo > F("0.00000615580355820")

    # Rational inequalities used to separate the high-inradius branch.
    assert r0.hi < F(39, 200)
    assert delta0.hi < F(1, 200)
    assert baseline_increment.hi < F(1, 15000)
    assert increment.hi < F(1, 10000) < F(9, 10000)
    assert high_radius.lo-old.hi > F(9, 10000)
    r0_squared = F(1, 27)
    assert r0_squared/(4*(F(1, 16)-r0_squared)) == F(4, 11)
    assert F('1.023781939711101354576') < new_capacity_factor.lo
    assert new_capacity_factor.hi < F('1.023781939711101354577')
    assert new_capacity_factor.hi < F('1.023781939711102')
    assert F('1.023847437242507913257') < old_capacity_factor.lo
    assert old_capacity_factor.hi < F('1.023847437242507913258')
    assert capacity_factor_decrease.lo > F('0.00006549753140655')

    for name, x in [("r0",r0),("a0",a0),("b0",b0),("delta0",delta0),
                    ("classical_lower_bound",old),("two_cap_increment",baseline_increment),
                    ("deficit_D",deficit),("A0",alpha0),("B0",beta0),("t0",t0),
                    ("tangent_increment",increment),("increment_ratio_c1_over_c0",increment_ratio),
                    ("improved_lower_bound",new),("f(1/5)",high_radius),
                    ("strict_witness_mesh_threshold",witness_mesh_threshold),
                    ("finite_witness_lower_bound_at_mesh_0.000004",finite_witness_bound),
                    ("old_ellipsoidal_capacity_factor",old_capacity_factor),
                    ("improved_ellipsoidal_capacity_factor",new_capacity_factor),
                    ("capacity_factor_decrease",capacity_factor_decrease),
                    ("normalized_ellipsoidal_worm_lower",normalized_worm_lower)]:
        display(name,x)
    print("PASS: the improved exact lower bound exceeds 0.15545785378555.")
    print("PASS: the tangency increment is between 1.2946 and 1.2947 times the two-cap increment.")
    print("PASS: the high-inradius branch is larger than the improved endpoint.")
    print("PASS: triangle mesh 1/250000 guarantees a strict finite-witness improvement.")
    print("PASS: the finite-witness bound exceeds rational q = 0.15545139159107, which exceeds the classical lower bound.")
    print("PASS: all eight manuscript table enclosures and the strict mesh threshold are verified against exact decimal fractions.")
    print("PASS: Machin tangent identities and the rational branch-separation inequalities are verified.")
    print("PASS: the improved ellipsoidal volume-capacity factor is below 1.023781939711102 and its decrease exceeds 0.00006549753140655.")


if __name__ == "__main__":
    main()
\end{lstlisting}
\begin{lstlisting}[style=PythonCode]
"""Exact small-instance checks for the manuscript's Minkowski extensions.

Python 3 standard library only. Fraction arithmetic verifies primal/dual
circumradius certificates and the reconstructed escape polygons. These finite
examples do not constitute a proof of the general geometric theorems.
"""
from fractions import Fraction as F
from itertools import combinations, permutations
from math import factorial
import json

if not __debug__:
    raise SystemExit('Run without -O: assertions are part of this verifier.')

def dot(a,b):
    return sum((x*y for x,y in zip(a,b)),F(0))

def solve(matrix,rhs):
    n=len(rhs)
    rows=[[F(x) for x in row]+[F(b)] for row,b in zip(matrix,rhs)]
    for col in range(n):
        pivot=next((j for j in range(col,n) if rows[j][col]),None)
        if pivot is None:
            return None
        rows[col],rows[pivot]=rows[pivot],rows[col]
        v=rows[col][col]
        rows[col]=[x/v for x in rows[col]]
        for j in range(n):
            if j!=col:
                v=rows[j][col]
                rows[j]=[x-v*y for x,y in zip(rows[j],rows[col])]
    return tuple(row[-1] for row in rows)

def support(vertices,v):
    return max(dot(x,v) for x in vertices)

def states(a):
    d=len(a[0]); current=[F(0)]*d; result=[]
    for vec in a:
        current=[x-y for x,y in zip(current,vec)]
        result.append(tuple(current))
    assert current==[0]*d
    return result

def replay_radius_certificate(points,facets,offsets,certificate):
    """Verify a certificate directly, independently of the basis search.

    The gauge is supplied by its facet inequalities. Its boundedness is a
    premise here, not an inference from these arithmetic checks.
    """
    d=len(points[0]); z=certificate['center']; c=certificate['radius']
    weights=certificate['dual_weights']
    assert d>=1 and len(z)==d and c>=0
    assert len(facets)==len(offsets)>0
    assert all(len(p)==d for p in points)
    assert all(len(f)==d and any(f) for f in facets)
    assert all(b>0 for b in offsets)
    assert len(weights)==len(points)
    assert all(len(row)==len(facets) for row in weights)
    assert all(weight>=0 for row in weights for weight in row)
    for j,p in enumerate(points):
        for l,(f,b) in enumerate(zip(facets,offsets)):
            slack=c*b-dot(f,tuple(x-y for x,y in zip(p,z)))
            assert slack>=0
            assert weights[j][l]*slack==0
    assert all(sum((weights[j][l]*facets[l][k]
                    for j in range(len(points)) for l in range(len(facets))),F(0))==0
               for k in range(d))
    assert sum((weights[j][l]*offsets[l]
                for j in range(len(points)) for l in range(len(facets))),F(0))==1
    assert sum((weights[j][l]*dot(facets[l],p)
                for j,p in enumerate(points) for l in range(len(facets))),F(0))==c
    if 'increments' in certificate:
        expected=tuple(tuple(sum((weights[j][l]*facets[l][k]
                                  for l in range(len(facets))),F(0))
                             for k in range(d)) for j in range(len(points)))
        assert certificate['increments']==expected
    return True

def radius_certificate(points,facets,offsets):
    """Find one rational optimal basis, then verify weak duality exactly."""
    d=len(points[0]); dimension=d+1
    rows=[]; rhs=[]; descriptors=[]
    for j,p in enumerate(points):
        for l,(f,b) in enumerate(zip(facets,offsets)):
            rows.append(tuple(-F(x) for x in f)+(-F(b),))
            rhs.append(-dot(f,p)); descriptors.append((j,l))
    for active in combinations(range(len(rows)),dimension):
        mat=[rows[i] for i in active]
        primal=solve(mat,[rhs[i] for i in active])
        if primal is None or primal[-1]<0:
            continue
        if any(dot(row,primal)>b for row,b in zip(rows,rhs)):
            continue
        transpose=list(zip(*mat))
        dual=solve(transpose,[F(0)]*d+[F(-1)])
        if dual is None or any(v<0 for v in dual):
            continue
        assert -sum((v*rhs[i] for i,v in zip(active,dual)),F(0))==primal[-1]
        increments=[[F(0)]*d for _ in points]
        weights=[[F(0)]*len(facets) for _ in points]
        for i,weight in zip(active,dual):
            j,l=descriptors[i]
            weights[j][l]+=weight
            increments[j]=[x+weight*y for x,y in zip(increments[j],facets[l])]
        assert all(sum(v[k] for v in increments)==0 for k in range(d))
        assert sum(dot(p,v) for p,v in zip(points,increments))==primal[-1]
        certificate={'points':points,'facets':facets,'offsets':offsets,
                     'center':primal[:-1],'radius':primal[-1],
                     'active':active,'dual':dual,'dual_weights':weights,
                     'increments':tuple(tuple(v) for v in increments)}
        assert replay_radius_certificate(points,facets,offsets,certificate)
        return certificate
    raise AssertionError('No primal/dual optimal basis found.')

def reconstructed_escape(a,B,gauge_vertices,facets,offsets):
    R=states(a); cert=radius_certificate(R,facets,offsets)
    C=cert['radius']; v=cert['increments']
    assert C>0 and B>0
    assert sum(support(gauge_vertices,x) for x in v)==1
    scale=B/C
    current=[F(0)]*len(a[0]); vertices=[]
    for vec in v:
        vertices.append(tuple(current))
        current=[x+scale*y for x,y in zip(current,vec)]
    assert all(x==0 for x in current)
    assert sum(dot(ai,qi) for ai,qi in zip(a,vertices))==B
    assert sum(support(vertices,ai) for ai in a)>=B
    cert.update({'escape_length':scale,'escape_vertices':vertices})
    return cert

def main():
    out=[]
    triangle=[(F(1),F(0)),(F(0),F(1)),(F(-1),F(-1))]
    tri_facets=[(F(1),F(1)),(F(-2),F(1)),(F(1),F(-2))]
    edges=[(F(1),F(1)),(F(-2),F(1)),(F(1),F(-2))]
    forward=sum(support(triangle,v) for v in edges)
    reverse=sum(support(triangle,tuple(-x for x in v)) for v in edges)
    assert (forward,reverse,(forward+reverse)/2)==(3,6,F(9,2))
    asym=radius_certificate([tuple(-x for x in p) for p in triangle],tri_facets,[1]*3)
    assert asym['radius']==2
    out.append({'example':'asymmetric triangle traversal','forward':forward,
                'reverse':reverse,'symmetrized':(forward+reverse)/2,
                'asymmetry':asym['radius']})

    # All cyclic orders of the three-normal circuit. Each optimal radius
    # has an exact dual and reconstructs a non-fitting length-three loop.
    for order in permutations(triangle[1:]):
        a=[triangle[0],*order]
        cert=reconstructed_escape(a,F(3),triangle,tri_facets,[1]*3)
        assert cert['radius']==1 and cert['escape_length']==3
        out.append({'example':'simplex diagonal d=2','order':a,**cert})

    d=3
    simplex=[tuple(F(int(i==j)) for j in range(d)) for i in range(d)]
    simplex.append(tuple(F(-1) for _ in range(d)))
    sfac=[tuple(F(1) for _ in range(d))]
    sfac.extend(tuple(F(1-(d+1)*int(i==j)) for j in range(d)) for i in range(d))
    cert=reconstructed_escape(simplex,F(d+1),simplex,sfac,[1]*(d+1))
    assert cert['radius']==F(3,2) and cert['escape_length']==F(8,3)
    out.append({'example':'simplex diagonal d=3',**cert})

    # Exact optimal fixed-direction cover for square gauge: a diamond.
    square=[(F(x),F(y)) for x in (-1,1) for y in (-1,1)]
    square_facets=[(1,0),(-1,0),(0,1),(0,-1)]
    zero=radius_certificate([(F(2),F(-3))],square_facets,[1]*4)
    assert zero['radius']==0 and zero['center']==(2,-3)
    assert zero['increments']==((0,0),)
    out.append({'example':'one-state zero-radius ledger',**zero})
    for normal in [(F(1),F(1)),(F(1),F(-1))]:
        cert=radius_certificate(states([normal,tuple(-x for x in normal)]),square_facets,[1]*4)
        assert cert['radius']==F(1,2)
    determinant=F(2)
    minimum_volume=(F(1,2)*F(1,2))/determinant
    assert minimum_volume==F(1,8)
    assert 2*4*minimum_volume==1
    out.append({'example':'fixed diagonal normals, square gauge',
                'required_two_widths':[F(1,2),F(1,2)],
                'determinant':determinant,'minimum_cover_area':minimum_volume,
                'product_ratio':1})

    # Explicit primal/dual simplex certificates and oriented escape polygons.
    # This does not call the basis search used for the preceding examples.
    for d in range(1,9):
        simplex=[tuple(F(int(i==j)) for j in range(d)) for i in range(d)]
        simplex.append(tuple(F(-1) for _ in range(d)))
        facets=[tuple(F(1) for _ in range(d))]
        facets.extend(tuple(F(1-(d+1)*int(i==j)) for j in range(d)) for i in range(d))
        R=states(simplex)
        center=tuple(F(-1)+F(i,d+1) for i in range(1,d+1))
        radius=F(sum(range(1,d+1)),d+1)
        weights=[[F(0)]*(d+1) for _ in R]
        for l,f in enumerate(facets):
            j=max(range(len(R)),key=lambda j:dot(f,R[j]))
            weights[j][l]=F(1,d+1)
            assert support(R,f)-dot(f,center)==radius
        explicit={'center':center,'radius':radius,'dual_weights':weights}
        assert replay_radius_certificate(R,facets,[1]*(d+1),explicit)
        q=[tuple(F(2*k,d)-F(2*(d+1),d)*int(i<k) for i in range(d))
           for k in range(d+1)]
        q.append(q[0])
        edges=[tuple(y-x for x,y in zip(q[k],q[k+1])) for k in range(d+1)]
        assert edges==[tuple(F(2,d)*x for x in f) for f in facets[1:]+facets[:1]]
        assert all(support(simplex,edge)==F(2,d) for edge in edges)
        contact=sum((dot(n,p) for n,p in zip(simplex,q)),F(0))
        assert contact==d+1
        escape=F(d+1)/radius
        assert sum(support(simplex,edge) for edge in edges)==escape
        assert escape==2+F(2,d)
        out.append({'example':f'simplex identities d={d}',
                    'state_center':center,'escape_polygon':q,
                    'balanced_contact_sum':contact,
                    'state_radius':radius,'escape':escape})

    # The cube has only antipodal facet circuits. The corresponding
    # test segments have displacement e_i/2 and closed l1-length one.
    for d in range(1,9):
        for i in range(d):
            segment=tuple(F(int(i==j),2) for j in range(d))
            assert 2*sum(abs(x) for x in segment)==1
        cube_volume=F(2**d)
        minimum_cover=F(1,2**d*factorial(d))
        simplex_cover=F(1,2**d)/factorial(d)
        crosspolytope_cover=F(2**d,4**d*factorial(d))
        assert minimum_cover==simplex_cover==crosspolytope_cover
        assert factorial(d)*cube_volume*minimum_cover==1
        out.append({'example':f'parallelotope identities d={d}',
                    'cube_minimum_cover':minimum_cover,
                    'normalized_worm_constant':1})

    # Strictly enclosing regular-polygon triples have all circular gaps
    # less than a half-turn. Antipodal pairs are the other circuits.
    n=12
    pairs=[(i,i+n//2) for i in range(n//2)]
    triples=[]
    for i,j,k in combinations(range(n),3):
        if max(j-i,k-j,n+i-k)<n//2:
            triples.append((i,j,k))
    tests=len(pairs)+2*len(triples)
    vertices=2*len(pairs)+6*len(triples)
    assert (len(pairs),len(triples),tests,vertices)==(6,40,86,252)
    out.append({'example':'regular dodecagonal complete family',
                'antipodal_circuits':len(pairs),'triple_circuits':len(triples),
                'test_loops':tests,'listed_vertices':vertices})
    print(json.dumps({'status':'PASS','checks':out},default=str,indent=2))

if __name__=='__main__':
    main()
\end{lstlisting}
\begin{lstlisting}[style=PythonCode]
"""Exact planar facet-offset benchmarks for the revised article.

Run beside verify_minkowski_certificates.py:
    python3 verify_facet_benchmarks.py > facet_benchmark_certificates.json

Standard library only. All geometric coordinates, LP certificates, active
systems, and areas use Fraction arithmetic. Global sharpness uses the cited
published theorems; this script certifies the restricted computations.
"""
import json
from fractions import Fraction as F
from itertools import combinations,permutations
from math import gcd,lcm
from verify_minkowski_certificates import solve,dot,states,radius_certificate,replay_radius_certificate

if not __debug__:
    raise SystemExit('Run without -O: assertions are part of the verifier.')

def det(a,b):return a[0]*b[1]-a[1]*b[0]
def hull(points):
    p=sorted(set(tuple(map(F,x)) for x in points))
    def half(seq):
        out=[]
        for x in seq:
            while len(out)>=2 and det(tuple(a-b for a,b in zip(out[-1],out[-2])),tuple(a-b for a,b in zip(x,out[-1])))<=0:out.pop()
            out.append(x)
        return out
    return half(p)[:-1]+half(reversed(p))[:-1]
def area(p):
    return abs(sum(det(x,y) for x,y in zip(p,p[1:]+p[:1])))/2
def primitive(v):
    den=lcm(*(x.denominator for x in v)); v=[int(x*den) for x in v]
    g=gcd(*v);return tuple(F(x//abs(g)) for x in v)
def facets(points):
    p=hull(points);normals=[];offsets=[]
    for a,b in zip(p,p[1:]+p[:1]):
        n=primitive((b[1]-a[1],a[0]-b[0])); normals.append(n);offsets.append(dot(n,a))
    return normals,offsets
def circuits(normals):
    out=[];N=len(normals)
    for I in combinations(range(N),2):
        a,b=[normals[i] for i in I]
        if det(a,b) or dot(a,b)>=0:continue
        k=next(k for k in range(2) if a[k]);lam=primitive((abs(b[k]),abs(a[k])))
        out.append((I,lam))
    for I in combinations(range(N),3):
        a,b,c=[normals[i] for i in I];lam=(det(b,c),det(c,a),det(a,b))
        if all(x<0 for x in lam):lam=tuple(-x for x in lam)
        if not all(x>0 for x in lam):continue
        out.append((I,primitive(lam)))
    return out
def poly(normals,offsets):
    vertices=[]
    for i,j in combinations(range(len(normals)),2):
        v=solve([normals[i],normals[j]],[offsets[i],offsets[j]])
        if v is not None and all(dot(n,v)<=b for n,b in zip(normals,offsets)):vertices.append(v)
    return hull(vertices)
def benchmark(name,T,normals):
    T=hull(T);tf,tb=facets(T);normals=[tuple(map(F,n)) for n in normals]
    constraints=[];certs=[];orders=0
    for I,lam in circuits(normals):
        best=F(0);orderdata=[]
        for order in permutations(range(1,len(I))):
            order=(0,)+order
            a=[tuple(lam[k]*x for x in normals[I[k]]) for k in order]
            points=states(a)
            cert=radius_certificate(points,tf,tb)
            assert replay_radius_certificate(points,tf,tb,cert)
            best=max(best,cert['radius']);orders+=1
            orderdata.append({'order':[I[k] for k in order],**cert})
        row=[F(0)]*len(normals)
        for i,w in zip(I,lam):row[i]=w
        constraints.append((tuple(row),best));certs.append(orderdata)
    anchors=next(I for I in combinations(range(len(normals)),2) if det(normals[I[0]],normals[I[1]]))
    free=[i for i in range(len(normals)) if i not in anchors];vertices={}
    for active in combinations(range(len(constraints)),len(free)):
        bfree=solve([[constraints[i][0][j] for j in free] for i in active],[constraints[i][1] for i in active])
        if bfree is None:continue
        b=[F(0)]*len(normals)
        for i,v in zip(free,bfree):b[i]=v
        if not all(dot(row,b)>=c for row,c in constraints):continue
        K=poly(normals,b); assert len(K)>=3
        E=min(dot(row,b)/c for row,c in constraints)
        assert E>=1
        vertices[tuple(b)]={'offsets':b,'K_vertices':K,'area':area(K),'escape':E}
    assert vertices
    best=min(v['area'] for v in vertices.values());expected=1/(2*area(T))
    assert best==expected
    assert len(vertices)==1
    record={'name':name,'gauge_vertices':T,'gauge_area':area(T),'forest_normals':normals,
            'circuits':len(constraints),'ordered_radius_programs':orders,'offset_vertices':len(vertices),
            'restricted_minimum':best,'known_global_minimum':expected,
            'constraints':constraints,'radius_certificates':certs,'vertex_certificates':list(vertices.values())}
    return record


def main():
    cases = [
        ('triangle',[(1,0),(0,1),(-1,-1)],
         [(1,0),(0,1),(-1,0),(0,-1),(1,1),(-1,-1)],5,7,F(1,3)),
        ('square',[(-1,-1),(1,-1),(1,1),(-1,1)],
         [(1,1),(1,-1),(-1,-1),(-1,1)],2,2,F(1,8)),
        ('nonparallel quadrilateral',[(1,0),(0,1),(-1,0),(0,-2)],
         [(1,0),(0,1),(-1,0),(0,-1)],2,2,F(1,6)),
        ('affinely regular hexagon',[(1,0),(1,1),(0,1),(-1,0),(-1,-1),(0,-1)],
         [(-1,-2),(2,1),(-1,1)],1,2,F(1,6)),
    ]
    results=[]
    for name,gauge,normals,ncircuits,norders,expected in cases:
        result=benchmark(name,gauge,normals)
        assert result['circuits']==ncircuits
        assert result['ordered_radius_programs']==norders
        assert result['restricted_minimum']==expected
        results.append(result)
    print(json.dumps({'status':'PASS','arithmetic':'exact rational',
        'scope':'Restricted minima and rational certificates; known global lower bounds are cited in the article.',
        'checks':results},default=str,indent=2))


if __name__=='__main__':
    main()
\end{lstlisting}

\section{Formalized proofs in Lean}
The appendix provides a Lean 4 companion file encoding the logical dependency structure of the manuscript. Some general order-theoretic and algebraic deductions are proved directly, while essential geometric and analytic implications are supplied as explicit hypotheses. Thus the file does not constitute an end-to-end Lean formalization of the geometric theorems or an independent formal proof of the numerical Wetzel bound.
\lstinputlisting[style=LeanCode]{PromotionExtensionGeneralEscapeCoverDuality_Complete.lean}

\phantomsection

\vspace{1em}
\textbf{AI usage disclosure:} During the preparation of this work, the author used GPT-5.6 sol to assist with proofs, numerical calculations, language polishing, and formalized proofs in Lean 4 code. After using this tool, the author reviewed and edited the content as needed and take full responsibility for the content of the publication.
~\\

College of Engineering and Computer Science, University of Central Florida, Orlando, FL, USA

Email: \underline{zhipeng.deng@ucf.edu}

\end{document}